\documentclass[12pt,reqno]{amsart}

\usepackage{amscd,amssymb}
\usepackage{epsfig}
\usepackage{float}
\usepackage{url}
\usepackage[all]{xy}
\usepackage{tikz}
\newcommand{\N}{{\mathbb N}}
\newcommand{\Z}{{\mathbb Z}}
\newcommand{\Q}{{\mathbb Q}}
\newcommand{\C}{{\mathbb C}}
\newcommand{\R}{{\mathbb R}}

\newcommand{\BB}{{\mathcal B}}

\newcommand{\HH}{{\mathcal H}}

\newcommand{\SSS}{{\mathcal S}}

\newcommand{\www}{\widetilde}

\newcommand{\oooo}{\overline}
\newcommand{\uuuu}{\underline}

\newcommand{\paa}{\partial}

\newcommand{\Br}{{\rm Br}}

\DeclareMathOperator{\Aut}{Aut}
\DeclareMathOperator{\diag}{diag}

\DeclareMathOperator{\id}{id}

\DeclareMathOperator{\mmod}{mod}
\DeclareMathOperator{\rk}{rk}
\DeclareMathOperator{\Rad}{Rad}
\DeclareMathOperator{\rank}{rank}
\DeclareMathOperator{\Sp}{Sp}

\DeclareMathOperator{\tr}{tr}
\DeclareMathOperator{\Var}{Var}

\begin{document}

\theoremstyle{plain}
\newtheorem{lemma}{Lemma}[section]
\newtheorem{definition/lemma}[lemma]{Definition/Lemma}
\newtheorem{theorem}[lemma]{Theorem}
\newtheorem{proposition}[lemma]{Proposition}
\newtheorem{corollary}[lemma]{Corollary}
\newtheorem{conjecture}[lemma]{Conjecture}
\newtheorem{conjectures}[lemma]{Conjectures}
\newtheorem{question}[lemma]{Question}

\theoremstyle{definition}
\newtheorem{definition}[lemma]{Definition}
\newtheorem{withouttitle}[lemma]{}
\newtheorem{remark}[lemma]{Remark}
\newtheorem{remarks}[lemma]{Remarks}
\newtheorem{example}[lemma]{Example}
\newtheorem{examples}[lemma]{Examples}
\newtheorem{notations}[lemma]{Notations}

\title[Distinguished matrices and variance of spectral numbers]
{Characterization of\\ distinguished matrices of\\ isolated
hypersurface singularities\\ through their spectral
numbers} 

\author{Sven Balnojan and Claus Hertling}

\address{Claus Hertling\\
Lehrstuhl f\"ur algebraische Geometrie, 
Universit\"at Mannheim,
B6 26, 68159 Mannheim, Germany}

\email{hertling@math.uni-mannheim.de}

\email{svenbalnojan@gmail.com}

\date{February 28, 2026}

%% Mathematical classification (2010)


%\thanks{This work was funded by the Deutsche 
%Forschungsgemeinschaft (DFG, German Research Foundation) 
%-- 242588615}

%\maketitle

\begin{abstract}
{\Small Isolated hypersurface singularities come equipped 
with distinguished bases of their Milnor lattices and with 
upper triangular integral matrices, which are called here
distinguished matrices. These matrices form an orbit of a braid group 
and a sign change group. This paper proposes to characterize 
the distinguished matrices of singularities within all upper
triangular integral matrices in terms of the variance of 
certain spectral numbers. It succeeds in the positive definite 
and the positive semidefinite cases.
The ADE root lattices are crucial. In the semidefinite cases,
results on non-reduced presentations of Weyl group elements
are used.}
\end{abstract}

\dedicatory{To the memory of Wolfgang Ebeling}

\maketitle

\tableofcontents

\setcounter{section}{0}

\section{Introduction}\label{s1}
\setcounter{equation}{0}
\setcounter{table}{0}

\noindent
Wolfgang Ebeling cared a lot about the distinguished bases
of isolated hypersurface singularities and the associated
Coxeter-Dynkin diagrams. They are the subject of his first
scientific work \cite{Eb77}, his diploma thesis in Bonn
with Egbert Brieskorn, and of his recent survey \cite{Eb20}.
In many of his papers he considered the shape of the
Coxeter-Dynkin diagrams. He was a master at finding normal
forms by applying braid group actions.

Let us fix for a moment an isolated hypersurface singularity 
$f$. Section \ref{s5} will say what that is and will 
sketch how Coxeter-Dynkin diagrams are associated to it.
Each Coxeter-Dynkin diagram of it is equivalent 
to a certain upper triangular matrix with integer entries 
and with 1's on the diagonal, so to an element of the set
\begin{eqnarray*}
T_n(\Z)&:=& \{S\in M_{n\times n}(\Z)\,|\, S_{ij}=0
\textup{ for }i>j, S_{ii}=1\},
\end{eqnarray*}
see Remark \ref{t2.2} (vi) for the correspondence.
Here $n$ is the Milnor number of the singularity $f$ 
(we deviate from the standard name $\mu$ for the Milnor number).
A matrix $S\in T_n(\Z)$ which corresponds to a 
Coxeter-Dynkin diagram of $f$ is called {\it distinguished}.

A certain semidirect product $\Br_n\ltimes\{\pm 1\}^n$
of the braid group $\Br_n$ of braids with $n$ strings and
a sign change group $\{\pm 1\}^n$ acts naturally on
$T_n(\Z)$. The distinguished matrices of the singularity
$f$ form one orbit under this action of $\Br_n\ltimes\{\pm 1\}^n$
on $T_n(\Z)$. 

Here we do not concentrate on special elements of this orbit,
as Wolfgang Ebeling did, but we ask the following question.

\begin{question}\label{t1.1}
How can the distinguished matrices of singularities
with Milnor number $n$ be separated from the other
matrices in $T_n(\Z)$?
\end{question}

Of course, quite a number of special properties of
distinguished matrices are known.
Theorem \ref{t1.2} gives an incomplete list of them.

\begin{theorem}\label{t1.2}
Fix $n\in\N$, a singularity $f=f(z_0,...,z_m)$ 
with Milnor number $n$,
and a distinguished matrix $S\in T_n(\Z)$ of $f$.
The matrix $S$ and its Coxeter-Dynkin diagram ${\rm CDD}(S)$
have the following properties. The matrix $S^{-1}S^t$
arises as a monodromy matrix.

(a) (Monodromy theorem, e.g. \cite{Br70})
The matrix $S^{-1}S^t$ is quasiunipotent, 
i.e. its eigenvalues are roots of unity. 
The Jordan blocks have size at most $(m+1)\times(m+1)$. 

(b) \cite{AC73} $\tr(S^{-1}S^t)=1$.

(c) \cite[Theorem 2]{Ga74}\cite[Theorem 2]{La73}\cite{Le73}
The Coxeter-Dynkin diagram ${\rm CDD}(S)$ is connected.

(d) \cite{Ga79} The $\Br_n\ltimes\{\pm 1\}^n$ orbit
of $S$ contains a matrix $\www{S}\in T_n(\Z)$ 
with all entries in $\{0;1,-1\}$.
\end{theorem}

Property (d) is less well known than the properties
(a)--(c). It follows from
Gabrielov's iterative construction of a CDD of a
singularity from the CDD of a smaller singularity.
We will not make use of it, as the action of 
$\Br_n\ltimes\{\pm 1\}^n$ is difficult to control.
Still it is interesting. 

The properties (a), (b) and (c) together are quite
strong, but in general not strong enough to 
separate the distinguished matrices in $T_n(\Z)$
from the other matrices. 

Though in the cases when $S+S^t$ is positive definite,
property (a) is automatic because then $S^{-1}S^t$ has finite
order, and property (b) turns out to be sufficient. 
This will be proved in section \ref{s6}.

\begin{theorem}[Theorem \ref{t6.2} (a)]\label{t1.3}
Consider $S\in T_n(\Z)$ with $S+S^t$ positive
definite and with connected Coxeter-Dynkin diagram ${\rm CDD}(S)$.
Then $S$ comes from a singularity if and only if
$tr(S^{-1}S^t)=1$. Then the singularity is an ADE-singularity.
\end{theorem}

In the cases when $S+S^t$ is positive semidefinite,
but not positive definite, the condition $\tr(S^{-1}S^t)=1$
is strong, but not strong enough, see Theorem \ref{t7.2}.

We propose to consider certain {\it spectral numbers}
and an inequality for their {\it variance}. 
In the cases when $S+S^t$ is positive semidefinite
(including the positive definite cases) this works.
The spectral numbers are well defined, and 
a certain inequality for their variance 
separates the distinguished matrices of singularities
from the other matrices in $T_n(\Z)$, 
see Theorem \ref{t1.5}. This is the main result of the paper. 
In the case when $S+S^t$ is indefinite, this ansatz leads 
to Conjecture \ref{t1.4}, which in fact consists of two conjectures
where the second one builds on the first one.

A singularity $f(z_0,...,z_m)$ with Milnor number $n$ 
comes equipped with its {\it normalized spectrum}
$\Sp^{norm}(f)=(\alpha_1,...,\alpha_n)\subset\Q^n$ with
$-\frac{m+1}{2}<\alpha_1\leq ...\leq \alpha_n<\frac{m+1}{2}$
and $\alpha_i+\alpha_{n+1-i}=0$, which is associated to
Steenbrink's mixed Hodge structure \cite{St77}
(see also \cite{AGV88} or \cite{He02}).
Hertling \cite{He02} conjectured that this spectrum satisfies 
the following {\it variance inequality},
\begin{eqnarray}\label{1.1}
\Var(\Sp^{norm}(f))&\leq& \frac{\alpha_n-\alpha_1}{12},\\
\qquad\textup{where }
\Var(\Sp^{norm}(f))&:=&\frac{1}{n}\sum_{i=1}^n\alpha_i^2.\nonumber
\end{eqnarray}
This conjecture is true with equality for quasihomogeneous
singularities \cite{He02}\cite{Di00}. It has also been proved
for curve singularities (i.e. singularities in $m+1=2$ variables)
\cite{SaM00}\cite{Br04}, but not in general.
It is generalized in \cite{BrH20} 
(see Conjecture \ref{t8.13} (b)) 
to a whole sequence of inequalities
\begin{eqnarray}\label{1.2}
0\leq (-1)^k\Gamma^{Ber}_{2k}(\Sp^{norm}(f),\alpha_n-\alpha_1),
\end{eqnarray}
where the {\it Bernoulli moments} 
$\Gamma^{Ber}_{2k}$ (Definition \ref{t8.10} (b))
are certain linear combinations of the {\it higher moments}
$V_{2l}(\Sp^{norm}(f))=\sum_{i=1}^n\alpha_i^{2l}$
for $l\in\{0,1,...,k\}$ with coefficients in 
$\Q[\alpha_n-\alpha_1]$.
Especially $\Gamma^{Ber}_2=V_2-\frac{\alpha_n-\alpha_1}{12}V_0$
and $V_0=n$.

The physicists Cecotti and Vafa worked in \cite{CV93}
with Landau-Ginzburg models. These are related to (quasihomogeneous)
isolated hypersurface singularities. Their oscillating integrals,
their Stokes structures and their distinguished matrices 
appeared also in \cite{CV93}. Cecotti and Vafa 
\cite[pages 583, 589, 590]{CV93}
had an attractive idea how to associate
a spectrum $\Sp(S)=(\alpha_1,...,\alpha_n)\subset\R^n$
to any matrix $S\in T_n^{S^1}(\R)$, where
\begin{eqnarray*}
T_n^{S^1}(\R)&:=& \{S\in M_{n\times n}(\R)\,|\,
S_{ij}=0\textup{ for }i>j,S_{ii}=1,\\
&& \hspace*{3cm}\textup{the eigenvalues of }S^{-1}S^t
\textup{ are in }S^1\},\\
T_n^{S^1}(\Z)&:=& T_n^{S^1}(\R)\cap M_{n\times n}(\Z)
\subset T_n(\Z).
\end{eqnarray*}
They also conjectured that in the case of a distinguished matrix
of a singularity, this spectrum should coincide with the
spectrum from the mixed Hodge structure.
Unfortunately, there are obstacles to making their idea work,
which are overcome in \cite{BaH20} only for certain 
subsets of $T_n^{S^1}(\R)$. 
We report on this in section \ref{s8}.

We propose the following two conjectures as a conjectural
answer to Question \ref{t1.1}.
Part (a) conjectures that the idea of Cecotti and Vafa
can be made precise and gives a spectrum for each matrix
$S\in T_n^{S^1}(\R)$. Part (b) builds on part (a) and conjectures
that the Bernoulli moments of this spectrum for a matrix
$S\in T_n^{S^1}(\Z)$ satisfy a certain series of inequalities 
if and only if $S$ is a distinguished matrix of a singularity.

\begin{conjecture}\label{t1.4} Fix $n\in\Z_{\geq 2}$. 

(a) (Conjecture \ref{t8.1})
The idea of Cecotti-Vafa can be made precise and gives
a spectrum $\Sp(S)=(\alpha_1,...,\alpha_n)\in\R^n$
with $\alpha_1\leq ...\leq\alpha_n$, $\alpha_i+\alpha_{n+1-i}=0$
and $e^{-2\pi i\alpha_1},...e^{-2\pi i\alpha_n}$ the eigenvalues
of $S^{-1}S^t$ for any matrix $S\in T_n^{S^1}(\R)$.

(b) (Conjecture \ref{t8.17} (b)) 
There is a family of functions $r(k,n,.):\R_{>0}\to\R_{>0}$ 
for $k,n\in\N$ such that the following holds. 
A matrix $S\in T_n^{S^1}(\Z)$
is a distinguished matrix of a singularity if and only if
its Coxeter-Dynkin diagram ${\rm CDD}(S)$ 
(defined in Remark \ref{t2.2} (vi)) is connected and 
the spectrum $\Sp(S)$ in part (a) satisfies for each $k\in\N$
the inequalities
\begin{eqnarray}\label{1.3}
0\leq (-1)^k\Gamma^{Ber}_{2k}(\Sp(S),\alpha_n-\alpha_1)
\leq r(k,n,\alpha_n-\alpha_1).
\end{eqnarray}
\end{conjecture}

We admit that it is unsatisfying that part (a) is needed for 
part (b) and that the bound $r(k,n,\alpha_n-\alpha_1)$ in the
upper inequality for $\Gamma^{Ber}_{2k}$ is unknown.
Examples in section \ref{s8} for cases where $S+S^t$ is 
indefinite and where the idea of Cecotti-Vafa works 
will show the necessity of such upper inequalities
and also the necessity to go beyond $\Gamma^{Ber}_2$, so beyond
the variance.

Though in the cases where $S+S^t$ is positive semidefinite 
(including the positive definite cases) the situation simplifies
a lot, to the satisfying state in Theorem \ref{t1.5}.
There part (a) of Conjecture \ref{t1.4} is not a problem.
There the natural spectrum satisfies 
$-\frac{1}{2}\leq\alpha_1\leq ...
\leq \alpha_n\leq\frac{1}{2}$ and is uniquely defined by this and
the properties in part (a) of Conjecture \ref{t1.4}. 
In the positive definite cases one needs from \eqref{1.3}
only the lower inequality for $\Gamma_2^{Ber}$, so the
variance inequality. In the positive semidefinite, 
but not positive definite cases one needs the variance
inequality and the condition $\tr(S^{-1}S^t)=1$, which can
be considered as a limit for $k\to\infty$ of both 
inequalities in \eqref{1.3} (see Theorem \ref{t8.12} (c) and
Conjecture \ref{t8.17} (c)).
Theorem \ref{t1.5} is the main result of the paper. 

\begin{theorem}\label{t1.5}
Consider $S\in T_n(\Z)$ 
with connected Coxeter-Dynkin diagram ${\rm CDD}(S)$.

(a) (Theorem \ref{t6.2} (b)) 
Suppose that $S+S^t$ is positive definite.
Then the variance inequality 
\begin{eqnarray}\label{1.4}
\Var(\Sp(S)) \leq \frac{\alpha_n-\alpha_1}{12}
\end{eqnarray}
holds if and only if $S$ is a distinguished matrix
of an ADE-singularity.

(b) (Theorem \ref{t7.4})
Suppose that $S+S^t$ is positive semidefinite,
but not positive definite. Then the variance inequality 
\eqref{1.4} holds if and only if $S$ is in one of four
$\Br_n\ltimes\{\pm 1\}^n$ orbits, namely either
$S$ is a distinguished matrix of one of the three families 
$\www{E_6}$, $\www{E_7}$ and 
$\www{E_8}$ of simple elliptic singularities, 
or it is in one $\Br_6\ltimes\{\pm 1\}^6$ orbit called
$\www{D_4}$. In the fourth case $\tr(S^{-1}S^t)=2(\neq 1)$.
\end{theorem}

Section \ref{s2} recalls structures which are induced
by a matrix $S\in T_n(\Z)$ and which are explained
in detail in \cite{HL24}.

Section \ref{s3} recalls classical facts \cite{Bo68} \cite{Dy57}
\cite{Ca72} and more recent facts \cite{BaH20} around the 
irreducible  simply laced root lattices. The more recent facts 
from \cite{BaH20} concern nonreduced presentations of Weyl group 
elements by products of reflections such that the corresponding 
roots generate the whole root lattice. These nonreduced 
presentations are needed for the semidefinite cases in section 
\ref{s7}.

Section \ref{s4} recalls the definition of the 
{\it tubular elliptic
root systems} $D_4^{(1,1)},E_6^{(1,1)},E_7^{(1,1)}$ and 
$E_8^{(1,1)}$ and their Coxeter elements. It cites
facts from \cite{SaK85}, \cite{Kl87} and \cite{BaWY18}.
It gives a more conceptual characterization of the Coxeter
elements, which is needed in section \ref{s7}.

Section \ref{s5} recalls some facts on isolated hypersurface
singularities, their distinguished bases and their distinguished
matrices.

Section \ref{s6} treats the case of matrices $S\in T_n(\Z)$
with positive definite matrix $S+S^t$.

Section \ref{s7} treats the case of matrices $S\in T_n(\Z)$
with positive semidefinite, but not positive definite matrix $S+S^t$.

Section \ref{s8} discusses further the tower of conjectures
\ref{t1.4} and provides uncomfortable examples 
$S\in T_n^{S^1}(\Z)$ with indefinite matrix $S+S^t$.

\section{Unimodular bilinear lattices with triangular bases}
\label{s2}
\setcounter{equation}{0}
\setcounter{table}{0}
\setcounter{figure}{0}

The conceptual object behind a matrix $S\in T_n(\Z)$
is a {\it unimodular bilinear lattice with a triangular basis}.
These objects are studied systematically in \cite{HL24}.
In this section we recall their definition and a number 
of important properties and induced structures. 
The notions introduced here are fundamental and will be needed
in the later sections. 
For $S\in T_n^{S^1}(\Z)$ with $S+S^t$ positive (semi)definite,
Lemma \ref{t2.11} (c) gives a well-defined spectrum.
We are motivated by the case of isolated hypersurface 
singularities, where these objects are known since
long time.

\begin{definition}
[{\cite[Definition 2.3 and Definition 2.10]{HL24}}]\label{t2.1}

(a) A unimodular bilinear lattice is a pair $(H_\Z,L)$
where $H_\Z$ is a $\Z$-lattice, so a free $\Z$-module of
some finite rank $n\in\N=\{1,2,3,...\}$, and where 
$L:H_\Z\times H_\Z\to\Z$ is a $\Z$-bilinear form with
$\det L(\uuuu{e}^t,\uuuu{e})=1$ for some (or any)
$\Z$-basis $\uuuu{e}=(e_1,...,e_n)\in M_{1\times n}(H_\Z)$
of $H_\Z$. 

(b) For $(H_\Z,L)$ as in (a), a $\Z$-basis $\uuuu{e}$
is a {\it triangular basis} if $L(\uuuu{e}^t,\uuuu{e})^t
\in T_n(\Z)$.
The set of triangular bases is called $\BB^{tri}$. 

(c) Let $(H_\Z,L,\uuuu{e})$ be a unimodular bilinear lattice
with a triangular basis. The triple is called {\it reducible}
if a decomposition $\{1,...,n\}=I_1\ \dot\cup\ I_2$ with
$I_1\neq\emptyset$, $I_2\neq\emptyset$ and 
$$L(e_i,e_j)=L(e_j,e_i)=0\quad\textup{for }i\in I_1,j\in I_2,$$
exists. Otherwise the triple is called {\it irreducible}.
\end{definition}

\begin{remarks}\label{t2.2}
(i) Hubery and Krause \cite{HK16} 
defined the notion of a {\it bilinear lattice}.
There $(H_\Z,L)$ is as above except that the condition
$\det L(\uuuu{e}^t,\uuuu{e})=1$ is weakened to the condition
$\det L(\uuuu{e}^t,\uuuu{e})\neq 0$.

(ii) In Definition \ref{t2.1} (c), the choice
$L(\uuuu{e}^t,\uuuu{e})^t\in T_n(\Z)$ and not 
$L(\uuuu{e}^t,\uuuu{e})\in T_n(\Z)$ is motivated
by the case of isolated hypersurface singularities.

(iii) In the whole paper $(H_\Z,L)$ will denote a unimodular
bilinear lattice, and $(H_\Z,L,\uuuu{e})$ will denote a 
unimodular bilinear lattice with a triangular basis.

(iv) By far not every unimodular bilinear lattice has 
triangular bases. But here we care only about those which have.

(v) For fixed $n\in\N$ there is an obvious 1-1 correspondence
between the set 
$$\{(H_\Z,L,\uuuu{e})\,|\, (H_\Z,L,\uuuu{e})
\textup{ as in (iii)}\}/ \textup{isomorphism}$$
and the set $T_n(\Z)$, which is given by the map
$[(H_\Z,L,\uuuu{e})]\mapsto L(\uuuu{e}^t,\uuuu{e})^t$. 
To a given matrix $S\in T_n(\Z)$ we always associate
the corresponding triple $(H_\Z,L,\uuuu{e})$. 

(vi) The information in a matrix $S\in T_n(\Z)$ can
be encoded in a {\it Coxeter-Dynkin diagram}. 
This is often more convenient than writing down the matrix.
The Coxeter-Dynkin diagram ${\rm CDD}(S)$ is a graph with
$n$ vertices which are numbered from $1$ to $n$ and with
weighted edges which are defined as follows:
Between the vertices $i$ and $j$ with $i<j$ one has no edge
if $S_{ij}=0$ and an edge with weight $S_{ij}$ if $S_{ij}\neq 0$.
Alternatively, one draws $|S_{ij}|$ edges if $S_{ij}<0$ and
$S_{ij}$ dotted edges if $S_{ij}>0$. 

(vii) Obviously, a triple $(H_\Z,L,\uuuu{e})$ is irreducible
if and only if its Coxeter-Dynkin diagram 
${\rm CDD}(S)$ is connected 
(here $S=L(\uuuu{e}^t,\uuuu{e})^t$). 
\end{remarks}

[HL24] associates many objects to a triple $(H_\Z,L,\uuuu{e})$
as above. Here we need only a selection
of these objects, the following.

\begin{definition}[{\cite[Definition 2.3]{HL24}}]
\label{t2.3} 
Let $(H_\Z,L,\uuuu{e})$ be a unimodular bilinear lattice
of rank $n\in\N$ with a triangular basis $\uuuu{e}$
and with $S:=L(\uuuu{e}^t,\uuuu{e})^t\in T_n(\Z)$. 
The following objects are associated to it.

(a) A symmetric bilinear form 
$I^{(0)}:H_\Z\times H_\Z\to\Z$ and a skew-symmetric bilinear
form $I^{(1)}:H_\Z\times H_\Z\to\Z$ with 
\begin{eqnarray*}
I^{(0)}&=&L^t+L,\quad\textup{so } 
I^{(0)}(\uuuu{e}^t,\uuuu{e})=S+S^t,\\
I^{(1)}&=&L^t-L,\quad\textup{so } 
I^{(1)}(\uuuu{e}^t,\uuuu{e})=S-S^t,
\end{eqnarray*}
which are called {\it even} respectively {\it odd 
intersection form}.

(b) An automorphism $M:H_\Z\to H_\Z$ which is defined by
\begin{eqnarray*}
L(Ma,b)=L(b,a),\quad\textup{so }M(\uuuu{e}) = \uuuu{e}\cdot
S^{-1}S^t,
\end{eqnarray*}
and which is called {\it the monodromy}. It respects $L$
(and $I^{(0)}$ and $I^{(1)}$) because
$L(Ma,Mb)=L(Mb,a)=L(a,b)$. 

(c) Two automorphism groups
\begin{eqnarray*}
O^{(k)}&:=& \Aut(H_\Z,I^{(k)})\qquad\textup{for }k\in\{0;1\}.
\end{eqnarray*}

(d) The set of {\it roots}
\begin{eqnarray*}
R^{(0)}:=\{a\in H_\Z\,|\, L(a,a)=1\}=\{a\in H_\Z\,|\, 
I^{(0)}(a,a)=2\}, 
\end{eqnarray*}
and the set
\begin{eqnarray*}
R^{(1)}:= H_\Z.
\end{eqnarray*}

(e) For $k\in\{0;1\}$ and $a\in R^{(k)}$ the 
{\it reflection} (if $k=0$) respectively {\it transvection} 
(if $k=1$) $s^{(k)}_a\in O^{(k)}$ with 
\begin{eqnarray*}
s^{(k)}_a(b):= b-I^{(k)}(a,b)a\quad\textup{for }b\in H_\Z.
\end{eqnarray*}

(f) For $k\in\{0;1\}$ the {\it even} (if $k=0$) respectively 
{\it odd} (if $k=1$) {\it monodromy group}
\begin{eqnarray*}
\Gamma^{(k)}:=\langle s^{(k)}_{e_1},...,s^{(k)}_{e_n}
\rangle \subset O^{(k)}.
\end{eqnarray*}

(g) For $k\in\{0;1\}$ the set of {\it even} (if $k=0)$ 
respectively {\it odd} (if $k=1)$ {\it vanishing cycles}
\begin{eqnarray*}
\Delta^{(k)}:=\Gamma^{(k)}\{\pm e_1,...,\pm e_n\}
\subset R^{(k)}.
\end{eqnarray*}
\end{definition}

The following lemma is well known in the theory of 
isolated hypersurface singularities and in 
Picard-Lefschetz theory, see e.g. \cite{AGV88}, \cite{Eb01} or 
\cite[Lemma 2.6 and Theorem 2.7]{HL24}.

\begin{lemma}\label{t2.4}
Let $(H_\Z,L,\uuuu{e})$ be as in Definition \ref{t2.3}.
Let $k\in\{0;1\}$.

(a)
\begin{eqnarray}
I^{(k)}(a,b)&=& L((M+(-1)^k\id)(a),b)\quad\textup{for }
a,b\in H_\Z,\nonumber\\
\Rad I^{(k)}&=& \ker((M+(-1)^k\id):H_\Z\to H_\Z).
\label{2.1}
\end{eqnarray}

(b) For $g\in O^{(k)}$ and $a\in R^{(k)}$ 
\begin{eqnarray}\label{2.2}
g s_a^{(k)} g^{-1}=s_{g(a)}^{(k)}.
\end{eqnarray}

(c) $s_a^{(k)}\in \Gamma^{(k)}$ for $a\in \Delta^{(k)}$.

(d) Let $\uuuu{e}$ be a triangular basis. Then 
\begin{eqnarray}\label{2.3}
s^{(k)}_{e_1}...s^{(k)}_{e_n}=(-1)^{k+1}M.
\end{eqnarray} 
\end{lemma}

\begin{remarks}\label{t2.5}
The parts (a)--(e) of Definition \ref{t2.3} and the
parts (a)+(b) of Lemma \ref{t2.4} work for any
unimodular bilinear lattice $(H_\Z,L)$
and any $\Z$-basis $\uuuu{e}$, so $S\in SL_n(\Z)$.

The parts (f)+(g) of Definition \ref{t2.3} and part
(c) of Lemma \ref{t2.4} require $\uuuu{e}\in (R^{(0)})^n$
if $k=0$.  

Only Part (d) of Lemma \ref{t2.4} requires that $\uuuu{e}$
is a triangular basis. 
\end{remarks}

Let $n\in\Z_{\geq 2}$. The braid group $\Br_n$ of braids with
$n$ strings was introduced by E. Artin in \cite{Ar25}. 
He showed that
$\Br_n$ is generated by $n-1$ elementary braids
$\sigma_1,...,\sigma_{n-1}$ and that all relations between
them come from the relations 
\begin{eqnarray*}
\sigma_i\sigma_j&=&\sigma_j\sigma_i\quad\textup{for }
i,j\in\{1,...,n\}\textup{ with }|i-j|\geq 2,\\
\sigma_i\sigma_{i+1}\sigma_i&=&\sigma_{i+1}\sigma_i\sigma_{i+1}
\quad\textup{for }i\in\{1,...,n-2\}.
\end{eqnarray*}

An important action of $\Br_n$ is the 
{\it \index{Hurwitz action}Hurwitz action}
on the $n$-th power $G^n$ for any group $G$. 
The braid group $\Br_n$ acts via
\begin{eqnarray*}
\sigma_i(g_1,...,g_n)&:=& (g_1,...,g_{i-1},g_ig_{i+1}g_i^{-1},g_i, 
g_{i+2},...,g_n),\\
\sigma_i^{-1}(g_1,...,g_n)&:=& (g_1,...,g_{i-1},g_{i+1},
g_{i+1}^{-1}g_ig_{i+1},g_{i+2},...,g_n).
\end{eqnarray*}
The fibers of the map
\begin{eqnarray*}
\pi_n:G^n\to G,\quad \uuuu{g}=(g_1,...,g_n)\mapsto g_1...g_n,
\end{eqnarray*}
are invariant under this action.

In the case of a unimodular bilinear lattice $(H_\Z,L)$ 
and $k\in\{0;1\}$, the Hurwitz action on $(O^{(k)})^n$ 
restricts because of \eqref{2.2} to an action on the subset 
$(\{s_a^{(k)}\,|\, a\in R^{(k)}\})^n$ of $(O^{(k)})^n$.

If $\BB^{tri}\neq\emptyset$, 
this action has a natural lift, which is described in
Lemma \ref{t2.6}. Also this lemma is well known in the theory 
of isolated hypersurface singularities, 
see e.g. \cite{AGV88}, \cite{Eb01} or 
\cite[Section 3.2]{HL24}.

\begin{lemma}\label{t2.6}
Let $(H_\Z,L)$ be a unimodular bilinear lattice of rank
$n\in\Z_{\geq 2}$ with $\BB^{tri}\neq\emptyset$. Recall 
$\BB^{tri}\subset (R^{(0)})^n(\subset(R^{(1)})^n=H_\Z^n$).

(a) $\Br_n$ acts on $\BB^{tri}$ by
\begin{eqnarray*}
\sigma_j(\uuuu{v})&=& (v_1,...,v_{j-1},v_{j+1}-L(v_{j+1},v_j)v_j,
v_j,v_{j+2},...,v_n)\\
&=& (v_1,...,v_{j-1},s_{v_j}^{(k)}(v_{j+1}),
v_j,v_{j+2},...,v_n)\quad\textup{for }k\in\{0;1\},\\
\sigma_j^{-1}(\uuuu{v})&=& (v_1,...,v_{j-1},v_{j+1}, 
v_j-L(v_{j+1},v_j)v_{j+1},
v_{j+2},...,v_n)\\
&=& (v_1,...,v_{j-1},v_{j+1},(s_{v_{j+1}}^{(k)})^{-1}(v_j),
v_{j+2},...,v_n)\quad\textup{for }k\in\{0;1\}.
\end{eqnarray*}

(b) The multiplicative group $\{\pm 1\}^n$ is called 
{\sf sign group}. It is generated by the elements
$\delta_j=((-1)^{\delta_{ij}})_{i=1,...,n}$ for 
$j\in\{1,...,n\}$. It acts on $\BB^{tri}$ by 
\begin{eqnarray*}
\delta_j(\uuuu{v})&=& (v_1,...,v_{j-1},-v_j,v_{j+1},...,v_n).
\end{eqnarray*}

(c) The actions in (a) and (b) unite to an action on $\BB^{tri}$
of the semidirect product $\Br_n\ltimes\{\pm 1\}^n$ with
$\{\pm 1\}^n$ as normal subgroup, where
\begin{eqnarray*}
\sigma_j\delta_i\sigma_j^{-1}&=& \delta_i\quad\textup{for }
i\in\{1,...,n\}-\{j,j+1\},\\
\sigma_j\delta_j\sigma_j^{-1}&=&\delta_{j+1},\quad
\sigma_j\delta_{j+1}\sigma_j^{-1}=\delta_j.
\end{eqnarray*}
In the following $\Br_n\ltimes\{\pm 1\}^n$ always means
this semidirect product.

(d) Fix $k\in\{0;1\}$. The map
\begin{eqnarray*}
\BB^{tri}\to (O^{(k)})^n,\quad \uuuu{v}\mapsto
(s_{v_1}^{(k)},...,s_{v_n}^{(k)}),
\end{eqnarray*}
is $\Br_n$ equivariant, where the action on $\BB^{tri}$
is the one in part (a) and the action on $(O^{(k)})^n$ is the
Hurwitz action. The action of $\{\pm 1\}^n$ on $\BB^{tri}$
maps to the trivial action on $(O^{(k)})^n$. 

(e) The group $\Br_n\ltimes\{\pm 1\}^n$ acts on the set
$T_n(\Z)$ by
\begin{eqnarray*}
\sigma_j(A)&=& C_{n,j}(-A_{j,j+1})\cdot A\cdot C_{n,j}(-A_{j,j+1})
\quad\textup{for }j\in\{1,...,n-1\},\\
\sigma_j^{-1}(A)&=& C_{n,j}^{-1}(A_{j,j+1})\cdot A\cdot C_{n,j}^{-1}(A_{j,j+1})
\quad\textup{for }j\in\{1,...,n-1\},\\
\delta_j(A)&=& \diag(((-1)^{\delta_{ij}})_{i=1,...,n})\cdot A\cdot 
\diag(((-1)^{\delta_{ij}})_{i=1,...,n})\\
&&\hspace*{4cm} \textup{for }j\in\{1,...,n\},
\end{eqnarray*}
for $A\in T_n(\Z)$. 
Here $C_{n,j}(a)$ for $a\in\Z$ is the $n\times n$ matrix
\begin{eqnarray*}
C_{n,j}(a)&=&
\begin{pmatrix}
1 & & & & & \\ & \ddots & & & & \\ 
& & a & 1 & & \\
& & 1 & 0 & & \\ & & & & \ddots & \\ & & & & & 1 \end{pmatrix}
\end{eqnarray*}
which differs from the unit matrix only in the positions
$(j,j),(j,j+1),(j+1,j),(j+1,j+1)$. 

(f) The map 
\begin{eqnarray*}
\BB^{tri}\to T_n(\Z),\quad \uuuu{v}\mapsto 
L(\uuuu{v}^t,\uuuu{v})^t,
\end{eqnarray*}
is $\Br_n\ltimes\{\pm 1\}^n$ equivariant. 
\end{lemma}

\begin{definition}[{\cite[Definition 3.18]{HL24}}]\label{t2.7}
Let  $(H_\Z,L,\uuuu{e})$ be a unimodular bilinear lattice
of rank $n\in\Z_{\geq 2}$ 
with a triangular basis $\uuuu{e}$ and associated matrix
$S:=L(\uuuu{e}^t,\uuuu{e})^t\in T_n(\Z)$.
The $\Br_n\ltimes\{\pm 1\}^n$ orbit
$$\BB^{dist}:=\Br_n\ltimes\{\pm 1\}^n(\uuuu{e})\subset\BB^{tri}$$
is called {\sf set of distinguished bases} of the triple
$(H_\Z,L,\uuuu{e})$. 
The $\Br_n\ltimes\{\pm 1\}^n$ orbit
$$\SSS^{dist}:=\Br_n\ltimes\{\pm 1\}^n(S)\subset T_n(\Z)$$
is called {\sf set of distinguished matrices}.
\end{definition}

\begin{remarks}\label{t2.8}
(i) Let $(H_\Z,L,\uuuu{e})$ and $S$ be as in Definition \ref{t2.7}.
Any distinguished basis $\uuuu{v}$ is a triangular basis
and thus satisfies $s_{v_1}^{(k)}...s_{v_n}^{(k)}=(-1)^{k+1}M$
for $k\in\{0;1\}$, and it is in  
$(\Delta^{(0)})^n$ and in $(\Delta^{(1)})^n$.
Especially, there are inclusions
\begin{eqnarray}\label{2.4}
\BB^{dist}&\subset& \{\uuuu{v}\in (\Delta^{(0)})^n\,|\, 
s_{v_1}^{(0)}...s_{v_n}^{(0)}=-M,
\sum_{i=1}^n\Z\cdot v_i=H_\Z\},\\
\BB^{dist}&\subset& \{\uuuu{v}\in (\Delta^{(1)})^n\,|\, 
s_{v_1}^{(1)}...s_{v_n}^{(1)}=M,
\sum_{i=1}^n\Z\cdot v_i=H_\Z\}.\label{2.5}
\end{eqnarray}
It is interesting to ask when these inclusions are equalities.
Equality in \eqref{2.4} holds, even without the condition 
$\sum_{i=1}^n\Z\cdot v_i=H_\Z$, for the ADE-singularities
\cite{De74}. This generalizes on the level of tuples of
reflections greatly, namely 
to all Coxeter groups \cite{IS10}\cite{BaDSW14}.
Also in the cases of the simple elliptic and the hyperbolic
singularities \eqref{2.4} holds \cite[VI 1.2 Theorem]{Kl87}
\cite[Theorem 1.5]{BaWY21} (see Theorem \ref{t4.11} and the
Remarks \ref{t5.6} (iv) and (v)).
Equality in \eqref{2.5} holds  for the 
singularities $A_2$ and $A_3$. 
It is unknown whether \eqref{2.4} or \eqref{2.5} hold 
for all isolated hypersurface singularities. 
\cite{HL24} studies the inclusions \eqref{2.4} and 
\eqref{2.5} systematically for all cases with $n=2$ and $n=3$.

(ii) A unimodular bilinear lattice $(H_\Z,L)$ with
$\BB^{tri}\neq\emptyset$ arises in a canonical way in 
geometric situations (e.g. isolated hypersurface singularities).
But a basis $\uuuu{e}\in\BB^{tri}$ comes from some choice.
Though we claim that the $\Br_n\ltimes\{\pm 1\}^n$ orbit
$\BB^{dist}=\Br_n\ltimes\{\pm 1\}^n(\uuuu{e})\subset\BB^{tri}$
is again a rather canonical object. This claim is strengthened
by the fact that $\Gamma^{(k)}$ and $\Delta^{(k)}$
depend on the set $\BB^{dist}$, but not on the choice of a
basis $\uuuu{e}$ in $\BB^{dist}$. This fact is a consequence of
Lemma \ref{t2.4} (b)+(c). 
\end{remarks}

The following observation fixes a situation where the inclusion
in \eqref{2.4} respectively \eqref{2.5} is not an equality.
It was used in  \cite[Example 3.23 (ii)]{HL24} for the elliptic 
root lattice $A_1^{(1,1)*}$ (which is called $\HH_{1,2}$ there). 
Below in Example \ref{t4.13} and in Theorem \ref{t7.4} 
it is used for the tubular elliptic root system $D_4^{(1,1)}$. 

\begin{lemma}\label{t2.9}
Let $(H_\Z,L,\uuuu{e})$ and $S$ be as in Definition \ref{t2.7}.
Fix $k\in\{0;1\}$. Let $g$ be an automorphism of $H_\Z$ which
respects $I^{(k)}$ and which commutes with $M$, but which does
not respect $L$, so
$$g\in\Aut(H_\Z,I^{(k)},M)-\Aut(H_\Z,L).$$
Then $g(\uuuu{e}):=(g(e_1),...,g(e_n))$ is in the right hand
side of \eqref{2.4} (if $k=0$) respectively \eqref{2.5} 
(if $k=1$), but not in $\BB^{tri}$, and especially not in
the left hand side $\BB^{dist}$ of \eqref{2.4} respectively
\eqref{2.5}.
\end{lemma}

{\bf Proof:}
The condition $\sum_{i=1}^n\Z\cdot g(e_i)=H_\Z$ is satisfied
because $g$ is an automorphism of $H_\Z$. We also have
\begin{eqnarray*}
(-1)^{k+1}M&\stackrel{gM=Mg}{=}&g(-1)^{k+1}Mg^{-1}\\
&=&gs^{(k)}_{e_1} ... s^{(k)}_{e_n}g^{-1}\\
&=&(gs^{(k)}_{e_1}g^{-1}) ... (gs^{(k)}_{e_n}g^{-1})\\
&\stackrel{\eqref{2.2}}{=}& s^{(k)}_{g(e_1)}...
s^{(k)}_{g(e_n)}.
\end{eqnarray*}
That $g$ respects $I^{(k)}$, gives the third of the following
equalities.
\begin{eqnarray*}
S+(-1)^kS^t&=& (L^t+(-1)^kL)(\uuuu{e}^t,\uuuu{e})
=I^{(k)}(\uuuu{e}^t,\uuuu{e})\\
&=& I^{(k)}(g(\uuuu{e})^t,g(\uuuu{e}))
=(L^t+(-1)^kL)(g(\uuuu{e})^t,g(\uuuu{e}))\\
&=& L(g(\uuuu{e})^t,g(\uuuu{e}))^t
+(-1)^kL(g(\uuuu{e})^t,g(\uuuu{e})).
\end{eqnarray*}
But because $g$ does not respect $L$, we have
$$S=L(\uuuu{e}^t,\uuuu{e})^t\neq L(g(\uuuu{e})^t,g(\uuuu{e}))^t.
$$
Therefore the matrix $L(g(\uuuu{e})^t,g(\uuuu{e}))^t$ is not
in $T_n(\Z)$, and thus $g(\uuuu{e})$ is not in $\BB^{tri}$.
\hfill$\Box$

\bigskip
The following Lemma of Ebeling applies to the case
of isolated hypersurface singularities because of
Theorem \ref{t1.2} (d). 

\begin{lemma}[{\cite[Proposition 2.2]{Eb83}}]
\label{t2.10} 
Let $(H_\Z,L,\uuuu{e})$ and $S$ be as in Definition \ref{t2.7}.
If all entries of $S$ are in the set $\{0;1;-1\}$
then $\BB^{dist}=\Br_n(\uuuu{e})$ and 
$\SSS^{dist}=\Br_n(S)$, so any sign change can be carried out
by a suitable braid.
\end{lemma}

Here are some first facts on the cases when
$S+S^t$ is positive definite or positive semidefinite.
Especially, for such $S$ Lemma \ref{t2.11} (c)
gives a well-defined spectrum of $S$.

\begin{lemma}\label{t2.11}
Let $(H_\Z,L,\uuuu{e})$ be a unimodular bilinear lattice
with a triangular basis and let $S:=L(\uuuu{e}^t,\uuuu{e})^t
\in T_n(\Z)$ ¸be the associated matrix.

(a) Suppose that $S+S^t$ is positive definite. 

(i) $S_{ij}\in\{0;\pm 1\}$. The sets $R^{(0)}$, 
$\BB^{dist}$ (if $n\geq 2$) and $\SSS^{dist}$ are finite.

(ii) All eigenvalues of $M$ are roots of unity. 
$-1$ is not an eigenvalue of $M$.

(b) Suppose that $S+S^t$ is positive semidefinite.

(i) $S_{ij}\in\{0;\pm 1;\pm 2\}$. 
The set $\SSS^{dist}$ is finite.

(ii) All eigenvalues of $M$ are roots of unity. 
The multiplicity if $-1$ as a generalized eigenvalue of $M$ 
(that means, the multiplicity of $-1$ as a zero of the 
characteristic polynomial) is even.

(c) Let $S+S^t$ be positive definite or positive semidefinite.
There are unique numbers $\alpha_1,...,\alpha_n\in 
[-\frac{1}{2},\frac{1}{2}]\cap\Q$ with
\begin{eqnarray*} 
\alpha_1\leq ...\leq \alpha_n,\\
e^{-2\pi i\alpha_1},...,e^{-2\pi i\alpha_n}
\quad\textup{the eigenvalues of }M,\\ 
|\{i\,|\, \alpha_i=-\frac{1}{2}\}|
=|\{i\,|\, \alpha_i=\frac{1}{2}\}|.
\end{eqnarray*}
They satisfy the symmetry $\alpha_{n+1-i}=\alpha_i$. 
The tuple $(\alpha_1,...,\alpha_n)$ is called 
{\sf spectrum of $S$} and denoted $\Sp(S)$. 
The numbers $\alpha_1,...,\alpha_n$ are called 
{\sf spectral numbers} of $S$.
\end{lemma}

{\bf Proof:} The parts (a)(i) and (b)(i): 
If $S+S^t$ is positive semidefinite, also any submatrix
$\begin{pmatrix}2 & S_{ij} \\ S_{ij} & 2\end{pmatrix}$
(with $i<j$) is positive semidefinite. Therefore then 
$S_{ij}\in\{0,\pm 1,\pm 2\}$. Therefore the set 
$\SSS^{dist}$ is finite.

If $S+S^t$ is positive definite, the same argument gives
$S_{ij}\in\{0,\pm 1\}$. Because $I^{(0)}$ is positive 
definite, the set $R^{(0)}$ is finite. The set 
$\BB^{dist}$ is a subset of the finite set $(R^{(0)})^n$.

Part (a)(ii): As $M$ maps $R^{(0)}$ to itself and $R^{(0)}$ is
finite, $M$ has finite order, so its eigenvalues are roots of unity. 
\eqref{2.1}, namely $\{0\}=\Rad I^{(0)}=\ker(M+\id)$, 
shows that $-1$ is not an eigenvalue of $M$. 

Part (b)(ii): $I^{(0)}$ induces a positive definite pairing
on the quotient $H_\Z/\Rad I^{(0)}$. Therefore $M$ acts on it 
as an automorphism of finite order, so its eigenvalues on it
are roots of unity.
The pairing $I^{(0)}$ on $H_\C:=H_\Z\otimes_\Z\C$
is $M$-invariant, so the generalized eigenspaces with
eigenvalues $\lambda$ and $\lambda^{-1}$ are dual and thus
have the same dimension. As $1=\det M=\det S^{-1}S^t$, 
the product of all eigenvalues with their multiplicities 
(as generalized eigenvalues) is 1. Therefore the multiplicity 
of $-1$ as a generalized eigenvalue of $M$ is even.

Part (c): The existence and uniqueness of $\alpha_1,...,\alpha_n$
follows from the parts (a)(ii) and (b)(ii).
The symmetry follows from the proof of part (b)(ii).
\hfill$\Box$ 

\bigskip

In the case when $S+S^t$ is indefinite, the first half of
section \ref{s8} discusses an incomplete recipe of
Cecotti and Vafa \cite{CV93} and a solution in special
cases in \cite{BaH20} how to associate a spectrum
to a matrix $S\in T_n^{S^1}(\R)$.

We conclude this section with a comparison of 
unimodular bilinear lattices and {\it simply laced
generalized root systems}, which were defined essentially
by K. Saito \cite{SaK85} \cite[(5.3)]{SaK98}. 
The following Definition \ref{t2.12} 
is taken from \cite{STW16}. Lemma \ref{t2.13} and Remark
\ref{t2.14} make the comparison.

\begin{definition}[{\cite[Definition 2.1]{STW16}}]\label{t2.12}
A {\it simply laced generalized root system} is a tuple
$(H_\Z,I,\Delta_{re},c)$ with the following ingredients.
$H_\Z$ is a $\Z$-lattice of some rank $n\in\N$.
$I:H_\Z\times H_\Z\to\Z$ is a symmetric bilinear form.
The set $\Delta_{re}$ of {\it real roots} is a subset of the
set $\{a\in H_\Z\,|\, I(a,a)=2\}$, and $c\in\Aut(H_\Z,I)$ is
a certain automorphism (called {\it Coxeter transformation} in
\cite{STW16}). These data are subject to the following conditions.
Here $s_a\in\Aut(H_\Z,I)$ for $a\in H_\Z$ with $I(a,a)=2$ 
is the reflection along $a$ and is defined as in Definition 
\ref{t2.3} (e) (with $k=0$ and $I$ instead of $I^{(0)}$).

A $\Z$-basis $\uuuu{e}=(e_1,...,e_n)\subset\Delta_{re}^n$ of
$H_\Z$ exists such that 
\begin{eqnarray*}
c&=&s_{e_1}...s_{e_n},\\
W:= \langle s_a\,|\, a\in \Delta_{re}\rangle 
&\stackrel{!}{=}&\langle s_{e_1},...,s_{e_n}\rangle,\\
\Delta_{re}&=& W\{\pm e_1,...,\pm e_n\}.
\end{eqnarray*}
\end{definition}

\begin{lemma}\label{t2.13}
(a) Let $(H_\Z,L)$ be a unimodular bilinear lattice with 
$\BB^{tri}\neq\emptyset$. Then the tuple 
$(H_\Z,I^{(0)},\Delta^{(0)},-M)$ is a simply laced generalized
root system.

(b) Let $(H_\Z,I,\Delta_{re},c)$ be a simply laced generalized
root system. Choose a $\Z$-basis $\uuuu{e}=(e_1,...,e_n)$
with the properties in Definition \ref{t2.12}. Define
a bilinear form $L:H_\Z\times H_\Z\to\Z$ by 
$L(\uuuu{e}^t,\uuuu{e})^t=S$ where $S\in T_n(\Z)$ is defined by
$S+S^t=I(\uuuu{e}^t,\uuuu{e})$. Then
$(H_\Z,L,\uuuu{e})$ is a unimodular bilinear lattice with 
triangular basis $\uuuu{e}$ and $I^{(0)}=I$, $\Gamma^{(0)}=W$, 
$\Delta^{(0)}=\Delta_{re}$ and $-M=c$.

(c) Let $(H_\Z,I,\Delta_{re},c)$ be a simply laced generalized
root system. 

(i) A bilinear form $L$ as in part (b) satisfies
\begin{eqnarray}\label{2.6}
I=L+L^t,\quad L(-ca,b)=L(b,a)\textup{ for }a,b\in\Z,\\
I(a,b)=L((\id-c)a,b)\textup{ for }a,b\in\Z.\label{2.7}
\end{eqnarray}

(ii) $\Rad(I)=\ker(c-\id)\subset H_\Z$.

(iii) If $\Rad(I)=\ker(c-\id)=\{0\}$ then a bilinear
form $L: H_\Z\times H_\Z\to\Z$ with \eqref{2.7} is unique.
So then the bilinear form in part (b) is unique.
\end{lemma}

{\bf Proof:} (a) Choose any triangular basis $\uuuu{e}$.
Then $-M=s^{(0)}_{e_1}...s^{(0)}_{e_n}$ by \eqref{2.3},
$\Gamma^{(0)}=\langle s^{(0)}_a\,|\, a\in\Delta^{(0)}\rangle$
by Lemma \ref{t2.4} (c) and $\Delta^{(0)}=\Gamma^{(0)}\{\pm e_1,...,
\pm e_n\}$ by definition.

(b) Obviously $I=L+L^t=I^{(0)}$, 
$\Gamma^{(0)}=\langle s_{e_1},...,s_{e_n}\rangle = W$, 
$\Delta^{(0)}=\Delta_{re}$ and 
$-M\stackrel{\textup{\eqref{2.3}}}{=}s_{e_1}...s_{e_n}=c$.
\hfill$\Box$ 

(c) (i) \eqref{2.6} follows from $I=I^{(0)}$ and $c=-M$ and 
the definitions of $I^{(0)}$ and $M$ in a unimodular bilinear
lattice. \eqref{2.6} implies \eqref{2.7}.

(ii) A bilinear form $L$ in part (b) exists and satisfies 
\eqref{2.7}. Because $L$ is nondegenerate, 
\eqref{2.7} implies $\Rad(I)=\ker(c-\id)$.

(iii) If $\ker(c-\id)=\{0\}$ then $I$, $c$ and \eqref{2.7}
determine $L$ uniquely.\hfill$\Box$

\begin{remarks}\label{t2.14}
Therefore the notions of a {\it unimodular bilinear lattice
with $\BB^{tri}\neq\emptyset$} and of a {\it simply laced
generalized root system} are closely related. 
By part (c)(iii) they are equivalent if $c$ 
in a simply laced generalized root system 
satisfies $\ker(c-\id)=\{0\}$. 

But if $\ker(c-\id)\neq\{0\}$, then several unimodular 
bilinear lattices with triangular bases might induce 
the same simply laced generalized root system.
In part (b) of Lemma \ref{t2.13} 
the $\Z$-basis $\uuuu{e}$ was chosen. A different choice
$\www{\uuuu{e}}$ (with the same properties) 
might lead to a different bilinear form $\www{L}$
if $\ker(c-\id)\neq\{0\}$. 

The proof of Lemma \ref{t2.9} gives an idea how to construct 
examples. If there $k=0$ then $\uuuu{e}$ and 
$\www{\uuuu{e}}:=g(\uuuu{e})$ satisfy
\begin{eqnarray*}
s^{(0)}_{e_1}...s^{(0)}_{e_n}&=& s^{(0)}_{\www{e}_1}... 
s^{(0)}_{\www{e}_n},\\
\www{\Gamma}^{(0)}&=& g\Gamma^{(0)}g^{-1},\\
\www{\Delta}^{(0)}&=& g(\Delta^{(0)}).
\end{eqnarray*}
We need $\www{\Gamma}^{(0)}= \Gamma^{(0)}$ and 
$\www{\Delta}^{(0)}= \Delta^{(0)}$.

In Example \ref{t4.13} we know 
$$\Delta^{(0)}\stackrel{\textup{Theorem \ref{t4.3}}}{=}
R^{(0)}=g(R^{(0)})\stackrel{\textup{Theorem \ref{t4.3}}}{=}
g(\Delta^{(0)})=\www{\Delta}^{(0)}$$ 
and thus $\www{\Gamma}^{(0)}=\Gamma^{(0)}$. So in this example
$\uuuu{e}$ and $\www{\uuuu{e}}:=g(\uuuu{e})$ both satisfy the
properties in Definition \ref{t2.12}. But $L$ and
$\www{L}$ are different. 

Example \ref{t4.13} is on the simply laced generalized
root system $D_4^{(1,1)}$ (with a chosen Coxeter element).
The similar simply laced generalized root systems
$E_6^{(1,1)},E_7^{(1,1)}$ and $E_8^{(1,1)}$ (each with a 
chosen Coxeter element $c$) behave better. By Theorem \ref{t4.11}
there $L$ is unique, although also there $\ker(c-\id)\neq\{0\}$. 
\end{remarks}

\section{Review of results on the simply laced root lattices}
\label{s3}
\setcounter{equation}{0}
\setcounter{table}{0}
\setcounter{figure}{0}

This section recalls several facts on the irreducible simply laced
root lattices, especially results of Dynkin \cite{Dy57}
on the classification of their subroot lattices,
results of Carter \cite{Ca72} on the classification of 
conjugacy classes of elements in the Weyl group and 
results of Balnojan and Hertling \cite{BaH19} on 
nonreduced presentations of these elements.

All the results are needed in the proofs in section 
\ref{s6} and section \ref{s7} of Theorem \ref{t1.3} and \ref{t1.5}.
The quasi-Coxeter elements (Definition \ref{t3.8} (c) and Theorem
\ref{t3.10}) are needed for the positive definite cases.
The subroot lattices (Definition \ref{t3.5} and Theorem \ref{t3.6})
and the non-reduced presentations of Weyl group elements
(Theorem \ref{t3.13}) are needed for the positive semidefinite 
cases. They arise in the quotient $H_\Z/\Rad I^{(0)}$
if $(H_\Z,L)$ is a unimodular bilinear lattice with triangular
bases and with positive semidefinite form $I^{(0)}$.

The initial data in the sections \ref{s2} and \ref{s3} 
are different. Therefore also some similar induced objects are 
denoted differently ($\Delta\subset R^{(0)}$ and $\Gamma^{(0)}$ 
in section \ref{s2}, $\Phi$ and $W$ in section \ref{s3}). 
A relation is established only in Lemma \ref{t6.1}.

\begin{definition}[{E.g. \cite{Bo68}}]\label{t3.1}
(a) A {\it simply laced root lattice} is a tuple $(H_\Z,I)$, 
where $H_\Z$ is free $\Z$-module ($=\Z$-lattice) of rank 
$n\in\N$, $I:H_\Z\times H_\Z\to\Z$ is a $\Z$-bilinear symmetric 
and positive definite form, and the finite set $
\Phi:=\{\delta\in H_\Z\,|\, I(\delta,\delta)=2\}$
generates $H_\Z$, so $H_\Z=\sum_{\delta\in\Phi}\Z\cdot\delta$. 
Then for $\delta\in\Phi$, the automorphism
\begin{eqnarray*}
s_\delta:H_\Z\to H_\Z\quad\textup{with}\quad 
s_\delta(b)&:=& b-I(\delta,b)\delta
\end{eqnarray*} 
is a reflection, it maps $H_\Z$ to $H_\Z$ and $\Phi$ to $\Phi$. 
The elements of $\Phi$
are the {\it roots}, and $\Phi$ ist the {\it root system}.
The finite group $$W:=\langle s_\delta\, |\, \delta\in\Phi\rangle\subset
\textup{Aut}(H_\Z,I)$$
is the {\it Weyl group}.

(b) A {\it root basis} of a simply laced root lattice $(H_\Z,I)$
is a set $\{\alpha_1,...,\alpha_n\}\subset\Phi$ such that
$H_\Z=\bigoplus_{j=1}^n\Z\cdot \alpha_j$ and such that any root
is a linear combination of $\alpha_1,...,\alpha_n$ with either
all coefficients nonnegative or all coefficients nonpositive.
\end{definition}

\begin{lemma}[{E.g. \cite{Bo68}}]\label{t3.2}
The orthogonal sum of two simply laced 
root lattices $(H_\Z^{(i)},I^{(i)})$
$(i=1,2)$ is a root lattice 
$(H^{(1)}_\Z\oplus H^{(2)}_\Z,I^{(1)}\oplus I^{(2)})$
with root system $\Phi^{(1)}\cup \Phi^{(2)}$. 
\end{lemma}

\begin{definition}[{E.g. \cite{Bo68}}]\label{t3.3}
A simply laced root lattice is {\it irreducible} if it is not
isomorphic to an orthogonal sum of two simply laced root lattices.
\end{definition}

\begin{theorem}[{E.g. \cite{Bo68}}]\label{t3.4}
(a) Any simply laced root lattice is isomorphic to a unique orthogonal
sum of irreducible simply laced root lattices.

\medskip
(b) The irreducible simply laced root lattices come in two series
and three separate cases, they are denoted by
$A_n\ (n\geq 1)$, $D_n\ (n\geq 4)$, $E_6,E_7,E_8$.

\medskip
(c) Any irreducible simply laced root lattice has root bases.
The Weyl group acts simply transitively on them.
If $\{\alpha_1,...,\alpha_n\}$ is a root basis, then
its {\sf Cartan matrix} is the matrix
$C=(c_{ij})=(I(\alpha_i,\alpha_j))$. It satisfies 
$$c_{ii}=2\quad\textup{and}\quad c_{ij}\in\{0,-1\}\textup{ for }i\neq j.$$
The Dynkin diagram of the Cartan matrix $C$ is by definition
the diagram ${\rm CDD}(S)$ where $C=S+S^t$ with $S\in T_n(\Z)$
(see Remark \ref{t2.2} (vi) for ${\rm CDD}(S)$).
The Dynkin diagram is up to the numbering of the vertices the
same for all root bases. It is given (without numbering
of the vertices) in Figure \ref{fig3.1}.
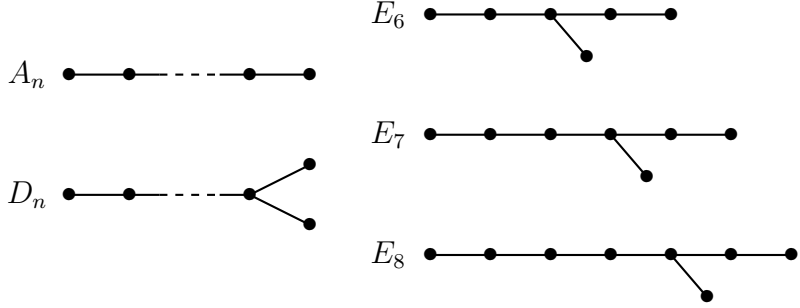
\begin{figure}
\begin{tikzpicture}[scale=0.8]
\node at (0,2) {$\bullet$};
\node at (1,2) {$\bullet$};
\draw [thick] (0,2)--(1.5,2);
\draw [thick,dashed] (1.5,2)--(2.5,2);
\draw [thick] (2.5,2)--(4,2);
\node at (3,2) {$\bullet$};
\node at (4,2) {$\bullet$};
\node at (-0.7,2) {$A_n$};
\node at (0,0) {$\bullet$};
\node at (1,0) {$\bullet$};
\draw [thick] (0,0)--(1.5,0);
\draw [thick,dashed] (1.5,0)--(2.5,0);
\draw [thick] (2.5,0)--(3,0);
\node at (3,0) {$\bullet$};
\node at (4,0.5) {$\bullet$};
\node at (4,-0.5) {$\bullet$};
\draw [thick] (3,0)--(4,0.5);
\draw [thick] (3,0)--(4,-0.5);
\node at (-0.7,0) {$D_n$};
\draw [thick] (6,3)--(10,3);
\node at (6,3) {$\bullet$};
\node at (7,3) {$\bullet$};
\node at (8,3) {$\bullet$};
\node at (9,3) {$\bullet$};
\node at (10,3) {$\bullet$};
\node at (8.6,2.3) {$\bullet$};
\draw [thick] (8,3)--(8.6,2.3);
\node at (5.3,3) {$E_6$};
\draw [thick] (6,1)--(11,1);
\node at (6,1) {$\bullet$};
\node at (7,1) {$\bullet$};
\node at (8,1) {$\bullet$};
\node at (9,1) {$\bullet$};
\node at (10,1) {$\bullet$};
\node at (11,1) {$\bullet$};
\node at (9.6,0.3) {$\bullet$};
\draw [thick] (9,1)--(9.6,0.3);
\node at (5.3,1) {$E_7$};
\draw [thick] (6,-1)--(12,-1);
\node at (6,-1) {$\bullet$};
\node at (7,-1) {$\bullet$};
\node at (8,-1) {$\bullet$};
\node at (9,-1) {$\bullet$};
\node at (10,-1) {$\bullet$};
\node at (11,-1) {$\bullet$};
\node at (12,-1) {$\bullet$};
\node at (10.6,-1.7) {$\bullet$};
\draw [thick] (10,-1)--(10.6,-1.7);
\node at (5.3,-1) {$E_8$};
\end{tikzpicture}
\caption[Figure 3.1]{Dynkin diagrams of the ADE root lattices}
\label{fig3.1}
\end{figure}
\end{theorem}

\begin{definition}\label{t3.5}
Let $(H_\Z,I)$ be a simply laced root lattice of rank $n\in\N$ 
with root system $\Phi$.

(a) A {\it subroot lattice} $U\subset H_\Z$ is a sublattice
of $H_\Z$ which is generated by roots. 

Remark: Then $(U,I|_{U\times U})$ is a root lattice with
root system $\Phi\cap U$.

(b) Let $U\subset H_\Z$ be a subroot lattice. Define
\begin{eqnarray*}
k_1&:=&k_1(H_\Z,U):=\min(k\,|\, \textup{the group }H_\Z/U
\textup{ has }k\textup{ generators}\},\\
i_1&:=&i_1(H_\Z,U):=[(H_\Z\cap U_\Q):U]\in\N,
\end{eqnarray*}
where $U_\Q:=U\otimes_\Z\Q$. The number 
$i_1$ is called the {\it index of $U$ in $H_\Z$}.
It is the size of the torsion part of $H_\Z/U$. 
\end{definition}

A recipe of Borel and de Siebenthal \cite{BdS49} and of Dynkin 
\cite{Dy57} allows to determine for each irreducible
simply laced root lattice $(H_\Z,I)$ the isomorphism classes
of all pairs $((H_\Z,I),U)$ with $U$ a subroot lattice of 
$(H_\Z,I)$. 
The recipe is recalled in \cite[Remark 3.2 (iv)]{BaH19}.
Dynkin \cite{Dy57} gives also a table and additional 
information. In \cite{BaH19} this information is enriched
by a calculation of the quotient groups $H_\Z/U$ and the
numbers $k_1$ and $i_1$.
Theorem \ref{t3.6} gives Dynkin's table and the additional
information $H_\Z/U$ and $k_1$ from \cite[Theorem 3.3]{BaH19}.
For complete information one has to go through
the construction of Borel-de Siebenthal and Dynkin.
This is documented in \cite{BaH19} in the Tables 3.2-3.6 
and in Lemma 3.7.

\begin{theorem}\label{t3.6}
Let $(H_\Z,I)$ be an irreducible simply laced root system.
The Tables \ref{tab3.1}, \ref{tab3.2}, \ref{tab3.3} and 
\ref{tab3.4} give names and some information for all
isomorphism classes of pairs $((H_\Z,I),U)$ with $U$
a subroot lattice. Table \ref{tab3.1} treats $(H_\Z,I)$
of type $A_n$ and $D_n$, the Tables \ref{tab3.2},
\ref{tab3.3} and \ref{tab3.4} treat $(H_\Z,I)$ of type
$E_6$ respectively $E_7$ respectively $E_8$. 

(a) $(H_\Z,I)$ of type $A_n$ or $D_n$:
The entries in Table \ref{tab3.1} are explained after
the table.
\begin{table}[H]
\begin{eqnarray*}
\begin{array}{l|l|l|ll}
H_\Z & U & H_\Z/U & k_1 & \\ \hline 
A_n & \sum_{i=1}^{r_1}A_{c_i} & \Z^{r_1+r_2} & r_1+r_2\\
D_n & \sum_{i=1}^{r_1}A_{c_i}+\sum_{j=1}^{r_3}D_{b_j} & 
\Z^{r_1+r_2}\times \Z_2^{r_3-1} & \sum_{i=1}^3r_i-1
&\textup{if }r_3\geq 1\\
 & & \Z^{r_1+r_2} & r_1+r_2 & \textup{if }r_3=0
\end{array}
\end{eqnarray*}
\caption[Table 3.1]{Isomorphism classes of pairs
$((H_\Z,I),U)$ for $(H_\Z,I)$ of type $A_n$ and $D_n$}
\label{tab3.1}
\end{table}
\noindent
In the case of $A_n$ 
\begin{eqnarray*}
r_1\in\N,\ r_2\in\Z_{\geq -1},\ c_i\in\N, 
\ \sum_{i=1}^{r_1}c_i=n-r_1-r_2.
\end{eqnarray*}
In the case of $D_n$
\begin{eqnarray*}
r_1,r_2,r_3\in\Z_{\geq 0},\ r_1+r_3\geq 1,
\ c_i\in\N,\ b_i\in\Z_{\geq 2},
\ \sum_{i=1}^{r_1}c_i+\sum_{j=1}^{r_3}b_j=n-r_1-r_2.
\end{eqnarray*}
In the case of $D_n$, the notations $D_2$ and $D_3$ 
in the second sum $\sum_{j=1}^{r_3}D_{b_j}$ mean subroot 
lattices of type $2A_1$ respectively $A_3$.
But the possible summands $D_2$ and $D_3$
in the second sum $\sum_{j=1}^{r_3}D_{b_j}$ 
indicate an embedding different from that of possible summands
$2A_1$ and $A_3$ in the first sum $\sum_{i=1}^{r_1}A_{c_i}$.
For details see \cite[\S 5]{Dy57} or \cite[Theorem 3.3 and Lemma
3.6]{BaH19}.

\begin{table}
\begin{eqnarray*}
\begin{array}{l|l|l}
U & H_\Z/U & k_1\\ \hline
E_6 & \{0\} & 0 \\
A_5+A_1 & \Z_2 & 1 \\
3A_2 & \Z_3 & 1 \\ \hline
A_5 & \Z & 1 \\
2A_2+A_1 & \Z & 1 \\
A_4 + A_1 & \Z & 1 \\
D_5 & \Z & 1 \\ 
A_3 +2A_1 & \Z\times \Z_2 & 2 \\ \hline
A_4 & \Z^2 & 2 \\
A_3+A_1 & \Z^2 & 2
\end{array}
\hspace*{0.5cm}
\begin{array}{l|l|l}
U & H_\Z/U & k_1\\ \hline
2A_2 & \Z^2 & 2 \\
A_2+2A_1 & \Z^2 & 2 \\
4A_1 & \Z^2\times \Z_2 & 3 \\
D_4 & \Z^2 & 2 \\ \hline
A_3 & \Z^3 & 3 \\
A_2+A_1 & \Z^3 & 3 \\
3A_1 & \Z^3 & 3 \\ \hline
A_2 & \Z^4 & 4 \\
2A_1 & \Z^4 & 4 \\ \hline
A_1 & \Z^5 & 5
\end{array}
\end{eqnarray*}
\caption[Table 3.2]{Isomorphism classes of pairs
$((H_\Z,I),U)$ for $(H_\Z,I)$ of type $E_6$}
\label{tab3.2}
\end{table}

(b) $(H_\Z,I)$ of types $E_6$, $E_7$ and $E_8$: 
The Tables \ref{tab3.3} and \ref{tab3.4} 
contain pairs $[H]'$ and $[H]''$ with
$H\in\{A_5+A_1,A_5,A_3+2A_1,A_3+A_1,4A_1,3A_1\}$ for $E_7$
and with $H\in\{A_7,A_5+A_1,2A_3,A_3+2A_1,4A_1\}$ for $E_8$.
Here $[H]'$ and $[H]''$ denote subroot lattices which are 
isomorphic if one forgets the embedding into $L$. 
But for a subroot lattice $U_1\subset H_\Z$ of type $[H]'$ and a 
subroot lattice $U_2\subset H_\Z$ of type $[H]''$, 
the pairs $((H_\Z,I),U_1)$ and $((H_\Z,I),U_2)$ are not 
isomorphic.

\begin{table}
\begin{eqnarray*}
\begin{array}{l|l|l}
U & H_\Z/U & k_1 \\ \hline
E_7 & \{0\} & 0 \\
D_6+A_1 & \Z_2 & 1 \\
A_5+A_2 & \Z_3 & 1 \\
2A_3+A_1 & \Z_4 & 1 \\
A_7 & \Z_2 & 1 \\
D_4+3A_1 & \Z^2_2 & 2 \\
7A_1 & \Z^3_2 & 3 \\ \hline
E_6 & \Z & 1 \\
D_5+A_1 & \Z & 1 \\
A_4+A_2 & \Z & 1 \\
A_3\! +\! A_2\! +\! A_1 & \Z & 1 \\
{}[A_5+A_1]' & \Z\times\Z_2 & 2 \\
{}[A_5+A_1]'' & \Z & 1 \\
D_6 & \Z & 1 \\
D_4+2A_1 & \Z\times\Z_2 & 2 \\
A_3+3A_1 & \Z\times\Z_2 & 2 \\
3A_2 & \Z\times\Z_3 & 2 \\
2A_3 & \Z\times\Z_2 & 2 \\
A_6 & \Z & 1 \\
6A_1 & \Z\times\Z^2_2 & 3 \\ \hline
D_5 & \Z^2 & 2 \\
A_4+A_1 & \Z^2 & 2 \\
2A_2+A_1 & \Z^2 & 2 \\\end{array}
\hspace*{0.5cm}
\begin{array}{l|l|l}
U & H_\Z/U & k_1 \\ \hline
{[}A_5]' & \Z^2 & 2 \\
{[}A_5]'' & \Z^2 & 2 \\
D_4+A_1 & \Z^2 & 2 \\
A_3+A_2 & \Z^2 & 2 \\
5A_1 & \Z^2\times\Z_2 & 3 \\
A_2+3A_1 & \Z^2 & 2 \\
{[}A_3+2A_1]' & \Z^2\times\Z_2 & 3 \\
{[}A_3+2A_1]'' & \Z^2 & 2 \\ \hline
D_4 & \Z^3 & 3 \\
A_4 & \Z^3 & 3 \\
{[}A_3+A_1]' & \Z^3 & 3 \\
{[}A_3+A_1]'' & \Z^3 & 3 \\
2A_2 & \Z^3 & 3 \\
A_2+2A_1 & \Z^3 & 3 \\
{[}4A_1]' & \Z^3\times\Z_2 & 4 \\
{[}4A_1]'' & \Z^3 & 3 \\ \hline
A_3 & \Z^4 & 4 \\
A_2+A_1 & \Z^4 & 4 \\
{[}3A_1]' & \Z^4 & 4 \\
{[}3A_1]'' & \Z^4 & 4 \\ \hline
A_2 & \Z^5 & 5 \\
2A_1 & \Z^5 & 5 \\ \hline
A_1 & \Z^6 & 6 
\end{array}
\end{eqnarray*}
\caption[Table 3.3]{Isomorphism classes of pairs
$((H_\Z,I),U)$ for $(H_\Z,I)$ of type $E_7$}
\label{tab3.3}
\end{table}

\begin{table}
\begin{eqnarray*}
\begin{array}{l|l|l}
U & H_\Z/U & k_1 \\ \hline
E_8 & \{0\} & 0 \\
A_8 & \Z_3 & 1 \\
D_8 & \Z_2 & 1 \\
A_7+A_1 & \Z_4 & 1 \\
A_5\! +\! A_2\! +\! A_1 & \Z_6 & 1 \\
2A_4 & \Z_5 & 1 \\
4A_2 & \Z_3^2 & 2 \\
E_6+A_2 & \Z_3 & 1 \\
E_7+A_1 & \Z_2 & 1 \\
D_6+2A_1 & \Z_2^2 & 2 \\
D_5+A_3 & \Z_4 & 1\\
2D_4 & \Z_2^2 & 2 \\
D_4+4A_1 & \Z_2^3 & 3 \\
2A_3\! +\! 2A_1 & \Z_2\times \Z_4 & 2 \\
8A_1 & \Z_2^4 & 4 \\ \hline
A_6+A_1 & \Z & 1\\
A_4\! +\! A_2\! +\! A_1 & \Z & 1\\
A_5+A_2 & \Z\times\Z_3 & 2\\
3A_2+A_1 & \Z\times\Z_3 & 2\\
E_6+A_1 & \Z & 1\\
E_7 & \Z & 1\\
D_7 & \Z & 1\\
D_5+2A_1 & \Z\times\Z_2 & 2\\
D_4+3A_1 & \Z\times\Z_2^2 & 3\\
2A_3+A_1 & \Z\times\Z_4 & 2\\
7A_1 & \Z\times\Z_2^3 & 4\\
D_6+A_1 & \Z\times\Z_2 & 2\\
D_5+A_2 & \Z & 1\\
A_3\! +\! A_2\! +\! 2A_1 & \Z\times\Z_2 & 2\\
D_4+A_3 & \Z\times\Z_2 & 2\\
A_3+4A_1 & \Z\times\Z_2^2 & 3\\
A_4+A_3 & \Z & 1\\
A_5+2A_1 & \Z\times\Z_2 & 2\\
{}[A_7]' & \Z & 1\\
{}[A_7]'' & \Z\times\Z_2 & 2 \\  \hline 
3A_2 & \Z^2\times\Z_3 & 3\\
E_6 & \Z^2 & 2\\
D_6 & \Z^2 & 2\\\end{array}
\hspace*{0.5cm}
\begin{array}{l|l|l}
U & H_\Z/U & k_1 \\ \hline
D_4+2A_1 & \Z^2\times\Z_2 & 3\\
{}[2A_3]' & \Z^2 & 2\\
{}[2A_3]'' & \Z^2\times\Z_2 & 3\\
D_5+A_1 & \Z^2 & 2\\
A_3+3A_1 & \Z^2\times\Z_2 & 3\\
D_4+A_2 & \Z^2 & 2\\
6A_1 & \Z^2\times\Z_2^2 & 4\\
A_2+4A_1 & \Z^2\times\Z_2 & 3\\
A_4+2A_1 & \Z^2 & 2\\
A_6 & \Z^2 & 2\\
A_3\! +\! A_2\! +\! A_1 & \Z^2 & 2\\
{}[A_5+A_1]' & \Z^2 & 2\\
{}[A_5+A_1]'' & \Z^2\times\Z_2 & 3\\
A_4+A_2 & \Z^2 & 2\\
2A_2+2A_1 & \Z^2 & 2 \\ \hline
D_5 & \Z^3 & 3\\
{}[A_3+2A_1]' & \Z^3 & 3\\
{}[A_3+2A_1]'' & \Z^3\times\Z_2 & 4\\
A_3+A_2 & \Z^3 & 3\\
A_5 & \Z^3 & 3\\
5A_1 & \Z^3\times\Z_2 & 4\\
A_4+A_1 & \Z^3 & 3\\
D_4+A_1 & \Z^3 & 3\\
A_2+3A_1 & \Z^3 & 3\\
2A_2+A_1 & \Z^3 & 3 \\ \hline
D_4 & \Z^4 & 4\\
{}[4A_1]' & \Z^4 & 4\\
{}[4A_1]'' & \Z^4\times\Z_2 & 5\\
A_2\! +\! 2A_1 & \Z^4 & 4\\
2A_2 & \Z^4 & 4\\
A_3+A_1 & \Z^4 & 4\\
A_4 & \Z^4 & 4\\   \hline
A_3 & \Z^5 & 5\\
A_2+A_1 & \Z^5 & 5\\
3A_1 & \Z^5 & 5 \\ \hline
A_2 & \Z^6 & 6\\
2A_1 & \Z^6 & 6\\ \hline
A_1 & \Z^7 & 7 
\end{array}
\end{eqnarray*}
\caption[Table 3.4]{Isomorphism classes of pairs
$((H_\Z,I),U)$ for $(H_\Z,I)$ of type $E_8$}
\label{tab3.4}
\end{table}

(c) With one class of exceptions, the following holds.
If $((H_\Z,I),U_1)$ and $((H_\Z,I),U_2)$ are
isomorphic pairs as in (b), then a Weyl group element $w\in W$
with $w(U_1)=U_2$ exists.
The class of exceptions are the sublattices of $D_n$ 
of types $\sum_{i=1}^{r_1}A_{c_i}$ with all $c_i$ odd.
For each of those types there are two conjugacy classes
with respect to $W$.
\end{theorem}

The following supplement to Theorem \ref{t3.6} will be useful
in the proof in Remark \ref{t3.15} (vii) of a part of Theorem
\ref{t3.14}.

\begin{corollary}\label{t3.7}
Let $(H_\Z,I)$ be an irreducible simply laced root system
of some rank $n\in\N$, 
and let $U\subset H_\Z$ be a subroot lattice with
index $i_1(H_\Z,U)>1$. Then the primitive sublattice
$U_\Q\cap H_\Z$ of $H_\Z$ is also a subroot lattice.
Especially 
\begin{eqnarray}\label{3.1}
U\cap\Phi\subsetneqq (U_\Q\cap H_\Z)\cap\Phi.
\end{eqnarray}
\end{corollary}

{\bf Proof:} One has to go through the construction in 
\cite{Dy57} or \cite{BaH19} of the subroot lattices
in the Tables \ref{tab3.1}--\ref{tab3.4}. It is the procedure
of Borel and de Siebenthal \cite{BdS49} and Dynkin \cite{Dy57}
which is described e.g. in \cite[Remark 3.2 (iv)]{BaH19}.
It consists in applying repeatedly one of two possible
steps (BDdS1) and (BDdS2). Step (BDdS1) leads from $H_\Z$
to a subroot lattice of the same rank, usually strictly smaller
than $H_\Z$. Step (BDdS2) leads from $H_\Z$ to a subroot lattice
which is a primitive sublattice of $H_\Z$ of rank $n-1$.

By \cite{BdS49} and \cite{Dy57} a suitable sequence of these
steps leads from $H_\Z$ to any subroot lattice $U$.
For the proof of the Corollary, one has to see that one can
choose the sequence in such a way that first the step
(BDdS2) is carried out $n-\rk U$ many times and then the step
(BDdS1) is carried out finitely many times,
because the intermediate subroot lattice $U_1$ which one 
obtains after carrying out (BDdS2) $n-\rk U$ many times is 
a primitive sublattice of rank equal to $\rk U$ 
and a subroot lattice
which contains $U$, so it is precisely $U_\Q\cap H_\Z$.

Unfortunately, the sequences of steps which are documented in
the Tables \ref{tab3.2}-\ref{tab3.4} in \cite{BaH19} for $E_6,E_7,E_8$ 
give first the steps of type (BDdS1) and then the steps
of type (BDdS2). But it is not difficult to see that one
can choose suitable sequences which give first the steps of type
(BDdS2) and then the steps of type (BDdS1).
The details are left to the reader.
\hfill$\Box$ 

\bigskip
Carter \cite{Ca72} studied and classified 
the conjugacy classes of the elements of the Weyl groups
of the irreducible simply laced root lattices. 
First we review some basic properties.

\begin{definition}\label{t3.8}
Let $(H_\Z,I)$ be a simply laced root lattice of rank
$n\in\N$ with root system $\Phi$ and Weyl group $W$. 

(a) (E.g. \cite{Bo68}) For any element $w\in W$ any tuple 
$(\alpha_1,...,\alpha_k)\in\Phi^k$ with $k\in\Z_{\geq 0}$ and
$$w=s_{\alpha_1}\circ ...\circ s_{\alpha_k}$$
is a {\it presentation} of $w$. Its {\it length} is $k$.
The length $l(w)\in\Z_{\geq 0}$ of $w$ is the minimum of
lengths of all presentations of $w$. The subroot lattice
of a presentation is 
$U:=\sum_{i=1}^k\Z\cdot\alpha_i\subset H_\Z$.
The {\it index} of this presentation is the index 
$i_1(H_\Z,U)$. 
A presentation of $w$ is called {\it reduced} if its length
is $l(w)$. 

(b) (E.g. \cite{Bo68}) 
An element $w\in W$ is a {\it Coxeter element} if it has a 
reduced presentation $(\alpha_1,...,\alpha_n)$ 
which is a root basis.

(c) \cite[Def. 3.2.1]{Vo85}
An element $w\in W$ is a {\it quasi-Coxeter element}
if $l(w)=n$ and if it has a reduced presentation with index $1$.
Remark: Obviously any Coxeter element is a quasi-Coxeter element.

(d) For any element $w\in W$ and any root of unity $\lambda$
\begin{eqnarray*}
H_\lambda(w)&:=& \ker (w-\lambda\cdot\id)\subset 
H_\C:=H_\Z\otimes_\Z\C,\\
H_{\neq 1}(w)&:=& \bigoplus_{\lambda\neq 1}H_\lambda(w)
\subset H_\C,
\end{eqnarray*}
and analogously $H_{\neq 1,\Q}(w),H_{\neq 1,\Z}(w),
H_{1,\Q}(w),H_{1,\Z}(w)$.
Observe $H_{\neq 1,\Q}(w)=H_{1,\Q}^\perp(w)$. 
\end{definition}

\begin{lemma}\label{t3.9}
Let $(H_\Z,I)$ be a simply laced root lattice
of rank $n\in\N$ with root system $\Phi$ and Weyl group $W$.
Let $w\in W$.

(a) (E.g. \cite{Bo68}) The Coxeter elements form one conjugacy
class with respect to $W$.

(b) \cite[Lemmata 1, 2 and 3]{Ca72}) 
A presentation $(\alpha_1,...,\alpha_k)\in\Phi^k$ of $w$
is reduced if and only if $\alpha_1,...,\alpha_k$
are linearly independent. The subroot lattice $U$ of a 
reduced presentation $(\alpha_1,...,\alpha_k)$ satisfies
\begin{eqnarray*}
\Bigl(\bigoplus_{i=1}^k\Q\cdot\alpha_i =\Bigr) U_\Q 
&=& H_{\neq 1,\Q}(w)\Bigl(= H_{1,\Q}^\perp(w)\Bigr),\\
\textup{and especially }l(w)&=& \dim_\Q H_{\neq 1,\Q}(w).
\end{eqnarray*}

(c) (\cite[Satz 3.2]{Kl83}\cite[Satz 3.2.3]{Vo85})
All reduced presentations of $w$ have the same index,
and this index is $|(\det(w-\id)|_{H_{\neq 1,\Q}}|$.

(d) Because of (b) and (c), all elements in the conjugacy class 
of a quasi-Coxeter element are quasi-Coxeter elements. 

(e) For any $\delta\in \Phi\cap H_{\neq 1,\Q}(w)$, the element
$w$ has a reduced presentation $(\alpha_1,...,\alpha_{l(w)})$
with $\alpha_1=\delta$. 
\end{lemma}

{\bf Proof of part (e):}
This is proved in \cite[Theorem 2.4.7]{Arm09}, but there it is 
also ascribed to \cite{Ca72}. It is implied by part (b) in the
following way. 

Consider any reduced presentation 
$(\beta_1,...,\beta_{l(w)})$ of $w$. Then 
$(\delta,\beta_1,...,\beta_{l(w)})$ is a nonreduced presentation
of $\www{w}:=s_{\delta}\circ w$ because $\delta,\beta_1,...,
\beta_{l(w)}$ are linearly dependent. This and
$\det\www{w}=(-1)^{l(\www{w})}=(-1)^{l(w)-1}$ imply
$l(\www{w})\leq l(w)-1$. 

Consider a reduced presentation 
$(\alpha_2,...,\alpha_{l(\www{w})+1})$ of $\www{w}$.
Then $(\delta,\alpha_2,...,\alpha_{l(\www{w})+1})$ is a 
presentation of $w$. Therefore $l(\www{w})=l(w)-1$, and
$(\delta,\alpha_2,...,\alpha_{l(w)})$ is a reduced 
presentation of $w$.\hfill$\Box$ 

\bigskip

The next theorem recalls results of Carter \cite{Ca72}
on the classification of the conjugacy classes of 
quasi-Coxeter elements 
in the case of an irreducible simply laced root lattice.

\begin{theorem}[{\cite{Ca72}, especially Table 3}]\label{t3.10}
Let $(H_\Z,I)$ be an irreducible simply laced root 
lattice of rank $n\in \N$.

Table \ref{tab3.5} gives in column 1 the type of the
root lattice, in column 2 Carter's symbols for the conjugacy
classes of all quasi-Coxeter elements $w$, 
in column 3 their characteristic polynomials,
and in column 4 their traces. 
The symbols $A_n,D_n,E_6,E_7,E_8$
denote the conjugacy classes of Coxeter elements.
\begin{table}
\begin{eqnarray*}
\begin{array}{l|l|l|l}
(H_\Z,I) & \textup{conj. class of }w & 
\textup{characteristic polynomial of }w & \tr(w)\\ \hline
A_n & A_n & (t^{n+1}-1)\cdot(t-1)^{-1} & -1\\ \hline
D_n & D_n & (t^{n-1}+1)(t+1)& -1\\
% & D_n(a_1) & (t^{n-2}+1)(t^2+1)& 0\\
% & \vdots & \vdots & \vdots \\
 & D_n(a_j) \textup{ for}& (t^{n-1-j}+1)(t^{j+1}+1) & 0\\
 & 1\leq j<\lfloor\frac{n}{2}\rfloor & & \\ \hline
E_6 & E_6 & (t^6+1)(t^3-1)\cdot (t^2+1)^{-1}(t-1)^{-1}& -1\\
 & E_6(a_1) & (t^9-1)\cdot (t^3-1)^{-1}& 0\\
 & E_6(a_2) & (t^6-1)(t^3+1)\cdot (t^2-1)^{-1}(t+1)^{-1}& 1\\ \hline
E_7 & E_7 & (t^9+1)(t+1)\cdot (t^3+1)^{-1}& -1\\
 & E_7(a_1) & t^7+1& 0\\
 & E_7(a_2) & (t^6+1)(t^3+1)\cdot (t^2+1)^{-1}& 0\\
 & E_7(a_3) & (t^5+1)(t^3+1)\cdot (t+1)^{-1}& 1\\
 & E_7(a_4) & (t^3+1)^3\cdot (t+1)^{-2} & 2\\ \hline
E_8 & E_8 & (t^{15}+1)(t+1)\cdot (t^5+1)^{-1}(t^3+1)^{-1}& -1\\
 & E_8(a_1) & (t^{12}+1)\cdot (t^4+1)^{-1}& 0\\
 & E_8(a_2) & (t^{10}+1)\cdot (t^2+1)^{-1}& 0\\
 & E_8(a_3) & (t^6+1)^2\cdot (t^2+1)^{-2}& 0\\
 & E_8(a_4) & (t^9+1)\cdot (t+1)^{-1}& 1\\
 & E_8(a_5) & (t^{15}-1)(t-1)\cdot (t^5-1)^{-1}(t^3-1)^{-1}& 1\\
 & E_8(a_6) & (t^5+1)^2\cdot (t+1)^{-2}& 2\\
 & E_8(a_7) & (t^6+1)(t^3+1)^2\cdot (t^2+1)^{-1}(t+1)^{-2}& 2\\
 & E_8(a_8) & (t^3+1)^4\cdot (t+1)^{-4}& 4
\end{array}
\end{eqnarray*}
\caption[Table 3.5]{The conjugacy classes of the quasi-Coxeter
elements in the irreducible simply laced root lattices}
\label{tab3.5}
\end{table}
\end{theorem}

\begin{remarks}\label{t3.11}
(i) Column 4 in Table \ref{tab3.5} follows from column 3
and the fact that the trace of an endomorphism with
characteristic polynomial $t^b\pm 1$ is $0$ of $b\geq 2$
and $\mp 1$ if $b=1$.

(ii) (E.g. \cite{Bo68}) In the case of $A_n$, the root lattice 
$H_\Z$ can be realized as follows,
\begin{eqnarray*}
H_\Z=\{\sum_{i=1}^{n+1}x_ie_i\,|\, \sum_{i=1}^{n+1}x_i=0\}
\subset\Z^{n+1}\textup{ with basis }e_1,...,e_{n+1}.
\end{eqnarray*}
Then $\Phi=\{\pm (e_i-e_j)\,|\, 1\leq i<j\leq n\}$.
The Weyl group is $S_{n+1}$ where $\sigma\in S_{n+1}$
acts on $H_\Z$ and $\Z^{n+1}$ by the $\Z$-linear extension
of the permutation $e_i\mapsto e_{\sigma(i)}$ of 
$e_1,...,e_{n+1}$.
The Coxeter elements are the cycles of length $n+1$.

(iii) (E.g. \cite{Bo68}) In the case of $D_n$ the root lattice
$H_\Z$ can be realized as follows,
\begin{eqnarray*}
H_\Z=\{\sum_{i=1}^nx_ie_i\,|\, \sum_{i=1}^nx_i\equiv 0\mmod 2\}
\subset\Z^n\textup{ with basis }e_1,...,e_n.
\end{eqnarray*}
Then $\Phi=\{\pm e_i\pm e_j)\,|\, 1\leq i<j\leq n\}$.

An element $(\varepsilon_1,...,\varepsilon_n,\sigma)$ 
of a semidirect product $\{\pm 1\}^n\rtimes S_n$ acts
on $\{\pm e_1,...,\pm e_n\}$ by mapping $\pm e_i$ to
$\pm \varepsilon_ie_{\sigma(i)}$. It is called a 
{\it signed permutation}. This action determines
the semidirect product $\{\pm 1\}^n\rtimes S_n$.
It extends $\Z$-linearly to an action on $\Z^n$
and on $H_\Z$. 

A signed permutation is called positive if 
$\prod_{i=1}^n\varepsilon_i=1$, negative if 
$\prod_{i=1}^n\varepsilon_i=-1$. 
The Weyl group $W$ is the index 2 subgroup of positive 
signed permutations.

A quasi-Coxeter element of type $D_n(a_j)$ for 
$j\in\{0,1,...,\lfloor\frac{n}{2}\rfloor-1\}$ with
$D_n(a_0):=D_n$ is a product of two negative signed
cycles with disjoint supports and lengths $j+1$ and
$n-(j+1)$. 
\end{remarks}

\begin{remarks}\label{t3.12}
(i) Consider an irreducible simply laced root lattice, 
an element $w\in W$, a reduced presentation
$(\alpha_1,...,\alpha_k)$ of $w$ and the corresponding
subroot lattice $U=\bigoplus_{i=1}^n\Z\cdot \alpha_i$.
Then $w$ restricts to a quasi-Coxeter element on this
subroot lattice $U$. The conjugacy class of $w$ is
determined by this quasi-Coxeter element on $U$.
One associates to this conjugacy class of $w$ a sum of symbols
of types $A_{b_1},D_{b_2}(a_j),E_{b_3}(a_j)$ with $j\geq 0$
as follows.
One identifies $U$ in the Tables \ref{tab3.1}--\ref{tab3.4}.
In the corresponding sum of symbols one replaces the symbols 
$D_{b_2},E_{b_3}$ by those symbols $D_{b_2}(a_j),E_{b_3}(a_j)$ 
with $D_b(a_0):=D_b$ and $E_b(a_0):=E_b$, which 
correspond to the correct quasi-Coxeter element on the 
corresponding irreducible subroot lattice. 

(ii) For any irreducible simply laced root lattice
and any element $w\in W$, let $U$ run through the
subroot lattices of all reduced presentations of $w$.
Theorem 5.10 in \cite{BaH19} implies the following.

In the cases of $A_n$ and $E_6$, 
the pairs $((H_\Z,I),U)$ are isomorphic,
so part (i) associates only one sum of symbols to $w$.

In the case of $D_n$, at least the quotients $H_\Z/U$
are isomorphic and thus the indices $i_1(H_\Z,U)$ 
and $k_1(H_\Z,U)$ are equal.

In the cases of $E_7$ and $E_8$ the indices $i_1(H_\Z,U)$
are the same by Lemma \ref{t3.9} (c),
but there are elements $w\in W$ where the pairs
$((H_\Z,I),U)$ are not isomorphic, and thus (i)
associates several sums of symbols to $w$.
The number of such conjugacy classes is $3$ in the case
of $E_7$ and $11$ in the case of $E_8$. Table 5.3 in
\cite{BaH19} lists them respectively their sums of symbols.

(iii) Table \ref{tab3.6} extracts from Table 5.3 in
\cite{BaH19} those conjugacy classes where also the quotients
$H_\Z/U$ are not all isomorphic and the numbers $k_1(H_\Z,U)$
are not all equal. These are $1$ of the $3$ in the case of 
$E_7$ and $4$ of the $11$ in the case of $E_8$. 
See Theorem \ref{t3.13} for column 3 in Table \ref{tab3.6}.
\begin{table}[H]
\begin{eqnarray*}
\begin{array}{l|l|l}
(H_\Z,I) & \textup{symbols in (i) for a }w, \textup{ their }k_1
& k_5(H_\Z,w) \\ \hline
E_7 & 2A_3+A_1\ (k_1=1)\sim D_4(a_1)+3A_1\ (k_1=2) & 1 \\
E_8 & 2A_3+A_1\ (k_1=2)\sim D_4(a_1)+3A_1\ (k_1=3) & 2 \\
E_8 & 2A_3+2A_1 (k_1=2)\sim D_4(a_1)+4A_1\ (k_1=3) & 2 \\
E_8 & D_5(a_1)+A_3\ (k_1=1)\sim D_4+D_4(a_1)\ (k_1=2) & 1 \\
E_8 & A_7+A_1\ (k_1=1)\sim D_5+A_3\ (k_1=1) & \\
 & \hspace*{2cm} \sim D_6(a_1)+2A_1\ (k_1=2) & 1 
\end{array}
\end{eqnarray*}
\caption[Table 3.6]{Conjugacy classes of elements of $W$ with
so different reduced presentations that even the numbers
$k_1(H_\Z,U)$ of the corresponding subroot lattices $U$
differ}
\label{tab3.6}
\end{table}
\end{remarks}

The length $k$ of any presentation $(\alpha_1,...,\alpha_k)$
of an element $w\in W$ for a simply laced root lattice differs
by an even number from $l(w)$ because $\det w=\pm 1$
and $\det s_\alpha=-1$. 

The following result on nonreduced presentations 
with subroot lattice the full lattice $H_\Z$ 
is a main result of \cite{BaH19}.
It will be crucial in the proofs of Theorem \ref{t7.2}
and Theorem \ref{t7.3} for the semidefinite cases.
The notations $k_1$ and $k_5$ follow 
\cite{BaH19}.

\begin{theorem}\label{t3.13}
Let $(H_\Z,I)$ be a simply laced root lattice of rank
$n\in\N$ with root system $\Phi$ and Weyl group $W$,
and let $w\in W$. The number
\begin{eqnarray*}
k_5(H_\Z,w)&:=&\min\left\{k\left|\begin{array}{l}
\textup{a presentation }(\alpha_1,...,\alpha_{l(w)+2k})
\textup{ with }\\
\textup{subroot lattice the full lattice exists}
\end{array}\right\}\right.
\end{eqnarray*}
satisfies
\begin{eqnarray*}
k_5(H_\Z,w)&=&k_1(H_\Z,U)\textup{ for }U
\textup{ \it the subroot lattice of any}\\
&&\hspace*{2cm}\textup{\it reduced presentation of }w
\end{eqnarray*}
in almost all cases. The only exceptions are listed in
Table \ref{t3.8}. In those few cases
\begin{eqnarray*}
k_5(H_\Z,w)&=&\min (k_1(H_\Z,U)\,|\, U\textup{ \it is a subroot
lattice of a}\\ 
&&\hspace*{3cm}\textup{\it reduced presentation of }w).
\end{eqnarray*}
\end{theorem}

Obviously, $k_5(H_\Z,w)$ is zero if and only if $w$ is a 
quasi-Coxeter element. 

The last theorem in this section, Theorem \ref{t3.14},
gives a complete answer to the question posed in the Remarks
\ref{t2.8} (i) 
for the simply laced root lattices. The Remarks
\ref{t3.15} comment on the history and part of the proof
and on related results.

\begin{theorem}[{Essentially \cite{De74,Vo85,BaGRW17}}]\label{t3.14}
Let $(H_\Z,I)$ be a simply laced root lattice
(irreducible or reducible). Let $w\in W$. The set of reduced
presentations of $w$,
\begin{eqnarray}\label{3.2}
\{(\alpha_1,...,\alpha_{l(w)})\in\Phi^{l(w)}\,|\, 
s_{\alpha_1}\circ...\circ s_{\alpha_{l(w)}}=w\}
\end{eqnarray}
is a single $\Br_{l(w)}\ltimes \{\pm 1\}^{l(w)}$ orbit if and
only if $w$ has index 1, so if and only if for some
(or any) reduced presentation $(\alpha_1,...,\alpha_{l(w)})$
of $w$ the subroot lattice $U:=\sum_{i=1}^{l(w)}\Z\cdot \alpha_i$
is a primitive sublattice of $H_\Z$, so $H_\Z/U$ has no torsion.
\end{theorem}

\begin{remarks}\label{t3.15}
(i) In the Tables \ref{tab3.1}--\ref{tab3.4}, in many cases
$H_\Z/U$ has no torsion and in many cases $H_\Z/U$ has torsion.

(ii) In the case of $(H_\Z,I)$ an irreducible simply laced
root lattice and $w$ a Coxeter element, Deligne \cite{De74}
proved the a priori slightly stronger statement that the set in 
\eqref{3.2} is a single $\Br_n$ orbit (and not only a single
$\Br_n\ltimes\{\pm 1\}^n$ orbit). 
But with Lemma \ref{t2.10} and Lemma \ref{t2.11} 
(a)(i) his result follows from Theorem \ref{t3.14}. 

(iii) Voigt \cite{Vo85} showed more generally for any
simply laced root lattice without summands of type $A_1$ 
and any quasi-Coxeter element $w$ that the set in \eqref{3.2}
is a single $\Br_n$ orbit. He had to forbid summands of type
$A_1$ because sign changes of their roots cannot be carried
out by braids. 

(iv) Though if one cares only about a single 
$\Br_n\ltimes\{\pm 1\}^n$ orbit (and not a $\Br_n$ orbit), one
can allow summands of type $A_1$. Voigt's result implies 
Theorem \ref{t3.14} in this case.

Because an element $w$ of index 1 is a quasi-Coxeter element
in the subroot lattice $U=U_\Q\cap H_\Z$, Voigt's result
implies the {\it if} implication in Theorem \ref{t3.14}.

(v) The {\it if} implication in Theorem \ref{t3.14} was
reproved by Baumeister, Gobet, Roberts and Wegener 
\cite{BaGRW17}. The {\it only if} implication was first 
observed and proved by them.

(vi) Though they did not use the condition
$$(C1)\quad w\textup{ has index 1,}$$
but the condition
$$(C2)\quad w\textup{ is a quasi-Coxeter element in a 
parabolic subgroup}.$$
Their definition of condition (C2) is slightly involved.
The equivalence $(C1)\iff (C2)$ is not obvious.

(vii) Instead of going through an explanation of $(C2)$ and 
a proof of this equivalence $(C1)\iff(C2)$, 
we give a direct proof of the {\it only if} implication of 
Theorem \ref{t3.14}, which essentially follows the argument 
in \cite{BaGRW17}. 

Consider $w\in W$ with index $>1$. It has to be shown that the
set in \eqref{3.2} is not a single
$\Br_{l(w)}\ltimes \{\pm 1\}^{l(w)}$ orbit. 
Consider a reduced presentation $(\alpha_1,...,\alpha_{l(w)})$
of $w$ and its subroot lattice 
$$U=\sum_{i=1}^{l(w)}\Z\cdot \alpha_i\subsetneqq U_\Q\cap H_\Z.$$
By \eqref{3.1} in Corollary \ref{t3.7} a root 
$\delta\in (U_\Q\cap H_\Z)\cap\Phi-U\cap\Phi$ exists.
By Lemma \ref{t3.9} (e) $w$ has reduced presentation
$(\beta_1,\beta_2,...,\beta_{l(w)})$ with $\beta_1=\delta$.
Reduced presentations in the $\Br_{l(w)}\ltimes\{\pm 1\}^{l(w)}$ 
orbit of $(\alpha_1,...,\alpha_{l(w)})$ have all the same
subroot lattice $U$. But $\delta\notin U$, so the 
reduced presentation $(\beta_1,...,\beta_{l(w)})$ is not in the
orbit of $(\alpha_1,...,\alpha_{l(w)})$. \hfill$\Box$
\end{remarks}

\section{Tubular elliptic root systems and their Coxeter elements}
\label{s4}
\setcounter{equation}{0}
\setcounter{table}{0}
\setcounter{figure}{0}

The positive semidefinite case in Theorem \ref{t1.5} is
Theorem \ref{t7.4}, which builds on the Theorems \ref{t7.3},
\ref{t4.9} and \ref{t4.11} and the Example \ref{t4.13}. 
The proof of Theorem \ref{t7.3} leads to the four 
{\it tubular elliptic root systems} whose definition in 
\cite{BaWY18} we recall in Definition \ref{t4.5}. 
Theorem \ref{t4.9} below gives a conceptual characterization
of their Coxeter elements. Together with Theorem \ref{t7.3}
it allows to restrict in the proof of Theorem \ref{t7.4}
to the four tubular elliptic root systems and their Coxeter 
elements. Theorem \ref{t4.11} below, which is
due to Kluitmann \cite{Kl87} and Baumeister, Wegener and 
Yahiatene \cite{BaWY18}, treats $\Br_n\ltimes\{\pm 1\}^n$ orbits
of $\Z$-bases for these cases and helps to understand the relevant 
$\Br_n\ltimes\{\pm 1\}^n$ orbits of matrices in $T_n(\Z)$.

In section \ref{s5} the classical relation of three of the four
tubular elliptic root lattices and their Coxeter elements
to the simple elliptic singularities is recalled. 

We start with a more general definition, which is due to
K. Saito \cite[(2.1) Definition 1]{SaK85}. 
After the classical Theorem \ref{t4.3}, which gives some 
structural results for these more general objects, we come
to the tubular elliptic root systems.

\begin{definition}\label{t4.1}
An {\it irreducible simply laced elliptic root system}
consists of the following data $(H_\Z,I,\Phi)$.

$H_\Z$ is a $\Z$-lattice of some rank $n\in\Z_{\geq 2}$ 
with a symmetric positive semidefinite bilinear form 
$I:H_\Z\times H_\Z\to\Z$ with radical $\Rad(I)\subset H_\Z$ 
of rank 2.
$\Phi\subset \{\delta\in H_\Z\,|\, I(\delta,\delta)=2\}$
is a set of roots, and 
$W:=\langle s_v\,|\, v\in\Phi\rangle\subset\Aut(H_\Z,I)$
is the corresponding Weyl group. They have to satisfy 
\begin{eqnarray*}
\sum_{v\in \Phi}\Z\cdot v&=&H_\Z\quad\textup{and}\\ 
g(\Phi)&=&\Phi\quad\textup{for }g\in W,\\
\Phi=\Phi_1\, \dot\cup\, \Phi_2&\textup{with}&
\Phi_1\perp_I\Phi_2\ \Rightarrow\ \Phi_1=\emptyset\textup{ or }
\Phi_2=\emptyset.
\end{eqnarray*}
\end{definition}

\begin{remarks}\label{t4.2}
(i) \cite{SaK85} defined more generally {\it irreducible 
elliptic root systems}. Though in \cite{SaK85} they are
called {\it extended affine root systems}. In \cite{SaT97}
K. Saito changed the name to {\it elliptic root systems}.

We follow the name {\sf irreducible} {\it elliptic root systems}
in \cite[Definition 5.1]{BaWY18}. Beware that the 
{\it extended root systems} in section 2 in \cite{BaWY18}
are quite different objects.
Beware also that \cite{AABDP97} uses the name 
{\it extended affine root systems} for a generalization of
the notion in \cite{SaK85}, where the radical $\Rad I$ can
have arbitrary rank $>0$. 

(ii) One could live without the irreducibility condition
at the end of Definition \ref{t4.1}. Therefore \cite{BaWY18}
put {\it irreducible} into the name {\it elliptic root system}.

(iii) Denote by $\oooo{H_\Z}:=H_\Z/\Rad(I)$ the quotient
lattice. For any structure/object on/in $H_\Z$ denote by
an overline the induced structure/object on/in $\oooo{H_\Z}$.
The pairing induces a symmetric positive definite pairing
$\oooo{I}:\oooo{H_\Z}\times \oooo{H_\Z}\to\Z$.
Obviously $(\oooo{H_\Z},\oooo{I})$ is a simply laced
root lattice, $\oooo{\Phi}$ is its set of roots,
and $\oooo{W}$ is its Weyl group. The simply laced
root lattice $(\oooo{H_\Z},\oooo{I})$ is irreducible because
$(H_\Z,I,\Phi)$ is irreducible.
\end{remarks}

The structure of irreducible simply laced elliptic root systems
is not difficult for $n\geq 4$. The next theorem follows
from the classification in \cite{SaK85}. It is also proved
explicitly in \cite[Theorem 2.37]{AABDP97} 
(in much greater generality, for irreducible simply laced
extended affine root systems with $\Rad(I)$ of arbitrary rank
$>0$) and in \cite[Proposition 5.4]{BaWY18}.

\begin{theorem}\label{t4.3}
Let $(H_\Z,I,\Phi)$ be an irreducible simply laced elliptic
root system of rank $n\geq 4$. Then there is an isomorphism
of lattices with bilinear forms, which extends also as follows
to the set of roots,
\begin{eqnarray*}
(H_\Z,I)&\cong & (\oooo{H_\Z},\oooo{I})\oplus (\Rad(I),0),\\
\Phi &\cong & \oooo{\Phi}+\Rad(I).
\end{eqnarray*}
Especially $\Phi=\{\delta\in H_\Z\,|\, I(\delta,\delta)=2\}$.
\end{theorem}

\begin{remarks}\label{t4.4}
(i) In rank $n=3$ this is not true. There are two 
irreducible simply laced elliptic root systems of rank 3.
They are called $A_1^{(1,1)}$ and $A_1^{(1,1)*}$ in \cite{SaK85}.
The second one does not fit into Theorem \ref{t4.3}.
See also \cite[Ch. 2, Prop 4.2 and Table 4.5]{AABDP97} 
for descriptions of them.

(ii) The last statement follows from 
$\Phi\cong\oooo{\Phi}+\Rad(I)$ and 
$\oooo{\Phi}=\{\delta\in\oooo{H_\Z}\,|\, 
\oooo{I}(\delta,\delta)=2\}$. 

(iii) The classification of not simply laced extended affine
root systems (with $\Rad(I)$ of arbitrary rank $>0$) 
in \cite[Ch. 2]{AABDP97} is more complicated.

(iv) By Theorem \ref{t4.3}, an irreducible simply laced
elliptic root system of rank $n\geq 4$ is characterized 
uniquely (up to isomorphism) by the underlying irreducible 
simply laced root lattice $(\oooo{H_\Z},\oooo{I})$.
\end{remarks}

\begin{definition}[{\cite[Definition 5.7]{BaWY18}}]\label{t4.5}
A {\it tubular elliptic root system} is an irreducible simply laced
elliptic root system $(H_\Z,I,\Phi)$ with quotient lattice
$(\oooo{H_\Z},\oooo{I})$ of type $D_4,E_6,E_7$ or $E_8$. 
For each of the four types, there is only one such system,
by Theorem \ref{t4.3}. It is called 
$D_4^{(1,1)},E_6^{(1,1)},E_7^{(1,1)}$ respectively
$E_8^{(1,1)}$. 
\end{definition}

Definition \ref{t4.7} of a {\it Coxeter element} 
in a tubular elliptic root system uses a 
priori a $\Z$-basis of $H_\Z$ with a specific Cartan matrix 
respectively associated Dynkin diagram.
Theorem \ref{t4.9} below, which is the main result of this 
section, will improve this and give a more conceptual 
characterization of the Coxeter elements.

The classification of the irreducible simply laced elliptic
root systems in \cite{SaK85} uses a special notion of 
Dynkin diagram, which we will not recall here.
In the case of the tubular elliptic root systems, the classification
in \cite{SaK85} boils down to the following.
See also \cite[Example 5.12]{BaWY18}.

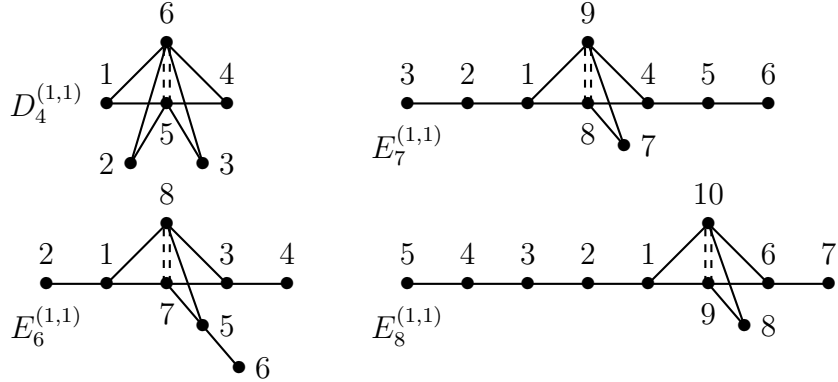
\begin{figure}
\begin{tikzpicture}[scale=0.8]
\node at (0,5) {$D_4^{(1,1)}$};
\node at (1,5) {$\bullet$};
\node at (2,5) {$\bullet$};
\node at (3,5) {$\bullet$};
\node at (2,6) {$\bullet$};
\draw [thick] (1,5)--(3,5);
\draw [thick] (1,5)--(2,6);
\draw [thick] (3,5)--(2,6);
\draw [thick] (1.4,4)--(2,5);
\draw [thick] (1.4,4)--(2,6);
\draw [thick] (2.6,4)--(2,5);
\draw [thick] (2.6,4)--(2,6);
\draw [thick,dashed] (1.95,5)--(1.95,6);
\draw [thick,dashed] (2.05,5)--(2.05,6);
\node at (1.4,4) {$\bullet$};
\node at (2.6,4) {$\bullet$};
\node at (1,5.5) {1};
\node at (3,5.5) {4};
\node at (1,4) {2};
\node at (3,4) {3};
\node at (2,4.5) {5};
\node at (2,6.5) {6};
\node at (0,1.3) {$E_6^{(1,1)}$};
\draw [thick] (0,2)--(4,2);
\node at (0,2) {$\bullet$};
\node at (1,2) {$\bullet$};
\node at (2,2) {$\bullet$};
\node at (3,2) {$\bullet$};
\node at (4,2) {$\bullet$};
\node at (2.6,1.3) {$\bullet$};
\draw [thick] (2,2)--(2.6,1.3);
\node at (3.2,0.6) {$\bullet$};
\draw [thick] (2.6,1.3)--(3.2,0.6);
\node at (2,3) {$\bullet$};
\draw [thick] (1,2)--(2,3);
\draw [thick] (3,2)--(2,3);
\draw [thick] (2.6,1.3)--(2,3);
\draw [thick,dashed] (1.95,2)--(1.95,3);
\draw [thick,dashed] (2.05,2)--(2.05,3);
\node at (0,2.5) {2};
\node at (1,2.5) {1};
\node at (3,2.5) {3};
\node at (4,2.5) {4};
\node at (2,1.5) {7};
\node at (3,1.3) {5};
\node at (3.6,0.6) {6};
\node at (2,3.5) {8};
\node at (6,4.3) {$E_7^{(1,1)}$};
\draw [thick] (6,5)--(12,5);
\node at (6,5) {$\bullet$};
\node at (7,5) {$\bullet$};
\node at (8,5) {$\bullet$};
\node at (9,5) {$\bullet$};
\node at (10,5) {$\bullet$};
\node at (11,5) {$\bullet$};
\node at (12,5) {$\bullet$};
\node at (9.6,4.3) {$\bullet$};
\draw [thick] (9,5)--(9.6,4.3);
\node at (9,6) {$\bullet$};
\draw [thick] (8,5)--(9,6);
\draw [thick] (10,5)--(9,6);
\draw [thick] (9.6,4.3)--(9,6);
\draw [thick,dashed] (8.95,5)--(8.95,6);
\draw [thick,dashed] (9.05,5)--(9.05,6);
\node at (6,5.5) {3};
\node at (7,5.5) {2};
\node at (8,5.5) {1};
\node at (10,5.5) {4};
\node at (11,5.5) {5};
\node at (12,5.5) {6};
\node at (10,4.3) {7};
\node at (9,4.5) {8};
\node at (9,6.5) {9};
\node at (6,1.3) {$E_8^{(1,1)}$};
\draw [thick] (6,2)--(13,2);
\node at (6,2) {$\bullet$};
\node at (7,2) {$\bullet$};
\node at (8,2) {$\bullet$};
\node at (9,2) {$\bullet$};
\node at (10,2) {$\bullet$};
\node at (11,2) {$\bullet$};
\node at (12,2) {$\bullet$};
\node at (13,2) {$\bullet$};
\node at (11.6,1.3) {$\bullet$};
\draw [thick] (11,2)--(11.6,1.3);
\node at (11,3) {$\bullet$};
\draw [thick] (10,2)--(11,3);
\draw [thick] (12,2)--(11,3);
\draw [thick] (11.6,1.3)--(11,3);
\draw [thick,dashed] (10.95,2)--(10.95,3);
\draw [thick,dashed] (11.05,2)--(11.05,3);
\node at (6,2.5) {5};
\node at (7,2.5) {4};
\node at (8,2.5) {3};
\node at (9,2.5) {2};
\node at (10,2.5) {1};
\node at (12,2.5) {6};
\node at (13,2.5) {7};
\node at (12,1.3) {8};
\node at (11,1.5) {9};
\node at (11,3.5) {10};
\end{tikzpicture}
\caption[Figure 4.1]{Dynkin diagrams of the tubular elliptic root
systems}
\label{fig4.1}
\end{figure}

\begin{theorem}\label{t4.6}
Let $(H_\Z,I,\Phi)$ be a tubular elliptic root system.
It has a $\Z$-basis $\uuuu{\delta}=(\delta_1,...,\delta_n)$
whose Cartan matrix $I(\uuuu{\delta}^t,\uuuu{\delta})$
has the Dynkin diagram in Figure \ref{fig4.1}.
Here the Dynkin diagram is defined as ${\rm CDD}(S)$ (see Remark 
\ref{t2.2} (vi)) where we write
$I(\uuuu{\delta}^t,\uuuu{\delta})=S+S^t$ with
$S\in T_n(\Z)$. Automatically 
$\delta_{n-1}-\delta_n\in\Rad(I)$, so 
$\oooo{\delta_{n-1}}=\oooo{\delta_n}$.
\end{theorem}

\begin{definition}\label{t4.7}
(\cite[(9.7) Definition]{SaK85} and 
\cite[Definition 2.4 (f)]{BaWY18})
Let $(H_\Z,I,\Phi)$ be a tubular elliptic root system.
An element $c\in W$ is a Coxeter element if there
is a $\Z$-basis $\uuuu{\delta}$ with Dynkin diagram of
its Cartan matrix as in Figure \ref{fig4.1} and
$c=s_{\delta_{\sigma(1)}}\circ ...\circ s_{\delta_{\sigma(n-2)}}
\circ s_{\delta_{n-1}}\circ s_{\delta_n}$ with 
$\sigma\in S_{n-2}$. 
\end{definition}

\begin{remarks}\label{t4.8}
(i) This Definition using a Dynkin diagram is somewhat 
analogous to the definition of a Coxeter 
element in a simply laced root lattice via a root basis. 
Though Theorem \ref{t4.9} below offers a more conceptual 
characterization of the Coxeter elements 
in a tubular elliptic root system. 

(ii) In fact, the definition in \cite[(9.7) Definition]{SaK85}
allows an almost arbitrary ordering of the reflections 
$s_{\delta_1},...,s_{\delta_n}$ in a product giving a 
Coxeter element. It only demands that 
$s_{\delta_{n-1}}$ is immediately followed by $s_{\delta_n}$.
Definition \ref{t4.7} is taken from 
\cite[Definition 2.4 (f)]{BaWY18}.

(iii) Lemma 1 in \cite[Ch. V, \S 6, Lemma 1]{Bo68} and 
the shape of the Dynkin diagrams in Figure \ref{fig4.1}
imply that the Coxeter
elements with respect to the same $\Z$-basis $\uuuu{\delta}$ 
are conjugate in $W$. A direct proof of this fact is offered in 
\cite[Lemma 4.1]{BaWY21}. 

(iv) All Coxeter elements are conjugate with respect to
$\Aut(H_\Z,I)$ \cite[Proposition 5.21 (c)]{BaWY18}.
\end{remarks}

Theorem \ref{t4.9} will be used in the proof of Theorem 
\ref{t7.4}.

\begin{theorem}\label{t4.9}
Let $(H_\Z,I,\Phi)$ be a tubular elliptic root system.
Let $(e_1,...,e_n)\in \Phi^n$ be a $\Z$-basis of $H_\Z$ and
$c:=s_{e_1}\circ...\circ s_{e_n}\in W$
such that $\oooo{c}\in\oooo{W}$ has a reduced presentation
with associated subroot lattice of the type in the third line
in Table \ref{tab4.1} if $(H_\Z,I,\Phi)$ is of the type in
the first line in Table \ref{tab4.1}.
\begin{table}[H]
\begin{eqnarray*}
\begin{array}{l|l|l|l}
D_4^{(1,1)} & E_6^{(1,1)} &E_7^{(1,1)} & E_8^{(1,1)}\\ \hline
D_4 & E_6 &E_7 & E_8\\ \hline
2D_2 & 3A_2 & 2A_3+A_1 & A_5+A_2+A_1
\end{array}
\end{eqnarray*}
\caption[Table 4.1]{A Coxeter element in a tubular elliptic root system
$H_\Z$ is determined by the induced element in $\oooo{\Phi}$}
\label{tab4.1}
\end{table}
Then $c$ is a Coxeter element.
\end{theorem}

This theorem is an almost immediate consequence of 
Theorem \ref{t4.10}, which states some facts for the
simply laced root lattices of type $D_4,E_6,E_7$ and $E_8$. 
After Theorem \ref{t4.10} we will first explain how it implies 
Theorem \ref{t4.9} and then prove Theorem \ref{t4.10}. 
The proof of Theorem \ref{t4.10} requires little own work. 
It builds on work done in \cite{LR16}, \cite{WY20} and 
\cite{BaWY21}.

\begin{theorem}\label{t4.10}
Let $(H_\Z,I)$ be a root lattice of a type in 
the second line in Table \ref{tab4.1}.
Let $c$ be an element of the Weyl group $W$ of
$(H_\Z,I)$ with a reduced presentation with
associated subroot lattice of the type in the third line in
Table \ref{tab4.1}. Let $(e_1,...,e_{n+2})$ be a
nonreduced presentation of $c$ 
(so $c=s_{e_1}\circ...\circ s_{e_{n+2}}$)
with $\sum_{i=1}^{n+2}\Z\cdot e_i=H_\Z$. 

Then there are a tuple 
$\uuuu{\delta}=(\delta_1,...,\delta_{n+2})\in\Phi^{n+2}$
with $\delta_{n+1}=\delta_{n+2}$ and a permutation
$\sigma\in S_n$ such that the $\Br_{n+2}\ltimes\{\pm 1\}^{n+2}$
orbit of $(e_1,...,e_{n+2})$ contains
$(\delta_{\sigma(1)},...,\delta_{\sigma(n)},\delta_{n+1},
\delta_{n+1})$ and such that the Dynkin diagram
of the matrix $I(\uuuu{\delta}^t,\uuuu{\delta})$ is as
in Figure \ref{fig4.1} (for the case in the first line
in Table \ref{tab4.1}). 
\end{theorem}

{\bf Proof of Theorem \ref{t4.9} using Theorem \ref{t4.10}:}
Start with the situation in Theorem \ref{t4.9}.
The quotient lattice $\oooo{H_\Z}:= H_\Z/\Rad(I)$ has rank
$\oooo{n}=n-2$. The induced data in $\oooo{H_\Z}$ are denoted
by overlining data in $H_\Z$. 
Theorem \ref{t4.10} applies to the data 
$(\oooo{H_\Z},\oooo{I})$ and $(\oooo{e_1},...,\oooo{e_n})$.
There are an element of $\Br_n\ltimes\{\pm 1\}^n$,
a tuple $\uuuu{\www{\delta}}=(\www{\delta_1},...,\www{\delta_n})$
with $\www{\delta_{n-1}}=\www{\delta_n}$ 
and a permutation $\sigma\in S_{n-2}$ such that
the element of $\Br_n\ltimes\{\pm 1\}^n$ maps 
$(\oooo{e_1},...,\oooo{e_n})$ to the tuple
$(\www{\delta_{\sigma(1)}},...,\www{\delta_{\sigma(n-2)}},
\www{\delta_{n-1}},\www{\delta_{n-1}})$ 
and such that the Dynkin
diagram of the matrix $\oooo{I}(\uuuu{\www{\delta}}^t,
\uuuu{\www{\delta}})$ is as in Figure \ref{fig4.1}.
The same element of $\Br_n\ltimes\{\pm 1\}^n$ 
maps $(e_1,...,e_n)$ to a $\Z$-basis
$(\delta_{\sigma(1)},...,\delta_{\sigma(n-2)},\delta_{n-1},
\delta_n)$ with $\oooo{\delta_i}=\www{\delta_i}$
and such that the Dynkin diagram of the matrix 
$(I(\uuuu{\delta}^t,\uuuu{\delta}))$ 
with $\uuuu{\delta}=(\delta_1,...,\delta_n)$ 
is as in Figure \ref{fig4.1}.
%Especially $\oooo{\delta_{n-1}}=\www{\delta_{n-1}}
%=\www{\delta_n}=\oooo{\delta_n}$. 
Therefore 
$$c=s_{e_1}\circ ...\circ s_{e_n}
=s_{\delta_{\sigma(1)}}\circ ...\circ s_{\delta_{\sigma(n-2)}}
\circ s_{\delta_{n-1}}\circ s_{\delta_n}$$
is a Coxeter element by Definition \ref{t4.7}.
\hfill$\Box$ 

\bigskip
{\bf Proof of Theorem \ref{t4.10}:}
The element $c$ has the maximal length $l(c)=n$, 
so the presentation by $(e_1,...,e_{n+2})$ of length $n+2$ 
is nonreduced.
Corollary 1.4 in \cite{LR16} (a second proof is given in
\cite[Theorem 1.1]{WY23}) implies that the Hurwitz orbit
by $\Br_{n+2}$ of the tuple $(s_{e_1},...,s_{e_{n+2}})$ contains
a tuple $(s_{\alpha_1},...,s_{\alpha_{n+2}})$
with $\alpha_{n+1}=\alpha_{n+2}$.
Therefore the $\Br_{n+2}\ltimes\{\pm 1\}^{n+2}$ orbit of
$(e_1,..,e_{n+2})$ contains 
$(\alpha_1,...,\alpha_{n+2})$.
Then $H_\Z=\sum_{i=1}^{n+2}\Z\cdot e_i
=\sum_{i=1}^{n+1}\Z\cdot\alpha_i$. 

Obviously $c=s_{\alpha_1}\circ...\circ s_{\alpha_n}$,
so the tuple $(\alpha_1,...,\alpha_n)$ is a reduced
presentation of $c$. We claim that the associated subroot lattice
is of the type in the third line in Table \ref{tab4.1}. 
In all cases except $E_7$ this is clear because
the type of the subroot lattice is unique by Table 
5.3 in \cite{BaH19}. In the case $E_7$
there are two possible subroot lattices $2A_3+A_1$ and 
$D_4(a_1)+3A_1$. If the associated subroot lattice
$U=\sum_{i=1}^7\Z\cdot\alpha_i$ were of type 
$D_4(a_1)+3A_1$ then $H_\Z/U\cong\Z_2^2$, but this does not allow 
$H_\Z=\sum_{i=1}^8\Z\cdot\alpha_i=U+\Z\cdot\alpha_8$.
Therefore $U$ is of type $2A_3+A_1$ in the case $E_7$. 

The construction of $U$ from $H_\Z$ by the procedure of 
Borel and de Siebenthal \cite{BdS49} and Dynkin \cite{Dy57}
requires a single step of type (BDdS1). 
One goes from the Dynkin diagram of type $D_4,E_6,E_7$ 
respectively $E_8$ to its extended Dynkin diagram in Figure
\ref{fig4.2} and then omits the center of this Dynkin diagram.

\begin{figure}
\begin{tikzpicture}[scale=0.8]
\node at (0,5) {$D_4$};
\node at (1,5) {$\bullet$};
\node at (2,5) {$\bullet$};
\node at (3,5) {$\bullet$};
\draw [thick] (1,5)--(3,5);
\draw [thick] (1.4,4)--(2,5);
\draw [thick] (2.6,4)--(2,5);
\node at (1.4,4) {$\bullet$};
\node at (2.6,4) {$\bullet$};
\node at (1,5.5) {1};
\node at (3,5.5) {4};
\node at (1,4) {2};
\node at (3,4) {3};
\node at (2,5.5) {5};
\node at (0,1.3) {$E_6$};
\draw [thick] (0,2)--(4,2);
\node at (0,2) {$\bullet$};
\node at (1,2) {$\bullet$};
\node at (2,2) {$\bullet$};
\node at (3,2) {$\bullet$};
\node at (4,2) {$\bullet$};
\node at (2.6,1.3) {$\bullet$};
\draw [thick] (2,2)--(2.6,1.3);
\node at (3.2,0.6) {$\bullet$};
\draw [thick] (2.6,1.3)--(3.2,0.6);
\node at (0,2.5) {2};
\node at (1,2.5) {1};
\node at (3,2.5) {3};
\node at (4,2.5) {4};
\node at (2,2.5) {7};
\node at (3,1.3) {5};
\node at (3.6,0.6) {6};
\node at (6,4.3) {$E_7$};
\draw [thick] (6,5)--(12,5);
\node at (6,5) {$\bullet$};
\node at (7,5) {$\bullet$};
\node at (8,5) {$\bullet$};
\node at (9,5) {$\bullet$};
\node at (10,5) {$\bullet$};
\node at (11,5) {$\bullet$};
\node at (12,5) {$\bullet$};
\node at (9.6,4.3) {$\bullet$};
\draw [thick] (9,5)--(9.6,4.3);
\node at (6,5.5) {3};
\node at (7,5.5) {2};
\node at (8,5.5) {1};
\node at (10,5.5) {4};
\node at (11,5.5) {5};
\node at (12,5.5) {6};
\node at (10,4.3) {7};
\node at (9,5.5) {8};
\node at (6,1.3) {$E_8$};
\draw [thick] (6,2)--(13,2);
\node at (6,2) {$\bullet$};
\node at (7,2) {$\bullet$};
\node at (8,2) {$\bullet$};
\node at (9,2) {$\bullet$};
\node at (10,2) {$\bullet$};
\node at (11,2) {$\bullet$};
\node at (12,2) {$\bullet$};
\node at (13,2) {$\bullet$};
\node at (11.6,1.3) {$\bullet$};
\draw [thick] (11,2)--(11.6,1.3);
\node at (6,2.5) {5};
\node at (7,2.5) {4};
\node at (8,2.5) {3};
\node at (9,2.5) {2};
\node at (10,2.5) {1};
\node at (12,2.5) {6};
\node at (13,2.5) {7};
\node at (12,1.3) {8};
\node at (11,2.5) {9};
\end{tikzpicture}
\caption[Figure 4.2]{Extended Dynkin diagrams}
\label{fig4.2}
\end{figure}
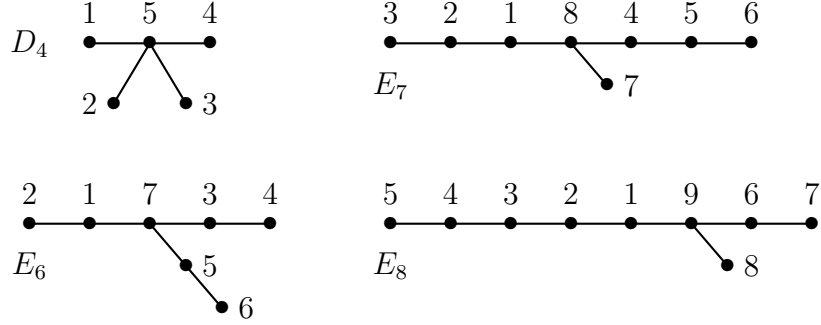

\begin{figure}
\begin{tikzpicture}[scale=0.8]
\node at (0,5) {$D_4$};
\node at (1,5) {$\bullet$};
\node at (2,5) {$\bullet$};
\node at (3,5) {$\bullet$};
\draw [thick] (1,5)--(3,5);
\draw [thick] (1.4,4)--(2,5);
\draw [thick] (2.6,4)--(2,5);
\node at (1.4,4) {$\bullet$};
\node at (2.6,4) {$\bullet$};
\node at (1,5.5) {1};
\node at (3,5.5) {1};
\node at (1,4) {1};
\node at (3,4) {1};
\node at (2,5.5) {2};
\node at (0,1.3) {$E_6$};
\draw [thick] (0,2)--(4,2);
\node at (0,2) {$\bullet$};
\node at (1,2) {$\bullet$};
\node at (2,2) {$\bullet$};
\node at (3,2) {$\bullet$};
\node at (4,2) {$\bullet$};
\node at (2.6,1.3) {$\bullet$};
\draw [thick] (2,2)--(2.6,1.3);
\node at (3.2,0.6) {$\bullet$};
\draw [thick] (2.6,1.3)--(3.2,0.6);
\node at (0,2.5) {1};
\node at (1,2.5) {2};
\node at (3,2.5) {2};
\node at (4,2.5) {1};
\node at (2,2.5) {3};
\node at (3,1.3) {2};
\node at (3.6,0.6) {1};
\node at (6,4.3) {$E_7$};
\draw [thick] (6,5)--(12,5);
\node at (6,5) {$\bullet$};
\node at (7,5) {$\bullet$};
\node at (8,5) {$\bullet$};
\node at (9,5) {$\bullet$};
\node at (10,5) {$\bullet$};
\node at (11,5) {$\bullet$};
\node at (12,5) {$\bullet$};
\node at (9.6,4.3) {$\bullet$};
\draw [thick] (9,5)--(9.6,4.3);
\node at (6,5.5) {1};
\node at (7,5.5) {2};
\node at (8,5.5) {3};
\node at (10,5.5) {3};
\node at (11,5.5) {2};
\node at (12,5.5) {1};
\node at (10,4.3) {2};
\node at (9,5.5) {4};
\node at (6,1.3) {$E_8$};
\draw [thick] (6,2)--(13,2);
\node at (6,2) {$\bullet$};
\node at (7,2) {$\bullet$};
\node at (8,2) {$\bullet$};
\node at (9,2) {$\bullet$};
\node at (10,2) {$\bullet$};
\node at (11,2) {$\bullet$};
\node at (12,2) {$\bullet$};
\node at (13,2) {$\bullet$};
\node at (11.6,1.3) {$\bullet$};
\draw [thick] (11,2)--(11.6,1.3);
\node at (6,2.5) {1};
\node at (7,2.5) {2};
\node at (8,2.5) {3};
\node at (9,2.5) {4};
\node at (10,2.5) {5};
\node at (12,2.5) {4};
\node at (13,2.5) {2};
\node at (12,1.3) {3};
\node at (11,2.5) {6};
\end{tikzpicture}
\caption[Figure 4.3]{Coefficients of the generating relation}
\label{fig4.3}
\end{figure}
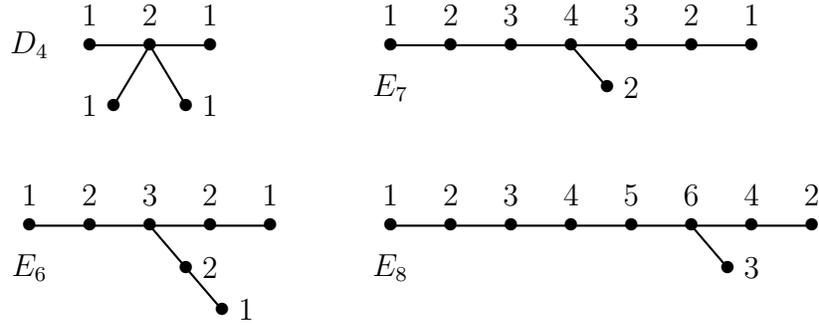

On the level of $\Z$-lattices, one has a generating tuple 
$\uuuu{\delta}=(\delta_1,...,\delta_{n+1})$ of $H_\Z$ as
$\Z$-lattice with one generating relation
$0=\sum_{i=1}^{n+1}b_i\delta_i$ with coefficients
$(b_1,...,b_{n+1})\in\Z^{n+1}$ in Figure \ref{fig4.3}
and such that the Dynkin diagram of the matrix
$I(\uuuu{\delta}^t,\uuuu{\delta})$ is the extended Dynkin 
diagram. Then 
\begin{eqnarray*}
U=\sum_{i=1}^n\Z\cdot\delta_i
\subset H_\Z=U+\Z\cdot\delta_{n+1},\\
\begin{array}{l|l|l|l|l}
H_\Z & D_4 & E_6 & E_7 & E_8 \\ \hline 
{}[H_\Z:U]=b_{n+1} & 2 & 3 & 4 & 6 
\end{array}
\end{eqnarray*} 
Because all irreducible components of $U$ are of type
$A_\bullet$, $c$ is a Coxeter element of the subroot lattice
$U$, and it has a reduced presentation 
$(\delta_{\sigma(1)},...,\delta_{\sigma(n)})$ for some
$\sigma\in S_n$ (e.g. \cite{Bo68}). 

Theorem \ref{t3.14} applies and says that the reduced 
presentations $(\alpha_1,...,\alpha_n)$ and
$(\delta_{\sigma(1)},...,\delta_{\sigma(n)})$ are in one
$\Br_n\ltimes\{\pm 1\}^n$ orbit. 

Therefore $(\alpha_1,....,\alpha_n,\alpha_{n+1},\alpha_{n+1})$
and $(\delta_{\sigma(1)},...,\delta_{\sigma(n)},
\alpha_{n+1},\alpha_{n+1})$ are in one 
$\Br_{n+2}\ltimes\{\pm 1\}^{n+2}$ orbit. 

It remains to see that 
$(\delta_{\sigma(1)},...,\delta_{\sigma(n)},
\alpha_{n+1},\alpha_{n+1})$ and
$(\delta_{\sigma(1)},...,\delta_{\sigma(n)},
\delta_{n+1},\delta_{n+1})$ are in one 
$\Br_{n+2}\ltimes\{\pm 1\}^{n+2}$ orbit.

Lemma \ref{t2.5} in \cite{WY20} says that for any $w\in 
\langle s_{\delta_1},...,s_{\delta_n}\rangle$
$(\delta_{\sigma(1)},...,\delta_{\sigma(n)},
\alpha_{n+1},\alpha_{n+1})$ and
$(\delta_{\sigma(1)},...,\delta_{\sigma(n)},
w(\alpha_{n+1}),w(\alpha_{n+1}))$ are in one 
$\Br_{n+2}\ltimes\{\pm 1\}^{n+2}$ orbit.

Therefore it is sufficient to show that 
$w\in \langle s_{\delta_1},...,s_{\delta_n}\rangle$ with
$w(\alpha_{n+1})=\pm \delta_{n+1}$ exists. 

The existence of such a $w$ follows from Corollary 5.3
in \cite{BaWY21}. \hfill$\Box$ 

\bigskip
The last theorem in this section, Theorem \ref{t4.11},
gives a positive answer to the question posed in the Remarks
\ref{t2.8} (i) 
for the tubular elliptic root systems and their Coxeter elements.

\begin{theorem}[{\cite[VI 1.2 Theorem]{Kl87}
\cite[Theorem 1.3]{BaWY18}}]\label{t4.11}
Let $(H_\Z,I,\Phi)$ be a tubular elliptic root system and let 
$c\in W$ be a Coxeter element in it. The set 
\begin{eqnarray}\label{4.1}
\{\uuuu{\alpha}\in\Phi^{n}\,|\, c=s_{\alpha_1}\circ ... \circ
s_{\alpha_n},H_\Z=\sum_{i=1}^n\Z\cdot \alpha_i\}
\end{eqnarray}
of reduced presentations of $c$ with subroot lattice the
full lattice $H_\Z$ is a single $\Br_n\ltimes\{\pm 1\}^n$
orbit in the cases $E_6^{(1,1)},E_7^{(1,1)},E_8^{(1,1)}$,
and it decomposes into two 
$\Br_n\ltimes\{\pm 1\}^n$ orbits in the case $D_4^{(1,1)}$.
\end{theorem}

\begin{remarks}\label{t4.12}
(i) Kluitmann \cite{Kl87} proved Theorem \ref{t4.11}
in the cases $E_6^{(1,1)},E_7^{(1,1)},E_8^{(1,1)}$.
In fact, he proved in these cases the slightly stronger statement
that the set in \eqref{4.1} is a single $\Br_n$ orbit.

(ii) \cite{BaWY18} reproves Theorem \ref{t4.11} for
the cases $E_6^{(1,1)},E_7^{(1,1)},E_8^{(1,1)}$, and proves
it for the case $D_4^{(1,1)}$.

(iii) \cite{BaWY18} does not consider tuples 
$(\alpha_1,...,\alpha_n)$ of roots, but the corresponding
tuples $(s_{\alpha_1},...,s_{\alpha_n})$ of reflections with
the Hurwitz action. Lemma 5.14 in \cite{BaWY18} translates
the condition in Theorem 1.3 in \cite{BaWY18} that
$(s_{\alpha_1},...,s_{\alpha_n})$ is {\it generating}
into the condition in Theorem 4.11 that 
$H_\Z=\sum_{i=1}^n\Z\cdot \alpha_i$.

(iv) This condition is needed in Theorem \ref{t4.11}
(contrary to Theorem \ref{t3.14} where it is not needed).
Example 5.22 in \cite{BaWY18} is an example of a reduced
presentation of a Coxeter element in $E_6^{(1,1)}$ whose
subroot lattice is not the full lattice $H_\Z$.
\end{remarks}

\begin{example}\label{t4.13}
Consider the tubular elliptic root system $D_4^{(1,1)}$ and a Coxeter
element $c=s_{e_1}...s_{e_6}$ for some basis 
$\uuuu{e}\in \Phi^6$
with Dynkin diagram as in Figure \ref{fig4.1}.
Define a unimodular pairing $L:H_\Z\times H_\Z\to\Z$ by 
$L(\uuuu{e}^t,\uuuu{e})^t=S$.
Then $(H_\Z,L,\uuuu{e})$ is a unimodular bilinear lattice
with triangular basis, $I=I^{(0)}$ and $-M=c$. 

Lemma \ref{t2.9} can be used to see that the set in \eqref{4.1}
decomposes in the case $D_4^{(1,1)}$ into at least two orbits.
As this is relevant for the proof of Theorem \ref{t7.4}, 
we give the details. 

It is sufficient to give an automorphism $g$ of $H_\Z$ which
respects $I$ and commutes with $c$, but which
does not respect $L$. Then $\uuuu{e}$ and $g(\uuuu{e})$ are both
in the set in \eqref{4.1}, but are not in the same
$\Br_n\ltimes\{\pm 1\}^n$ orbit. Then $g(\uuuu{e})$ is by
Lemma \ref{t2.9} not even an upper triangular basis 
with respect to $L$. 

$g$ will be constructed as follows. Let $H_{1,\Z}(c)$ and 
$H_{\neq 1,\Z}(c)$ be the sublattices of $H_\Z$ as in 
Definition \ref{t3.8}. For $L$ as defined above, they are
$L$-biorthogonal. Their (direct) sum has finite index in $H_\Z$.
One first calculates $H_{1,\Z}(c)$ using that it is equal to
$\Rad(I)$, and then one calculates $H_{\neq 1,\Z}(c)$ as
the $L$-biorthogonal complement. One finds
\begin{eqnarray*}
H_{1,\Z}(c)&=& \langle f_1,f_2\rangle_\Z,\quad\textup{where }
f_1:=e_5-e_6,\ f_2:=\sum_{i=1}^6e_i,\\
H_{\neq 1,\Z}(c)&=& \langle 2e_1+f_1,e_1-e_2,e_1-e_3,e_1-e_4
\rangle_\Z.
\end{eqnarray*}
Now define $g\in\Aut(H_\Q,I,c)$ as the unique
automorphism of $H_\Q$, which respects the subspaces
$H_{1,\Q}(c)$ and $H_{\neq 1,\Q}(c)$ of $H_\Q$, and which is
on each of them defined by 
\begin{eqnarray*}
g((f_1,f_2))=(f_1,f_2)\begin{pmatrix}-1&0\\0&1\end{pmatrix},
\qquad g|_{H_{\neq 1,\Q}}=\id.
\end{eqnarray*}
Obviously $g$ commutes with $c$ and respects $I$.
It remains to see that $g$ is an automorphism of $H_\Z$ as
$\Z$-lattice and that it does not respect $L$. Write
\begin{eqnarray*}
e_1&=& \frac{1}{2}(2e_1+f_1)-\frac{1}{2}f_1,\\
e_i&=& -(e_1-e_i)+\frac{1}{2}(2e_1+f_1)-\frac{1}{2}f_1
\quad\textup{for }i\in\{2,3,4\},\\
e_5&=& \frac{3}{2}f_1+\frac{1}{2}f_2-(2e_1+f_1)
+\frac{1}{2}((e_1-e_2)+(e_1-e_3)+(e_1-e_4)),\\
e_6&=& \frac{1}{2}f_1+\frac{1}{2}f_2-(2e_1+f_1)
+\frac{1}{2}((e_1-e_2)+(e_1-e_3)+(e_1-e_4)),
\end{eqnarray*}
and see
\begin{eqnarray*}
g(e_i)&=& e_i+f_1\quad\textup{for }i\in\{1,2,3,4\},\\
g(e_5)&=& e_5-3f_1=-2e_5+3e_6,\\
g(e_6)&=& e_6-f_1=-e_5+2e_6.
\end{eqnarray*}
Therefore $g\in\Aut(H_\Z,I^{(0)},c)$. It does not respect $L$ 
because
\begin{eqnarray*}
L((f_1,f_2)^t,(f_1,f_2))&=&
\begin{pmatrix}0&-2\\2&0\end{pmatrix},\\ 
L((g(f_1),g(f_2))^t,(g(f_1),g(f_2))
&=&\begin{pmatrix}0&2\\-2&0\end{pmatrix}.
\end{eqnarray*}
In fact, one calculates that 
$$L(g(\uuuu{e})^t,g(\uuuu{e}))^t\stackrel{!}{=}S^t$$
is just the transpose of the matrix 
$S=L(\uuuu{e}^t,\uuuu{e})^t$. 
\end{example}

\section{Isolated hypersurface singularities}\label{s5}
\setcounter{equation}{0}
\setcounter{table}{0}

An isolated hypersurface singularity 
gives rise to a unimodular bilinear lattice $(H_\Z,L)$ 
of some rank $n\in\N$ and a $\Br_n\ltimes\{\pm 1\}^n$ orbit of 
triangular bases. References for this statement 
and special properties of these data will be
presented in this section.

\begin{definition}\label{t5.1}
(E.g. \cite{AGV88}\cite{Eb01})
An isolated hypersurface singularity is a holomorphic
function germ $f:(\C^{m+1},0)\to(\C,0)$ (with $m\in\Z_{\geq 0}$) 
with $f(0)=0$ and with an isolated singularity at 0, 
i.e. the Jacobi ideal
$J_f:=\bigl(\frac{\paa f}{\paa z_0},...,\frac{\paa f}{\paa z_m}\bigr)
\subset \C\{z_0,...,z_m\}$ has an isolated zero at 0.
Its Milnor number $n$ (most often called $\mu$, but not here) is
\begin{eqnarray}\label{5.1}
n:=\dim\C\{z_0,...,z_m\}/J_f\in\N.
\end{eqnarray}
\end{definition}

See e.g. \cite{Lo84}, \cite{AGV88} or \cite{Eb01} 
for the construction of a good representative 
$f:Y\to\Delta_\eta$ where 
$\Delta_\eta:=\{\tau\in\C\, |\, |\tau|<\eta\}$ is a sufficiently
small disk. 

The relative homology groups (reduced if $m=0$)
$Ml(f,\zeta):=H_{m+1}(Y,f^{-1}(\zeta\eta),\Z)$ with 
$\zeta\in S^1$ are isomorphic  to $\Z^n$ 
\cite[(5.11)]{Lo84} \cite[I.2]{AGV88}, 
and some generators of them can be called (classes of)
{\it Lefschetz thimbles}. 
They form a flat $\Z$-lattice bundle on $S^1$. 
($Ml$ stands for {\it Milnor lattice}.)
An intersection form for Lefschetz thimbles
is well defined on relative homology groups with different 
boundary parts. It is for any $\zeta\in S^1$ a $(-1)^{m+1}$
symmetric unimodular bilinear form
\begin{eqnarray}\label{5.2}
I^{Lef}:Ml(f,\zeta)\times Ml(f,-\zeta)\to\Z
\end{eqnarray}

The intersection form $I^{Lef}$ leads to a 
bilinear form $L^{sing}$ on $Ml(f,\zeta)$ which is called
{\it Seifert form}. The definition below follows 
\cite[I.2]{AGV88}. It is also consistent with \cite{He05}.
But beware that there are different versions of the 
Seifert form in the literature.

Let $\gamma_\pi$ be the isomorphism 
$Ml(f,-\zeta)\to Ml(f,\zeta)$ by flat shift in
mathematically positive direction.
Then the classical Seifert form is given by
\begin{eqnarray}\label{5.3}
L^{sing}: Ml(f,\zeta)\times Ml(f,\zeta)\to\Z,\\
L^{sing}(a,b):=(-1)^{m+1}I^{Lef}(a,\gamma_{\pi} b).\nonumber
\end{eqnarray}
The classical monodromy $M^{sing}$ 
and the intersection form $I^{sing}$ on $Ml(f,\zeta)$ are given by
\begin{eqnarray}\label{5.4}
L^{sing}(M^{sing}a,b)&=&(-1)^{m+1}L^{sing}(b,a),\\
I^{sing}(a,b)&=& -L^{sing}(a,b)+(-1)^{m+1}L^{sing}(b,a)
\label{5.5}\\
&=&L^{sing}((M^{sing}-\id)a,b).\nonumber
\end{eqnarray}

\begin{theorem}[{\cite[I.2]{AGV88}\cite{Eb01}\cite{He05}}]
\label{t5.2}
Consider a singularity $f=f(z_0,...,z_m)$ and construct the
objects $Ml(f,\zeta)$, $I^{Lef}$, $L^{sing}$, $M^{sing}$ and
$I^{sing}$ above.
The {\sf normalized} Seifert form $L$ with
\begin{eqnarray}
L&:=&(-1)^{(m+1)(m+2)/2}\cdot L^{sing}.\label{5.6}
\end{eqnarray}
gives rise to a unimodular bilinear lattice $(H_\Z,L)$ by 
$$H_\Z:=H_\Z(f):=Ml(f,1).$$
Define $k\in\{0;1\}$ by
\begin{eqnarray}\label{5.7}
k:=m\textup{ mod}\, 2,\quad\textup{i.e. }\left\{\begin{array}{ll}
k:=0&\textup{if }m\textup{ is even,}\\
k:=1&\textup{if }m\textup{ is odd.}\end{array}\right.
\end{eqnarray}
The monodromy $M$ and the bilinear form $I^{(k)}$
of the unimodular bilinear lattice $(H_\Z,L)$ 
are related to $M^{sing}$ and $I^{sing}$ by 
\begin{eqnarray}
M&=& (-1)^{m+1}M^{sing},\label{5.8}\\
I^{(k)}&=&(-1)^{m(m+1)/2}I^{sing}.\label{5.9}
\end{eqnarray}
We call $M$ the {\it normalized} monodromy. 
\end{theorem}

An isolated hypersurface singularity comes equipped with a
$\Br_n\ltimes\{\pm 1\}^n$ orbit of triangular bases $\BB^{sing}$
of the unimodular bilinear lattice $(H_\Z,L)$. These bases
are called {\it distinguished bases}.
They can be constructed as follows.
One chooses a small deformation $\www f:\www Y\to\Delta_\eta$
such that $\textup{Sing}(\www f)=\{x^{(1)},...,x^{(n)}\}$
consists of $n$ $A_1$-singularities and their {\it critical 
values} $u_j:=\www f(x^{(j)})$ form a set 
$\Sigma=\{u_1,...,u_n\}\subset \Delta_\eta$ of $n$ pairwise
different points (it is also called {\it morsification}).
One chooses in $\Delta_\eta-\Sigma$ a {\it distinguished system
of paths} $(\gamma_1,...,\gamma_n)$ with 
$\gamma_1(0)=...=\gamma_n(0)=\eta$ and 
$\gamma_j(1)=u_{\sigma(j)}$ for some permutation 
$\sigma\in S_n$ (definition e.g.  in \cite[I.2]{AGV88} 
\cite[5.2]{Eb01} \cite[\S 2 1.5]{AGLV98}).
The family of vanishing cycles (with some chosen sign) over
the path $\gamma_j$ is a Lefschetz thimble 
and induces a unique element
\begin{eqnarray}\label{5.10}
\Gamma_j\in H_\Z=Ml(f,1)\cong Ml(\www f,1):=
H_{m+1}(\www Y,\www f^{-1}(\eta),\Z).
\end{eqnarray}
The family of vanishing cycles with the other possible sign
gives $-\Gamma_j$. 
The tuple $(\Gamma_1,...,\Gamma_n)$ turns out to be a triangular
basis of $(H_\Z,L)$. The set $\BB^{sing}$ consists of 
all in such a way constructed triangular bases.

\begin{remarks}\label{t5.3}
(i) The following isomorphisms hold, 
\begin{eqnarray}\label{5.11}
Ml(\www f,1)\cong H_m(\www f^{-1}(\eta),\Z)
\cong H_m(f^{-1}(\eta),\Z).
\end{eqnarray}
The set on the right hand side is usually called
{\it Milnor lattice}. A distinguished system
of paths induces a distinguished basis of the 
Milnor lattice by pushing vanishing cycles from the critical
points of $\www f$ above the (ordered) paths to $\www f^{-1}(\eta)$.
The proof of this was first fixed by Brieskorn 
\cite[Appendix]{Br70}, see also
\cite[I.2]{AGV88} \cite[Satz 5.5 in 5.4, and 5.6]{Eb01}
\cite[\S 1 1.5]{AGLV98}.
The family of vanishing cycles above one such path
forms a Lefschetz thimble. 

(ii) That $\BB^{sing}$ is one orbit of 
$\Br_n\ltimes\{\pm 1\}^n$ is rather trivial, as 
$\Br_n$ acts on the set of (homotopy classes of)
distinguished systems of paths, see \cite[I.2]{AGV88} 
\cite[5.7]{Eb01} \cite[\S 2 1.9]{AGLV98}.

(iii) In \cite[page 463, after Def. 8.4.5]{Eb20} it is claimed 
that the set $\BB^{sing}$ depends on the chosen 
morsification $\www{f}:\www{Y}\to\Delta_\eta$ and is thus not
canonical. That is not true. The set $\BB^{sing}$ is independent
of choices and canonical.

(iv) In \cite[5.6]{Eb01} the bilinear form $L^{sing}$ 
is called $v$. Korollar 5.3 (i) in \cite{Eb01} that
the matrix $V$ is upper triangular, is wrong.
It is lower triangular. 

(v) If one adds to $f(x_0,...,x_m)$ a square $x_{n+1}^2$ in a 
new variable, the new function $f^{stab}(x_0,...,x_{n+1})
=f(x_0,...,x_n)+x_{n+1}^2$ on $(\C^{m+2},0)$ 
is called a {\it suspension} of $f$. 
It is a special case of the Thom-Sebastiani 
construction, which associates to two singularities
$f$ and $g$ in different variables their sum $f+g$.

The Milnor lattices $Ml(f)$ and $Ml(f^{stab})$
are almost canonically isomorphic, namely there are two
canonical isomorphisms $Ml(f)\to Ml(f^{stab})$, which differ
only by the sign. Both isomorphisms map the normalized
Seifert form, the normalized monodromy and the set 
$\BB^{sing}$ of $f$ to the corresponding objects of $f^{stab}$
\cite[I.2.7]{AGV88}. This result for the monodromy is due to
Thom and Sebastiani, for the Seifert form \cite[I.2.7]{AGV88}
cites Deligne (without a reference), for the distinguished bases
it is due to Gabrielov \cite{Ga73}.
Therefore $f$ and $f^{stab}$ give rise to the same
unimodular bilinear lattice $(H_\Z,L)$ and the same set
$\BB^{sing}\subset\BB^{tri}$. 
\end{remarks}

The simple and the simple elliptic singularities are
the only singularities where $I^{(0)}$ is positive definite
respectively positive semidefinite. Because of this they are the subject
of Theorem \ref{t1.3} and Theorem \ref{t1.5}.
Therefore we recall here their definitions. 
We define them in the smallest number of variables.
As in Remark \ref{t5.3}, one can add squares in additional
variables, i.e. one can go over to suspensions.

\begin{definition}\label{t5.4}
(a) (E.g. \cite[II 15.1]{AGV85})
The {\it simple singularities} (or {\it ADE-singularities}) 
consist of two series and three exceptional cases. 
Normal forms of them are as follows.
\begin{eqnarray*}
\begin{array}{c|c|c|c|c}
A_n\ (n\geq 1) & D_n\ (n\geq 4) & E_6 & E_7 & E_8 \\ \hline
z_0^{n+1} & z_0^{n-1}+z_0z_1^2 & z_0^4+z_0^3 & z_0^3z_1+z_1^3 & 
z_0^5+z_1^3 
\end{array}
\end{eqnarray*}

(b) (E.g. \cite[1.9]{SaK74}) The {\it simple elliptic singularities}
are three 1-parameter families of isolated hypersurface 
singularities. The following normal forms are the {\it Legendre
normal forms}. The parameter is $\lambda\in\C-\{0;1\}$.
\begin{eqnarray*}
\begin{array}{c|c|c}
\textup{name} & n & \textup{Legendre normal form} \\ \hline 
\www{E}_6 & 8 & f_\lambda=z_1(z_1-z_0)(z_1-\lambda z_0)-z_0z_1^2\\
\www{E}_7 & 9 & f_\lambda=z_0z_1(z_1-z_0)(z_1-\lambda z_0)\\
\www{E}_8 & 10 & f_\lambda=z_1(z_1-z_0^2)(z_1-\lambda z_0^2)
\end{array}
\end{eqnarray*}
\end{definition}

Theorem \ref{t5.5} collects some known properties of
the unimodular bilinear lattice $(H_\Z(f),L)$ and its
set $\BB^{sing}(f)$ of distinguished bases.
The parts (a)--(e) and (h) are of a general nature.
The parts (f) and (g) make statements on those singularities
where $I^{(0)}$ is positive definite or positive 
semidefinite.

\begin{theorem}\label{t5.5}
Consider a singularity $f=f(z_0,...,z_m)$ 
and the associated unimodular bilinear lattice $(H_\Z,L)$ 
with set $\BB^{sing}$ of distinguished bases.

(a) \cite{Br70} The monodromy
$M$ is quasiunipotent, and the sizes of the
Jordan blocks are at most $m+1$. The sizes
of the Jordan blocks with eigenvalue 1 are at most $m$.

(b) \cite{AC73} $\tr(M)=1$.

(c) \cite{Ga74}\cite{La73}\cite{Le73} The {\rm CDD}
of any distinguished basis is connected. 

(d) (Classical, e.g. \cite{Ga79}\cite{Eb01})
For $k\in\{0;1\}$ the set $\Delta^{(k)}$ of vanishing
cycles is a single $\Gamma^{(k)}$ orbit (see Definition
\ref{t2.3} (f) and (g) for $\Gamma^{(k)}$ and $\Delta^{(k)}$).

(e) \cite{Ga79} The $\Br_n\ltimes\{\pm 1\}^n$ orbit
of $S$ contains a matrix $\www{S}\in T_n(\Z)$ 
with all entries in $\{0;1,-1\}$.

(f) \cite{Tj68}\cite{Ar73-1} The only
singularities with $I^{(0)}$ positive definite
are the {\sf ADE-singularities}.

(Classical, known long before \cite{Ar73-1})
For each of them $(H_\Z,I^{(0)})$ is an 
ADE root lattice, $-M$ is a Coxeter element
and $\BB^{sing}$ contains a root basis.

(g) \cite{Ar73-2} The only
singularities with $I^{(0)}$ positive semidefinite
and not positive definite 
are the three 1-parameter families $\www E_6,\www E_7$
and $\www E_8$ of {\sf simple elliptic singularities}.

\cite{Ga73} $\BB^{sing}$ contains a $\Z$-basis such that
the {\rm CDD} of its matrix $S$ is as in Figure \ref{fig4.1},
where $\www E_k\sim E_k^{(1,1)}$ for $k\in\{6,7,8\}$. 
Therefore $(H_\Z,I^{(0)},\Phi)$ with  
$\Phi=\{\delta\in H_\Z\,|\, I^{(0)}(\delta,\delta)=2\}$
is a tubular elliptic root system of type $E_k^{(1,1)}$
in the case of a simple elliptic singularity of type 
$\www{E_k}$, and $-M$ is a Coxeter element in it. 

(h) (Classical and \cite{Eb18}) The following holds.

\medskip
\begin{tabular}{l|l|l|l}
$f$ & $S_{ij}$ & $|\BB^{sing}|$ & 
$|\{\textup{CDD's}\}|$ \\ \hline 
ADE-singularity & in $\{0,\pm 1\}$ & finite & finite \\
simple elliptic sing. &  in $\{0,\pm 1,\pm 2\}$ & 
infinite & finite \\
any other singularity & unbounded & infinite & infinite 
\end{tabular}
\end{theorem}

\begin{remarks}\label{t5.6}
(i) In part (h) the third line is due to \cite{Eb18}.
$|\BB^{sing}|=\infty$ in the second line is proved 
for example in \cite{Kl87} and in \cite{HR21}.
The first line and the other statements in the second line 
follow from Lemma \ref{t2.11}.

(ii) Property (e) is not so well known. It follows from
Gabrielov's iterative construction of a CDD of a
singularity from the CDD of a smaller singularity.
We will not make use of it. Still it is interesting. 
Together with Lemma \ref{t2.10} it implies for any singularity
with $n\geq 2$, so for any singularity except the 
$A_1$ singularity, that $\BB^{sing}$ is a single $\Br_n$ orbit,
i.e. all sign changes can be carried out by braids.

(iii) Deligne's result Theorem \ref{t3.14} for the 
Coxeter elements in the irreducible simply laced root lattices
implies for the ADE-singularities
\begin{eqnarray}\label{5.12}
\BB^{sing} = \{\uuuu{\alpha}\in(\Delta^{(0)})^n\,|\, 
-M=s_{\alpha_1}\circ ...\circ s_{\alpha_n}\}.
\end{eqnarray}
Here $\Delta^{(0)}=\Phi=\{\delta\in H_\Z\,|,\, 
I^{(0)}(\delta,\delta)=2\}$. 

(iv) Kluitmann's result Theorem \ref{t4.11} for the Coxeter
elements in the tubular elliptic root systems $E_k^{(1,1)}$ 
$(k\in\{6,7,8\})$ implies for the simple elliptic singularities
\begin{eqnarray}\label{5.13}
\BB^{sing} = \{\uuuu{\alpha}\in(\Delta^{(0)})^n\,|\, 
-M=s_{\alpha_1}\circ ...\circ s_{\alpha_n}, 
H_\Z=\sum_{i=1}^n\Z\cdot \alpha_i\},
\end{eqnarray}
Here $\Delta^{(0)}=\Phi=\{\delta\in H_\Z\,|\, 
I^{(0)}(\delta,\delta)=2\}$ (see Theorem \ref{t4.3}). 

(v) Baumeister, Wegener and Yahiatene proved the analogous
statement \eqref{5.13} also for the {\sf hyperbolic}
singularities $T_{pqr}$ \cite[Theorem 1.5]{BaWY21}.

(vi) The results in (iii)--(v) for the ADE-singularities, the 
simple elliptic singularities and the hyperbolic singularities are the only cases, where one has a simple explicit
description of the set $\BB^{sing}$.

{\bf Open problem:} {\it Prove or disprove \eqref{5.13}
for other singularities}.
\end{remarks}

\section{The positive definite cases}\label{s6}
\setcounter{equation}{0}
\setcounter{table}{0}

The following 1-1 correspondence is fairly elementary, except
for the involvement of Voigt's part of Theorem \ref{t3.14}.

\begin{lemma}\label{t6.1}
(a) For each $n\in\N$ there is a 1-1 correspondence between the
following two sets.
\begin{list}{}{}
\item[(S1)] The $\Br_n\ltimes\{\pm 1\}^n$ orbits of matrices 
$S\in T_n(\Z)$
with $S+S^t$ positive definite and ${\rm CDD}(S)$ connected.
\item[(S2)] The pairs $((H_\Z,I),Cl(w))$ where $(H_\Z,I)$
is an irreducible simply laced root lattice, $w$ is a
quasi-Coxeter element and $Cl(w)$ denotes the conjugacy class
of $w$ with respect to the Weyl group.
\end{list}

(b) A matrix $S$ as in (S1) is a distinguished matrix of 
an isolated hypersurface singularity $\iff$ $w$ in (S2) is
a Coxeter element. Then $S$ is a distinguished matrix of an
ADE-singularity.

(c) The data in (S1) define a triple
$(H_\Z,L,\BB^{tri})$ (up to isomorphism) where $(H_\Z,L)$ is a
unimodular bilinear lattice and the set $\BB^{tri}$ of triangular
bases is a single $\Br_n\ltimes\{\pm 1\}^n$ orbit.
By Remark \ref{t2.8} (ii) this triple induces $\Gamma^{(k)}$
and $\Delta^{(k)}$ for $k\in\{0;1\}$. Then
\begin{eqnarray}\label{6.1}
((H_\Z,I^{(0)}),\Gamma^{(0)}\textup{ and }\Delta^{(0)}
\textup{ in }(S1))\cong
((H_\Z,I),W\textup{ and }\Phi\textup{ in }(S2)).
\end{eqnarray}
Furthermore here $\Delta^{(0)}=R^{(0)}$.
\end{lemma}

{\bf Proof:} (a) From (S1) to (S2):
A matrix $S$ in (S1) gives rise to a unimodular
bilinear lattice $(H_\Z,L,\uuuu{e})$ with a triangular base
$\uuuu{e}$. Recall $I^{(0)}=L+L^t$ and 
$I^{(0)}(\uuuu{e}^t,\uuuu{e})=S+S^t$. Therefore 
$(H_\Z,I^{(0)})$ is a simply laced root lattice,
and $\uuuu{e}$ is a $\Z$-basis of roots. 
Thus the element $-M=s_{e_1}^{(0)}\circ ...\circ s_{e_n}^{(0)}$
is a quasi-Coxeter element.
$(H_\Z,I^{(0)})$ is irreducible because $(H_\Z,L,\uuuu{e})$ is
irreducible. 
A different matrix $\www{S}$ in the $\Br_n$ orbit of $S$ gives
an isomorphic unimodular bilinear lattice 
$(\www{H_\Z},\www{L})$ with a triangular basis $\www{\uuuu{e}}$.
One can choose the isomorphism such that $\www{\uuuu{e}}$ is
obtained by the same element of $\Br_n\ltimes\{\pm 1\}^n$  
from $\uuuu{e}$ as $\www{S}$ is obtained from $S$. 
Then $-\www{M}=-M$.  An arbitrary isomorphism 
$(H_\Z,L)\to (\www{H_\Z},\www{L})$ maps the quasi-Coxeter element 
$-M$ to a conjugate of the quasi-Coxeter element $-\www{M}$. 

From (S2) to (S1): Let $(H_\Z,I)$ be an irreducible 
simply laced root lattice of rank $n$ and $w$ a quasi-Coxeter
element in it. By Voigt's part of Theorem \ref{t3.14} the set
$$\{\uuuu{\alpha}\in\Phi^n\,|\, w=s_{\alpha_1}\circ ...\circ
s_{\alpha_n}\}$$ 
is a single $\Br_n\ltimes\{\pm 1\}^n$ orbit 
(in fact, even a single $\Br_n$ orbit if $n\geq 2$, 
but we don't need that).
Each such $\uuuu{\alpha}$ is a $\Z$-basis of $H_\Z$ because
$w$ is a quasi-Coxeter element. Choose one such $\uuuu{\alpha}$.
It determines uniquely a matrix $S\in T_n(\Z)$ with
$I(\uuuu{\alpha}^t,\uuuu{\alpha})=S+S^t$. 
Obviously $S+S^t$ is positive definite. ${\rm CDD}(S)$ is connected
because the root lattice $(H_\Z,I)$ is irreducible. 
Another choice $\www{\uuuu{\alpha}}$ is obtained from 
$\uuuu{\alpha}$ by an element of 
$\Br_n\ltimes\{\pm 1\}^n$.
The same element maps $S$ to the matrix $\www{S}$ with
$I^{(0)}(\www{\uuuu{\alpha}}^t,\www{\uuuu{\alpha}})=
\www{S}+\www{S}^t$. Therefore one obtains a
$\Br_n\ltimes\{\pm 1\}^n$ orbit of matrices in (S1).

It is fairly obvious that the constructions
$(S1)\to(S2)$ and $(S2)\to(S1)$ are inverse.

(b) This follows from part (a) and Theorem \ref{t5.5} (f).

(c) If we start with (S1), we also
have the data in (S2). Then by \eqref{2.3} and Theorem \ref{t3.14}
the set $\BB^{tri}$ is a single $\Br_n\ltimes\{\pm 1\}^n$ orbit.

Now we start with (S2). Let $\uuuu{\alpha}\in\Phi^n$
be a reduced presentation of $w$. It is well known that
$\langle s_{\alpha_1},...,s_{\alpha_n}\rangle =W$ and that
$W\{\pm\alpha_1,...,\pm\alpha_n\}=\Phi$. Everything else follows
from the proof of part (a). 
\hfill$\Box$ 

\bigskip
Question \ref{t1.1} in the case of matrices $S$
with $S+S^t$ positive definite and ${\rm CDD}(S)$ connected
amounts thus only to the separation of the matrices $S$ 
which correspond to Coxeter elements from the matrices 
which correspond to quasi-Coxeter elements which are not 
Coxeter elements. 
Theorem \ref{t6.2} offers two ways to separate them,
first the trace, second the variance inequality. Both work here.

\begin{theorem}\label{t6.2}
Let $S\in T_n(\Z)$ for some $n\in\N$ with $S+S^t$ 
positive definite and ${\rm CDD}(S)$ connected, so as in (S1)
in Lemma \ref{t6.1} (a).
Let $(H_\Z,I,w)$ be the irreducible simply laced root lattice
with one quasi-Coxeter element $w$ which corresponds 
($w$ only up to conjugation) to $S$ by Lemma \ref{t6.1} (a).

(a) If $w$ is a Coxeter element then $\tr(S^{-1}S^t)=1$.\\
If $w$ is not a Coxeter element then $\tr(S^{-1}S^t)\leq 0$.

(b) Let $\Sp(S)=(\alpha_1,...,\alpha_n)$ be the spectrum
of $S$ with $\alpha_1\leq ...\leq \alpha_n$, 
which was defined in Lemma \ref{t2.11} (c).\\\
If $w$ is a Coxeter element then
$\Var(\Sp(S))=\frac{1}{12}(\alpha_n-\alpha_1)$.\\
If $w$ is not a Coxeter element then
$\Var(\Sp(S))>\frac{1}{12}(\alpha_n-\alpha_1)$.
\end{theorem}

{\bf Proof:}
(a) Let $(H_\Z,L,\uuuu{e})$ be the unimodular bilinear
lattice with triangular basis which is associated to $S$.
By the proof of Lemma \ref{t6.1} (a) the monodromy $M$ 
of the triple $(H_\Z,L,\uuuu{e})$ is $-w$, and 
$$\tr(S^{-1}S^t)=\tr(M)=\tr(-w)=-\tr(w).$$
Table \ref{tab3.5} shows that $\tr(w)=-1$ if $w$
is a Coxeter element and that $\tr(w)\geq 0$ if $w$
is  quasi-Coxeter element, but not a Coxeter element.

(b) The same Table \ref{tab3.5} shows the characteristic
polynomial of each quasi-Coxeter element.
Its spectrum and the variance of its spectrum are given
by Definition \ref{t6.3} and can be calculated with
Lemma \ref{t6.4}. 

Table \ref{tab6.1} gives the results.
Column 1 gives the irreducible simply laced root lattices.
Column 2 gives the conjugacy classes of the quasi-Coxeter 
elements.
Column 3 gives 12 times the non-normalized variance
$V_2(\Sp(S))$ (see Definition \ref{t6.3}). 
Column 4 gives $\alpha_1+\frac{1}{2}$.
Column 5 gives the difference 
$n(\alpha_n-\alpha_1)=n(1-2(\alpha_1+\frac{1}{2}))$.
Column 6 gives the equality or inequality between
$12\cdot V_2(\Sp(S))$ and $n(\alpha_n-\alpha_1)$.
One sees that the equality 
$\Var(\Sp(S))=\frac{1}{12}(\alpha_n-\alpha_1)$
holds for the Coxeter elements,
and the inequality 
$\Var(\Sp(S))>\frac{1}{12}(\alpha_n-\alpha_1)$ holds
for the other quasi-Coxeter elements. 

\begin{table}
\begin{eqnarray*}
\begin{array}{l|l|l|l|l|c}
 & & 12\cdot V_2(\Sp(S)) &\alpha_1+\frac{1}{2}& 
 n\cdot (\alpha_n-\alpha_1)& \\ \hline
A_n & A_n & \frac{n(n-1)}{n+1} & \frac{1}{n+1}& \frac{n(n-1)}{n+1}& = \\ \hline
D_n & D_n & \frac{n(n-2)}{n-1}& \frac{1}{2(n-1)}& \frac{n(n-2)}{n-1}& =\\
% & D_n(a_1) & \frac{(n-1)(n-3)}{n-2}+\frac{3\cdot 1}{2} & 
% \frac{1}{2(n-2)}& \frac{n(n-3)}{n-2}& >\\
% & \vdots & \vdots & \vdots & & \\
 & D_n(a_j) &  n-\frac{1}{n-1-j}-\frac{1}{j+1}& 
 \frac{1}{2(n-1-j)}& \frac{n(n-2-j)}{n-1-j}& >\\
 & \textup{for }1\leq j<[\frac{n}{2}]& & & & \\ \hline
E_6 & E_6 & 5 & { } {1}/{12}& 5& =\\
 & E_6(a_1) & { } {50}/{9}& { } {1}/{9}& { } {14}/{3}&>\\
 & E_6(a_2) & 6 & { } {1}/{6}& 4&>\\ \hline
E_7 & E_7 & { } {56}/{9} & { } {1}/{18}& { } {56}/{9}&= \\
 & E_7(a_1) & { } {48}/{7}& { } {1}/{14}& 6&>\\
 & E_7(a_2) & 7& { } {1}/{12}& { } {35}/{}6&>\\
 & E_7(a_3) & { } {112}/{15}& { } {1}/{10}& { } {28}/{5}&>\\
 & E_7(a_4) & 8& { } {1}/{6}& { } {14}/{3}&>\\ \hline
E_8 & E_8 & { } {112}/{15}& { } {1}/{30}& { } {112}/{15}&= \\
 & E_8(a_1) & { } {49}/{6}& { } {1}/{24}& { } {22}/{3}&>\\
 & E_8(a_2) & { } {42}/{5}& { } {1}/{20}& { } {36}/{5}&>\\
 & E_8(a_3) & { } {61}/{6}& { } {1}/{12}& { } {20}/{3}&>\\
 & E_8(a_4) & { } {80}/{9}& { } {1}/{18}& { } {64}/{9}&>\\
 & E_8(a_5) & { } {136}/{15}& { } {1}/{15}& { } {104}/{15}&>\\
 & E_8(a_6) & { } {48}/{5}& { } {1}/{10}& { } {32}/{5}&>\\
 & E_8(a_7) & { } {29}/{3}& { } {1}/{12}& { } {20}/{3}&>\\
 & E_8(a_8) & 8 & { } {1}/{6}& { } {16}/{3}&>
\end{array}
\end{eqnarray*}
\caption[Table 5.1]{Variance (in)equality for the spectra of 
the quasi-Coxeter elements of the irreducible simply laced
root lattices}
\label{tab6.1}
\end{table}
\hfill$\Box$

\begin{definition}\label{t6.3}
Let $g\in\C[t]$ be a unitary polynomial of degree
$n\in\N$ whose zeros are in $S^1$ and such that the multiplicity
of $1$ as a zero is even. 
Its spectrum $\Sp(g)$ is the unique tuple
$(\alpha_1,...,\alpha_n)\in [-\frac{1}{2},\frac{1}{2}]^n$ such 
that $e^{-2\pi i(\alpha_1+1/2)},...,e^{-2\pi i(\alpha_n+1/2)}$
are the roots of $g$, and $\alpha_1\leq...\leq \alpha_n$,
and $-\frac{1}{2}$ and $\frac{1}{2}$ turn up with the same 
multiplicity (which is half of the multiplicity of 1 as a zero 
of $g$).

The {\sf variance} $\Var(g)$ and the {\sf non-normalized variance}
$V_2(g)$ are defined as 
\begin{eqnarray}\label{6.2}
\Var(g):=\frac{1}{n}\sum_{j=1}^n\alpha_j^2,\qquad
V_2(g):=n\cdot \Var(g)=\sum_{j=1}^n\alpha_j^2.
\end{eqnarray}
\end{definition}

\begin{lemma}\label{t6.4}
(a) (Lemma) If $S\in T_n(\Z)$ is positive semidefinite 
and $g$ is the characteristic polynomial of $(-1)S^{-1}S^t$,
then $\Sp(S)=\Sp(g)$, $\Var(\Sp(S))=\Var(g)$
and $V_2(\Sp(S))=V_2(g)$. 

(b) (Lemma)
$V_2$ is additive, 
\begin{eqnarray*}
V_2(g_1\cdot g_1)&=&V_2(g_1)+V_2(g_2),\\
V_2(g_1/g_2)&=&V_2(g_1)-V_2(g_2)
\quad(\textup{\it in the case }g_2|g_1).
\end{eqnarray*}

(c) 
\begin{eqnarray}\label{6.3}
12\cdot V_2(t^n+1)&=& \frac{(n+1)(n-1)}{n},\\
12\cdot V_2(\frac{t^n-1}{t-1})&=& \frac{(n-1)(n-2)}{n}.\label{6.4}
\end{eqnarray}
\end{lemma}

{\bf Proof:}
(a) This follows with Lemma \ref{t2.11}. (b) Clear. (c) 
\begin{eqnarray*}
\Sp(t^n+1)+\frac{1}{2}&=&(\frac{2j-1}{2n}\, |\, j=1,...,n),\\
\Sp(\frac{t^n-1}{t-1})+\frac{1}{2}&=&(\frac{j}{n}\, |\, j=1,...,n-1).
\end{eqnarray*}
With the well-known formulas 
$$\sum_{j=1}^n(2j-1)^2=\frac{(2n-1)2n(2n+1)}{6}\quad\textup{and}\quad
\sum_{j=1}^n j^2=\frac{n(n+1)(2n+1)}{6}$$
one easily calculates
\begin{eqnarray*}
12\cdot V_2(t^n+1)
&=& 12\cdot \sum_{j=1}^n(\frac{2j-1}{2n}-\frac{1}{2})^2=...=\frac{(n+1)(n-1)}{n},\\
12\cdot V_2(\frac{t^n-1}{t-1})&=& 
12\cdot \sum_{j=1}^{n-1}(\frac{j}{n}-\frac{1}{2})^2=...=\frac{(n-1)(n-2)}{n}.
\qquad \hfill\Box
\end{eqnarray*}

\section{The positive semidefinite cases}\label{s7}
\setcounter{equation}{0}
\setcounter{table}{0}

This section will answer Question \ref{t1.1} for the
matrices $S\in T_n(\Z)$ with $S+S^t$ positive
semidefinite, but not positive definite, and 
${\rm CDD}(S)$ connected. Theorem \ref{t7.4}
is part (b) of Theorem \ref{t1.5}.  

Lemma \ref{t7.1} formulates first observations.
It starts with a unimodular bilinear lattice $(H_\Z,L,\uuuu{e})$
such that $I^{(0)}$ is positive semidefinite, but not positive
definite. It relates data on $H_\Z$ with induced data on the
quotient lattice $H_\Z/\Rad(I^{(0)})$.  

Theorem \ref{t7.2} is the analogue in the positive semidefinite
cases of Theorem \ref{t6.2} (a), which studied in the positive 
definite cases the condition $tr(S^{-1}S^t)=1$.
Though Theorem \ref{t7.2} is less satisfying. There are still
quite some families of matrices $S\in T_n(\Z)$ which satisfy
$\tr(S^{-1}S^t)\geq 1$, $S+S^t$ positive semidefinite and CDD
connected. The proof of Theorem \ref{t7.2} uses already
heavily results of section \ref{s3}, especially the Theorems
\ref{t3.6} and \ref{t3.13}. 

Theorem \ref{t7.3} uses the same machinery. It studies for 
which matrices $S\in T_n(Z)$ with CDD connected and
$S+S^t$ positive semidefinite their spectra satisfy the 
variance inequality \eqref{1.1} (=\eqref{7.3}). 
The classification is satisfying. Only the four tubular
elliptic root systems are relevant. 

But the final steps in Theorem \ref{t7.4}, which lead to the
four tubular elliptic root systems and their Coxeter elements,
use also the results of section \ref{s4}. At the end, the
trace condition $\tr(S^{-1}S^t)=1$ allows to discard 
$D_4^{(1,1)}$ as there $\tr(S^{-1}S^t)=2$.

\begin{lemma}\label{t7.1}
Let $S\in T_n(\Z)$ for some $n\in\N$ with $S+S^t$
positive semidefinite, but not positive definite, and
${\rm CDD}(S)$ connected. 
Let $(H_\Z,L,\uuuu{e})$ be the associated unimodular bilinear 
lattice with a triangular basis with 
$S=L(\uuuu{e}^t,\uuuu{e})^t$. So $(H_\Z,L,\uuuu{e})$ is
irreducible, and $I^{(0)}$ is positive semidefinite, but
not positive definite.

Then $\{0\}\subsetneqq\Rad(I^{(0)})=\ker(M+\id)\subset H_\Z$.
Denote by $\oooo{H_\Z}:=H_\Z/\Rad(I^{(0)})$ the quotient lattice. 
For any structure/object on/in $H_\Z$ denote by an overline 
the induced structure/object on/in $\oooo{H_\Z}$. 
The pairing $I^{(0)}$ induces a pairing $\oooo{I^{(0)}}$
on $\oooo{H_\Z}$. 
The pair $(\oooo{H_\Z},\oooo{I^{(0)}})$ is an irreducible 
simply laced root lattice with Weyl group $W=\oooo{\Gamma^{(0)}}$.
We write $s_\alpha:=s_\alpha^{(0)}\in \Gamma^{(0)}$ 
for $\alpha\in R^{(0)}$,
and $s_{\oooo{\alpha}}:=\oooo{s_{\alpha}^{(0)}}\in W$. 

Minus the monodromy 
$$-M=s_{e_1}\circ...\circ s_{e_n}\in\Gamma^{(0)}$$
acts as identity on $\Rad I^{(0)}$ and induces an automorphism
$$\oooo{-M}=s_{\oooo{e_1}}\circ...\circ s_{\oooo{e_n}}\in W.$$
$\oooo{-M}$ is semisimple. The difference 
$k_a:=\rank\oooo{H_\Z}-l(\oooo{-M})\in\Z_{\geq 0}$ is by 
Lemma \ref{t3.9} (b) the multiplicity of 1 as eigenvalue of
$\oooo{-M}$. There is a $k_b\in\Z_{\geq 0}$ such that 
$2(k_a+k_b)$ is the multiplicity of $1$ as generalized
eigenvalue of $-M$. Then
\begin{eqnarray}\label{7.1}
\oooo{n}:=\rank\oooo{H_\Z}=l(\oooo{-M})+k_a, \quad 
n=l(\oooo{-M})+2(k_a+k_b).
\end{eqnarray}

$(\oooo{e_1},...,\oooo{e_n})$ is a non-reduced presentation
of $\oooo{-M}$ whose subroot lattice is the full lattice 
$\oooo{H_\Z}$. 
Especially $k_a+k_b\geq k_5(\oooo{H_\Z},\oooo{-M})$. 

The trace of $-M$ is
\begin{eqnarray}\label{7.2}
\tr(\oooo{-M}|_{\oooo{H_{\neq 1,\Q}}})+2(k_a+k_b).
\end{eqnarray}
\end{lemma}

{\bf Proof:} \eqref{2.1} says that 
$\Rad I^{(0)}=\ker(M+\id)$. 
Obviously $(\oooo{H_\Z},\oooo{I^{(0)}})$ is a simply laced
root lattice. It is irreducible because 
$\oooo{e_1},....,\oooo{e_n}$ generate the full lattice 
$\oooo{H_\Z}$ and because ${\rm CDD}(S)$ is connected.

By Lemma \ref{t2.11} (b) the multiplicity of $1$ as a
generalized eigenvalue of $-M$ is even. $k_a$ is the number
of $2\times 2$ Jordan blocks with eigenvalue $1$ of $-M$.
Therefore the multiplicity of $1$ as generalized eigenvalue is
$2(k_a+k_b)$ for some $k_b\in\Z_{\geq 0}$,
and $\dim \ker(M+\id)=k_a+2k_b$. The rest is clear.
\hfill$\Box$ 

\bigskip
We would like to single out those matrices $S$ in Lemma 
\ref{t7.1} which come from a simple elliptic singularity.
First we will study which matrices satisfy $\tr(S^{-1}S^t)=1$
or $\geq 1$. Theorem \ref{t7.2} will tell that there are few, 
but more than those from the simple elliptic singularities.

\begin{theorem}\label{t7.2}
Let $S\in T_n(\Z)$ and $(H_\Z,L,\uuuu{e})$ be as in
Lemma \ref{t7.1}. Suppose additionally $\tr(S^{-1}S^t)\geq 1$.
The notations and results in Lemma \ref{t7.1} are used, 
especially $(\oooo{H_\Z},\oooo{I^{(0)}})$ is an irreducible
simply laced root lattice of rank $\oooo{n}=l(\oooo{-M})+k_a$.

Let $U\subset\oooo{H_\Z}$ be the subroot lattice associated
to a reduced presentation of $\oooo{-M}$.
Table \ref{tab7.1} lists information on all possibilities
for $((\oooo{H_\Z},\oooo{I^{(0)}}),\oooo{-M},U,k_a,k_b)$
and gives $\tr(S^{-1}S^t)$. 

\begin{table}
\begin{eqnarray*}
\begin{array}{l|l|l|l|l|l}
{}\oooo{H_\Z} & U & k_a & k_b & k_a+k_b & \tr(S^{-1}S^t) \\ 
\hline
A_{\oooo{n}} & - & - & - & - & -\\ \hline 
D_{\oooo{n}} & r_3D_2 \textup{ with} & 0 & r_3-1 & r_3-1 & 2 \\ 
 & \oooo{n}=2r_3,r_3\geq 2 & & & & \\ \hline 
D_{\oooo{n}} & r_4D_2+D_b\textup{ with} & 0 & r_4 & r_4 & 1 \\ 
 & \oooo{n}=2r_4+b,r_4\geq 1, b\geq 3 & & & & \\ \hline 
D_{\oooo{n}} & A_c+r_3D_2\textup{ with} & 1 & r_3-1 & r_3 & 1 \\ 
 & \oooo{n}=c+1+2r_3,r_3\geq 1,c\geq 1 & & & & \\ \hline 
E_6 & 3A_2 & 0 & 1 & 1 & 1 \\ \hline
E_6 & 2A_2+A_1 & 1 & 0 & 1 & 1 \\ \hline
E_7 & 2A_3+A_1 %\sim D_4(a_1)+3A_1 
& 0 & 1 & 1 & 1 \\ \hline
E_7 & 7A_1 & 0 & 3 & 3 & 1 \\ \hline
E_7 & A_3+A_2+A_1 & 1 & 0 & 1 & 1 \\ \hline
E_8 & A_5+A_2+A_1 & 0 & 1 & 1 & 1 \\ \hline
E_8 & A_4+A_2+A_1 & 1 & 0 & 1 & 1 
\end{array}
\end{eqnarray*}
\caption[Table 7.1]{Possibilities for 
$((\oooo{H_\Z},\oooo{I^{(0)}}),\oooo{-M},U,k_a,k_b)$
in Theorem \ref{t7.2}}
\label{tab7.1}
\end{table}

So $(\oooo{H_\Z},\oooo{I^{(0)}})$ is never of type 
$A_{\oooo{n}}$.

In all cases, the restriction of $\oooo{-M}$ to the
subroot lattice $U$ is a Coxeter element
of this subroot lattice, not another quasi-Coxeter element.
(This is an additional information only in the case
$D_{\oooo{n}}$ with $\oooo{n}=2r_4+b$ and $U$ of type
$r_4D_2+D_b$ and $b\geq 4$.)
\end{theorem}

{\bf Proof:}
Let $U\subset\oooo{H_\Z}$ be the subroot lattice associated
to a reduced presentation of $\oooo{-M}$.
Then $U_\Q=\oooo{H_{\neq 1,\Q}}$.
Recall formula \eqref{7.2} for $\tr(-M)$. 
We are interested in the cases where $\tr(-M)\leq -1$. 

Theorem \ref{t6.2} (a) implies the following.
Each irreducible component of $(U,\oooo{I^{(0)}}|_U)$
gives a contribution $-1$ to $\tr(-M)$, if
$\oooo{-M}|_{\textup{this component}}$ is a Coxeter
element, and a contribution $\geq 0$ if 
$\oooo{-M}|_{\textup{this component}}$ is a quasi
Coxeter element, which is not a Coxeter element.

First case, $(\oooo{H_\Z},\oooo{I^{(0)}})$ of type
$A_{\oooo{n}}$: Recall Theorem \ref{t3.6} (a).
Then $U$ is of type $\sum_{i=1}^{r_1}A_{c_i}$
with $r_1\in\N$, $r_2\in\Z_{\geq -1}$, $c_i\in\N$, 
$\oooo{n}=r_1+r_2+\sum_{i=1}^{r_1}c_i$. Then
\begin{eqnarray*}
\tr(-M)&=&r_1(-1)+2(k_a+k_b)
\geq -r_1+2k_5(\oooo{H_\Z},\oooo{-M})\\
&\stackrel{\textup{Theorem \ref{t3.11}}}{=}&
-r_1+2k_1 \stackrel{\textup{Table \ref{tab3.1}}}{=} 
-r_1+2(r_1+r_2)=r_1+2r_2.
\end{eqnarray*}
This is $\geq -1$. It equals $-1$ only if
$r_1=1,r_2=-1$ and $k_a+k_b=k_5$. But then $0=r_1+r_2=k_1=k_5=k_a+k_b$, so $k_a=k_b=0$. 
But $k_a+2k_b=n-\oooo{n}\geq 1$, a contradiction.

Second case, $(\oooo{H_\Z},\oooo{I^{(0)}})$ of type
$D_{\oooo{n}}$: Recall Theorem \ref{t3.6} (a).
Then $U$ is of type 
$\sum_{i=1}^{r_1}A_{c_i}+\sum_{j=1}^{r_3}D_{b_j}$ with 
$r_1,r_2,r_3\in\Z_{\geq 0}$, $r_1+r_3\geq 1$, 
$c_i\in\N$, $b_i\in\Z_{\geq 2}$,
$$\oooo{n}=r_1+r_2+\sum_{i=1}^{r_1}c_i+\sum_{j=1}^{r_3}b_j.$$
Table \ref{tab3.1} gives 
$$k_a=\oooo{n}-l(\oooo{-M})=\oooo{n}-\dim_\Q U_\Q=r_1+r_2.$$ 
Split $r_3=r_4+r_5+r_6$ with 
\begin{eqnarray*}
r_4&:=&|\{j\,|\, b_j=2\}|,\\
r_5&:=&|\{j\,|\, b_j\geq 3,
\oooo{-M}|_{(\textup{the component from }j)}
\textup{ is a Coxeter element}\}|,\\
r_6&=&|\{j\,|\, b_j\geq 3,
\oooo{-M}|_{(\textup{the component from }j)}
\textup{ is a quasi-Coxeter element},\\
&&\hspace*{6cm}\textup{but not a Coxeter element}\}|.
\end{eqnarray*}
Then
\begin{eqnarray*}
\tr(-M)&\geq&(r_1+2r_4+r_5)(-1)+2(k_a+k_b)\\
&\geq& -r_1-2r_4-r_5+2k_5(\oooo{H_\Z},\oooo{-M})\\
&\stackrel{\textup{Theorem \ref{t3.11}}}{=}&
-r_1-2r_4-r_5+2k_1\\
&\stackrel{\textup{Table \ref{tab3.1}}}{=}&
-r_1-2r_4-r_5+\left\{\begin{array}{ll}
2(\sum_{i=1}^3r_i-1)&\textup{ if }r_3\geq 1\\
2(r_1+r_2)&\textup{ if }r_3=0.\end{array}\right.\\
&=& \left\{\begin{array}{ll}
r_1+2r_2+r_5+2r_6-2 &\textup{ if }r_3\geq 1,\\
r_1+2r_2&\textup{ if }r_3=0\end{array}\right. 
\end{eqnarray*}
This is $\leq -1$ only if $r_3\geq 1$, $r_2=r_6=0$, 
$r_1+r_5\leq 1$ and $k_a+k_b=k_5=k_1=r_1+r_3-1$. Also, 
$k_a+2k_b=n-\oooo{n}\geq 1$ implies $k_a+k_b\geq 1$,
so $r_1+r_3-1\geq 1$. Now $k_a=r_1$ implies $k_b=r_3-1$. 
The following possibilities remain,
\begin{eqnarray*}
\begin{array}{l|l|l|l|l|l|l}
r_1 & r_5 & \textup{type of }U & \oooo{n} & k_a & k_b & 
\tr(-M)\\
\hline 
0 & 0 & r_3D_2 \textup{ with }r_3\geq 2& 2r_3 & 0 & r_3-1 & -2\\
0 & 1 & r_4D_2+D_b \textup{ with }b\geq 3,r_4\geq 1& 2r_4+b
& 0 & r_4 & -1\\
1 & 0 & A_c+r_3D_2\textup{ with }r_3\geq 1 & c+1+2r_3 & 1 
& r_3-1 & -1
\end{array}
\end{eqnarray*}

Third case, $(\oooo{H_\Z},\oooo{I^{(0)}})$ of type
$E_{\oooo{n}}$ with $\oooo{n}\in\{6,7,8\}$:
Recall Theorem \ref{t3.6} (b). Denote by $r_5+r_6$
the number of orthogonal irreducible simply laced root lattices
into which $U$ decomposes, where $r_5$ denotes the number of 
those on which $\oooo{-M}$ restricts to a Coxeter element,
and $r_6$ denotes the number of those on which $\oooo{-M}$
restricts to a quasi-Coxeter element which is not a Coxeter
element. Then 
\begin{eqnarray*}
\tr(-M)&\geq& r_5(-1)+2(k_a+k_b)\\
&\geq& -r_5+2k_5(\oooo{H_\Z},\oooo{-M})\\
&\stackrel{(1)}{\geq}& -r_5+2k_1\geq -(r_5+r_6)+2k_1.
\end{eqnarray*}
Here $\stackrel{(1)}{\geq}$ is almost always an equality.
The exceptions are made explicit in Theorem \ref{t3.13}.

The case $k_1=0$ is impossible
because always $k_a+k_b\geq 1$, and $k_1=0$ implies $r_5\leq 1$, 
so then $\tr(-M)\geq -1+2\cdot 1=1$.

Now compare the Tables \ref{tab3.2}, \ref{tab3.3} and 
\ref{tab3.4}. The inequalities $r_5+r_6>2k_1$ and $k_1\geq 1$
hold only in few cases, namely in the following cases:
%\begin{table}[H]
\begin{eqnarray*}
\begin{array}{l|l|l|l}
(\oooo{H_\Z},\oooo{I^{(0)}}) & U & r_5+r_6 & k_1 \\ \hline
E_6 & 3A_2 & 3 & 1 \\
    & 2A_2+A_1 & 3 & 1 \\ \hline
E_7 & 2A_3+A_1 & 3 & 1 \\
    & 7A_1 & 7 & 3 \\
    & A_3+A_2+A_1 & 3 & 1 \\ \hline
E_8 & A_5+A_2+A_1 & 3 & 1 \\
    & A_4+A_2+A_1 & 3 & 1 
\end{array}
\end{eqnarray*}
%\caption[Table 6.2]{..}
%\label{tab7.1}
%\end{table}
In these cases $k_5=k_1$ holds by Theorem \ref{t3.13}, 
and $\stackrel{(1)}{\geq}$ is an equality. 
$\tr(-M)\leq -1$ implies furthermore
$k_a+k_b=k_5=k_1$, $r_6=0$ and $\tr(-M)=-1$. 
The part of Table \ref{tab7.1}, which concerns
$(\oooo{H_\Z},\oooo{I^{(0)}})$ of type $E_6$, $E_7$
or $E_8$, follows.\hfill$\Box$

\bigskip
Now we will study for which matrices their spectra
satisfy the variance inequality. 
Theorem \ref{t7.3} will do the first step.
It will tell that the tuple $(\oooo{H_\Z},\oooo{I^{(0)}},
\oooo{-M})$ is in one of four families. 
Theorem \ref{t7.4} will do the second step.
It will tell that there are only four
$\Br_n\ltimes\{\pm 1\}^n$ orbits of such matrices, 
those from the three (families of) simple elliptic singularities 
and one from the tubular elliptic root system $D_4^{(1,1)}$. 
But the matrices $S$ from $D_4^{(1,1)}$ satisfy 
$tr(S^{-1}S^t)=2$.

So the trace condition $\tr(S^{-1}S^t)=1$ and the variance
inequality together are strong enough to single out within the
matrices in Lemma \ref{t7.1} those which come from the
simple elliptic singularities.

\begin{theorem}\label{t7.3}
Let $S\in T_n(\Z)$ and $(H_\Z,L,\uuuu{e})$ be as in
Lemma \ref{t7.1}. 
The notations and results in Lemma \ref{t7.1} are used, 
especially $(\oooo{H_\Z},\oooo{I^{(0)}})$ is an irreducible
simply laced root lattice of rank $\oooo{n}=l(\oooo{-M})+k_a$.
Suppose additionally that the variance inequality 
\begin{eqnarray}\label{7.3}
\Var(\Sp(S))\leq\frac{\alpha_n-\alpha_1}{12}
\end{eqnarray}
is satisfied. 

Let $U\subset\oooo{H_\Z}$ be the subroot lattice associated
to a reduced presentation of $\oooo{-M}$.
Table \ref{tab7.2} lists information on all possibilities
for $((\oooo{H_\Z},\oooo{I^{(0)}}),\oooo{-M},U,k_a,k_b)$.

\begin{table}
\begin{eqnarray*}
\begin{array}{l|l|l|l|l}
{}(\oooo{H_\Z},\oooo{I^{(0)}}) 
& U & k_a & k_b & k_a+k_b \\ 
\hline
D_{4} & 2D_2  & 0 & 1 & 1 \\ 
E_6 & 3A_2 & 0 & 1 & 1  \\ 
E_7 & 2A_3+A_1 %\sim D_4(a_1)+3A_1
& 0 & 1 & 1\\ 
E_8 & A_5+A_2+A_1 & 0 & 1 & 1 
\end{array}
\end{eqnarray*}
\caption[Table 7.2]{Possibilities for 
$((\oooo{H_\Z},\oooo{I^{(0)}}),\oooo{-M},U,k_a,k_b)$
in Theorem \ref{t7.3}}
\label{tab7.2}
\end{table}
\end{theorem}

{\bf Proof:} 
The first $k_a+k_b$ spectral numbers are equal to 
$-\frac{1}{2}$, the last $k_a+k_b$ spectral numbers are
equal to $\frac{1}{2}$. Let $g$ be the characteristic polynomial
of $\oooo{-M}|_{\oooo{H_{\neq 1,\Q}}}$. 
Recall Definition \ref{t6.3} of $V_2(g)$. Then 
\begin{eqnarray*}
12\cdot V_2(\Sp(S))
&=&12\cdot V_2(g)+12\cdot 2(k_a+k_b)\cdot \frac{1}{4}\\
&=&12\cdot V_2(g) + 6(k_a+k_b),\\
n(\alpha_n-\alpha_1)
&=& n = \oooo{n}+k_a+2k_b.
\end{eqnarray*}
Therefore the variance inequality \eqref{7.3} is equivalent to
the inequality  
\begin{eqnarray}\label{7.4}
12\cdot V_2(g)&\leq& \oooo{n}-5k_a-4k_b.
\end{eqnarray}
Let $U\subset\oooo{H_\Z}$ be the subroot lattice associated
to a reduced presentation of $\oooo{-M}$. Then
$$\dim U_\Q=l(\oooo{-M})=\deg g=\oooo{n}-k_a.$$

First case, $(\oooo{H_\Z},\oooo{I^{(0)}})$ of type $A_{\oooo{n}}$:
Recall Theorem \ref{t3.6} (a). Then $U$ is of type
$\sum_{i=1}^{r_1}A_{c_i}$ with $r_1\in\N$, $r_2\in\Z_{\geq -1}$,
$c_i\in\N$, $\oooo{n}=r_1+r_2+\sum_{i=1}^{r_1}c_i$,
$k_a=r_1+r_2$. 
Observe with the help of Table \ref{tab6.1} 
\begin{eqnarray}
12\cdot V_2(\Sp(A_{c_i}))
&=&\frac{c_i(c_i-1)}{c_i+1}=c_i-2+\frac{2}{c_i+1}.\label{7.5}
\end{eqnarray}
Therefore 
\begin{eqnarray*}
12\cdot V_2(g)&=& 
\sum_{i=1}^{r_1}\Bigl(c_i-2+\frac{2}{c_i+1}\Bigr)\\
&=& \oooo{n}-3r_1-r_2+\sum_{i=1}^{r_1}\frac{2}{c_i+1}.
\end{eqnarray*}
The inequality \eqref{7.4} is equivalent to the inequality

\begin{eqnarray*}
5k_a+4k_b+\sum_{i=1}^{r_1}\frac{2}{c_i+1}
&\leq& 3r_1+r_2=3k_a-2r_2,\\
\textup{so}\quad 
2k_a+4k_b+\sum_{i=1}^{r_1}\frac{2}{c_i+1}\leq -2r_2,
\end{eqnarray*}
which is impossible. So, $(\oooo{H_\Z},\oooo{I^{(0)}} )$
being of type $A_{\oooo{n}}$ is excluded. 

Second case, $(\oooo{H_\Z},\oooo{I^{(0)}})$ of type
$D_{\oooo{n}}$: Recall Theorem \ref{t3.6} (a).
Then the restriction $\oooo{-M}|_U$ is a quasi-Coxeter element
of type 
$\sum_{i=1}^{r_1}A_{c_i}+\sum_{i=1}^{r_3}D_{b_i}(a_{j_i})$ with 
$r_1,r_2,r_3\in\Z_{\geq 0}$, $r_1+r_3\geq 1$, 
$c_i\in\N$, $b_i\in\Z_{\geq 2}$, 
${j_i}=0$ if $b_i\in\{2;3\}$ and 
${j_i}\in\{0,...,\lfloor\frac{b_i}{2}\rfloor-1\}$ if $b_i\geq 4$
(and $D_b(a_0):=D_b$),  
$$\oooo{n}=r_1+r_2+\sum_{i=1}^{r_1}c_i+\sum_{i=1}^{r_3}b_i.$$
Table \ref{tab3.1} gives 
$$k_a=\oooo{n}-l(\oooo{-M})=\oooo{n}-\dim_\Q U_\Q=r_1+r_2.$$ 
Split $r_3=r_4+r_5+r_6$ as in the proof of Theorem \ref{t7.2}. 
Again \eqref{7.5} will be used. Table \ref{tab6.1} gives
for ${j_i}\in\{0,...,\lfloor\frac{b}{2}\rfloor-1\}$  
\begin{eqnarray}
12\cdot V_2(\Sp(D_{b}(a_j)))
&=& b-\frac{1}{b-j-1}-\frac{1}{j+1}.
\label{7.6}
\end{eqnarray}
Therefore 
\begin{eqnarray*}
12\cdot V_2(g)&=& \sum_{i=1}^{r_1}\Bigl(c_i-2+\frac{2}{c_i+1}\Bigr)
+\sum_{i=1}^{r_3}\Bigl(b_i-\frac{1}{b_i-{j_i}-1}-\frac{1}{{j_i}+1}\Bigr)\\
&=& \oooo{n}-3k_a+2r_2+\sum_{i=1}^{r_1}\frac{2}{c_i+1}
-\sum_{i=1}^{r_3}\Bigl(\frac{1}{b_i-{j_i}-1}+\frac{1}{{j_i}+1}\Bigr).
\end{eqnarray*}
The inequality \eqref{7.4} is equivalent to the inequality
\begin{eqnarray}\label{7.7}
2k_a+4k_b+2r_2+\sum_{i=1}^{r_1}\frac{2}{c_i+1}
&\leq& \sum_{i=1}^{r_3}\Bigl(\frac{1}{b_i-{j_i}-1}+\frac{1}{{j_i}+1}\Bigr).
\end{eqnarray}
The case $r_3=0$ is impossible, so $r_3\geq 1$. Also
\begin{eqnarray*}
k_a+k_b&\geq& k_5(\oooo{H_\Z},\oooo{-M}) = k_1 = 
\sum_{i=1}^3r_i-1=k_a+r_3-1,\\
\textup{so }k_b&\geq& r_3-1.
\end{eqnarray*}
One extracts from \eqref{7.7} the coarser inequality 
$$4(r_3-1)\leq 2r_3,$$
which implies $r_3\in\{1;2\}$. The case $r_3=1$ 
requires $k_b=0$, $k_a=1$, $r_1=r_2=0$ and thus  
$k_a=0$, which is impossible. 
The case $r_3=2$ requires $k_b=1$, $k_a=r_2=r_1=0$,
$D_{b_i}(a_{j_i})=D_2(a_0)=D_2$ for $i\in\{1;2\}$. 

Third case, $(\oooo{H_\Z},\oooo{I^{(0)}})$ of type
$E_{\oooo{n}}$ with $\oooo{n}\in\{6,7,8\}$:
The inequality \eqref{7.4} requires that its right hand side
$\oooo{n}-5k_a-4k_b$ is non-negative. The only possibilities
for $(k_a,k_b)$ are $(0,1)$, $(1,0)$ and in the case
$E_8$ also $(0,2)$. Though this last case would only work
if $V_2(g)=0$, which would imply that $\oooo{-M}$
is a quasi-Coxeter element of type $8A_1$. 
But then $0+2=k_a+k_b\geq k_5=k_1=4$, which is impossible.
So, only the two cases $(k_a,k_b)\in\{(0;1),(1;0)\}$ remain.
The inequality $k_a+k_b\geq k_5\stackrel{\textup{almost always}}{=}k_1$ 
implies that we can restrict to the cases with $k_1\in\{0;1\}$. 

The following Tables \ref{tab7.3}, \ref{tab7.4} and \ref{tab7.5}
go through the parts of the Tables
\ref{tab3.2}, \ref{tab3.3} and \ref{tab3.4} which concern
the cases with $k_1\in\{0;1\}$. 
Here the information on $U$ has to be refined
to the information on the type of $\oooo{-M}|_U$ as
quasi-Coxeter transformation on $U$. 

The tables show and compare the values
$$\www{L}:=12\cdot V_2(g)\quad\textup{and}\quad 
\www{R}:=\oooo{n}-5k_a-4k_b.$$
The values $\www{L}$ can be determined easily with Table 
\ref{tab6.1}
and the additivity of $V_2$ in Lemma \ref{t6.4} (b). 
The tables show that the variance inequality and the equivalent 
inequality \eqref{7.4}, which is $\www{L}\leq \www{R}$, 
are satisfied only in the cases in Table
\ref{tab7.2}, there with equality.\hfill$\Box$ 

\begin{table}
\begin{eqnarray*}
\begin{array}{l|c|r}
\oooo{-M}|_U & (k_a;k_b) & \www{L}\stackrel{?}{\leq}\www{R} 
\\ \hline
E_6 & (0;1) & 5>2 \\
E_6(a_1) & (0;1) & \frac{50}{9}>2 \\
E_6(a_2) & (0;1) & 6>2 \\
A_5+A_1 & (0;1)& \frac{10}{3}>2 \\
3A_2 & (0;1) & 2 = 2 
\end{array}
\hspace*{0.5cm}
\begin{array}{l|c|r}
\oooo{-M}|_U & (k_a;k_b) & \www{L}\stackrel{?}{\leq}\www{R}
\\ \hline
A_5 & (1;0) & \frac{10}{3}>1 \\
2A_2+A_1 & (1;0) & \frac{4}{3}>1\\
A_4 + A_1 & (1;0) & \frac{12}{5}>1 \\
D_5 & (1;0) & \frac{15}{4}>1\\ 
D_5(a_1) & (1;0) & \frac{25}{6}>1
\end{array}
\end{eqnarray*}
\caption[Table 7.3]{$(\oooo{H_\Z},\oooo{I^{(0)}})$ of type $E_6$:
variance of the spectrum of $\oooo{-M}|_U$ with 
$k_1\in\{0;1\}$}
\label{tab7.3}
\end{table}

\begin{table}
\begin{eqnarray*}
\begin{array}{l|c|r}
\oooo{-M}|_U & (k_a;k_b) & \www{L}\stackrel{?}{\leq}\www{R}
\\ \hline
E_7 & (0;1) & \frac{56}{9}>3 \\
E_7(a_1) & (0;1) & \frac{48}{7}>3 \\
E_7(a_2) & (0;1) & 7>3 \\
E_7(a_3) & (0;1) & \frac{112}{15}>3 \\
E_7(a_4) & (0;1) & 8>3 \\
D_6+A_1 & (0;1) & \frac{24}{5}>3 \\
D_6(a_1)+A_1 & (0;1) & \frac{21}{4}>3 \\
D_6(a_2)+A_1 & (0;1) & \frac{16}{3}>3 \\
A_5+A_2 & (0;1) & 4>3 \\
2A_3+A_1 & (0;1) & 3=3 \\
A_7 & (0;1) & \frac{21}{4}>3 
\end{array}
\hspace*{0.5cm}
\begin{array}{l|c|r}
\oooo{-M}|_U & (k_a;k_b) & L\stackrel{?}{\leq}R \\ \hline
E_6 & (1;0) & 5>2 \\
E_6(a_1) & (1;0) & \frac{50}{9}>2 \\
E_6(a_2) & (1;0) & 6>2 \\
D_5+A_1 & (1;0) & \frac{15}{4}>2 \\
D_5(a_1)+A_1 & (1;0) & \frac{25}{6}>2 \\
A_4+A_2 & (1;0) & \frac{46}{15}>2 \\
A_3\! +\! A_2\! +\! A_1 & (1;0) & \frac{13}{6}>2 \\
{}[A_5+A_1]'' & (1;0) & \frac{10}{3}>2 \\
D_6 & (1;0) & \frac{24}{5}>2 \\
D_6(a_1) & (1;0) & \frac{21}{4}>2 \\
D_6(a_2) & (1;0) & \frac{16}{3}>2 \\
A_6 & (1;0) & \frac{30}{7}>2 
\end{array}
\end{eqnarray*}
\caption[Table 7.4]{$(\oooo{H_\Z},\oooo{I^{(0)}})$ of type $E_7$:
variance of the spectrum of $\oooo{-M}|_U$ with $k_1\in\{0;1\}$}
\label{tab7.4}
\end{table}

\begin{table}
\begin{eqnarray*}
\begin{array}{l|c|r}
\oooo{-M}|_U & (k_a;k_b) & L\stackrel{?}{\leq}R \\ \hline
E_8 & (0;1) & \frac{112}{15}>4 \\
E_8(a_1) & (0;1) & \frac{49}{6}>4 \\
E_8(a_2) & (0;1) & \frac{42}{5}>4 \\
E_8(a_3) & (0;1) & \frac{61}{6}>4 \\
E_8(a_4) & (0;1) & \frac{80}{9}>4 \\
E_8(a_5) & (0;1) & \frac{136}{15}>4 \\
E_8(a_6) & (0;1) & \frac{48}{5}>4 \\
E_8(a_7) & (0;1) & \frac{29}{3}>4 \\
E_8(a_8) & (0;1) & 8>4 \\
A_8 & (0;1) & \frac{56}{9}>4 \\
D_8 & (0;1) & \frac{48}{7}>4 \\
D_8(a_1) & (0;1) & \frac{22}{3}>4 \\
D_8(a_2) & (0;1) & \frac{112}{15}>4 \\
D_8(a_3) & (0;1) & \frac{15}{2}>4 \\
A_7+A_1 & (0;1) & \frac{21}{4}>4 \\
A_5\! +\! A_2\! +\! A_1 & (0;1) & 4=4 \\
2A_4 & (0;1) & \frac{24}{5}>4 \\
E_6+A_2 & (0;1) & \frac{17}{3}>4 \\
E_6(a_1)+A_2 & (0;1) & \frac{56}{9}>4 \\
E_6(a_2)+A_2 & (0;1) & \frac{20}{3}>4 \\
E_7+A_1 & (0;1) & \frac{56}{9}>4 \\
E_7(a_1)+A_1 & (0;1) & \frac{48}{7}>4 
\end{array}
\hspace*{0.5cm}
\begin{array}{l|c|r}
\oooo{-M}|_U & (k_a;k_b) & L\stackrel{?}{\leq}R \\ \hline
E_7(a_2)+A_1 & (0;1) & 7>4 \\
E_7(a_3)+A_1 & (0;1) & \frac{112}{15}>4 \\
E_7(a_4)+A_1 & (0;1) & 8>4 \\
D_5+A_3 & (0;1) & \frac{21}{4}>4\\
D_5(a_1)+A_3 & (0;1) & \frac{17}{3}>4 \\ \hline
A_6+A_1 & (1;0) & \frac{30}{7}>3\\
A_4\! +\! A_2\! +\! A_1 & (1;0) & \frac{46}{15}>3\\
E_6+A_1 & (1;0) & 5>3\\
E_6(a_1)+A_1 & (1;0) & \frac{50}{9}>3\\
E_6(a_2)+A_1 & (1;0) & 6>3\\
E_7 & (1;0) & \frac{56}{9}>3\\
E_7(a_1) & (1;0) & \frac{48}{7}>3\\
E_7(a_2) & (1;0) & 7>3\\
E_7(a_3) & (1;0) & \frac{112}{15}>3\\
E_7(a_4) & (1;0) & 8>3\\
D_7 & (1;0) & \frac{35}{6}>3\\
D_7(a_1) & (1;0) & \frac{63}{10}>3\\
D_7(a_2) & (1;0) & \frac{77}{12}>3\\
D_5+A_2 & (1;0) & \frac{53}{12}>3\\
D_5(a_1)+A_2 & (1;0) & \frac{29}{6}>3\\
A_4+A_3 & (1;0) & \frac{39}{10}>3\\
{}[A_7]' & (1;0) & \frac{21}{4}>3
\end{array}
\end{eqnarray*}
\caption[Table 7.5]{$(\oooo{H_\Z},\oooo{I^{(0)}})$ of type $E_8$:
variance of the spectrum of $\oooo{-M}|_U$ with 
$k_1\in\{0;1\}$}
\label{tab7.5}
\end{table}

\begin{theorem}\label{t7.4}
Let $S\in T_n(\Z)$ be as in Theorem \ref{t7.3},
so $S+S^t$ is positive semidefinite, ${\rm CDD}(S)$ is connected,
and the variance inequality \eqref{7.3} is satisfied.
Then $S$ is either a distinguished matrix of a simple elliptic
singularity or it is in the $\Br_6\ltimes\{\pm 1\}^6$
orbit which is associated to Coxeter elements in the 
elliptic root system $D_4^{(1,1)}$. 
In the second case $\tr(S^{-1}S^t)=2$ (and not $=1$).
\end{theorem}

{\bf Proof:} Let $(H_\Z,L,\uuuu{e})$ be the unimodular
bilinear lattice with triangular basis which is associated to $S$
with $S=L(\uuuu{e}^t,\uuuu{e})$. 

We apply Theorem \ref{t7.3}. Table \ref{t7.2} gives
$\rk \Rad(I^{(0)})=2k_b=2$. Table \ref{t7.2} and 
Theorem \ref{t4.3} tell especially that 
$(H_\Z,I^{(0)})$ is the root lattice of a tubular elliptic root system
of type $D_4^{(1,1)},E_6^{(1,1)},E_7^{(1,1)}$ or 
$E_8^{(1,1)}$. 

Minus the monodromy is
$-M=s_{e_1}\circ...\circ s_{e_n}$, and the induced automorphism
$\oooo{-M}$ on $\oooo{H_\Z}$ is by Table \ref{t7.2} 
a Coxeter element in a subroot lattice of type 
$2D_2$, $3A_2$, $2A_3+A_1$ respectively $A_5+A_2+A_1$. 
Theorem \ref{t4.9} says that $-M$ is a Coxeter element in
the tubular elliptic root system. 

We see that $S$ arises from a tubular elliptic root system
$(H_\Z,I,\Phi)$ and a Coxeter element $c$ in it as the
upper triangular part $S$ of a matrix 
$I(\uuuu{\alpha}^t,\uuuu{\alpha})=S+S^t$ for some
$\uuuu\alpha$ in the set in \eqref{4.1}. 

Theorem \ref{t4.11} says that in the cases of 
$E_k^{(1,1)}$ with $k\in\{6,7,8\}$ all possible
tuples $\uuuu\alpha$ in the set in \eqref{4.1} form a
single $\Br_n\ltimes\{\pm 1\}^n$ orbit. Therefore also
the associated matrices form a single 
$\Br_n\ltimes\{\pm 1\}^n $ orbit. 
By Theorem \ref{t5.5} (g) it is the orbit of distinguished
matrices of a simple elliptic singularity of type $\www{E_k}$.

In the case of $D_4^{(1,1)}$, Theorem \ref{t4.11} says that 
the set in \eqref{4.1} decomposes into two
$\Br_6\ltimes\{\pm 1\}^6$ orbits of tuples $\uuuu{\alpha}$. 
But Example \ref{t4.13} gave a tuple $\uuuu{e}$ in one orbit
and a tuple $g(\uuuu{e})$ in the other orbit with the same
matrix
$$I(\uuuu{e}^t,\uuuu{e})=I(g(\uuuu{e})^t,g(\uuuu{e})).$$
Therefore also the upper triangular part $S\in T_6^{uni}(\Z)$
of the decomposition $S+S^t$ of this matrix is the same for 
$\uuuu{e}$ and for $g(\uuuu{e})$. 
Therefore the two orbits in the set in \eqref{4.1} give the
same $\Br_6\ltimes\{\pm 1\}^6$ orbit of matrices. 
So, on the level of matrices, we have only one orbit and not two.

Table \ref{tab7.1} shows that all these matrices satisfy 
$\tr(S^{-1}S^t)=2$.
\hfill$\Box$

\section{Conjectures for the indefinite cases}\label{s8}
\setcounter{equation}{0}
\setcounter{table}{0}

The main result Theorem \ref{t1.5} treats the cases 
$S\in T_n(\Z)$ with connected CDD where $S+S^t$ is positive 
(semi)definite.
In such a case Lemma \ref{t2.11} (c) provided a spectrum
in an easy way. In the case of $S+S^t$ indefinite, we rely 
on the idea of Cecotti and Vafa \cite{CV93} how to associate
a spectrum $\Sp(S)=(\alpha_1,...,\alpha_n)\in\R^n$ to a 
matrix $S\in T_n^{S^1}(\R)$. This idea is as follows.

\begin{conjecture}[{\cite[pages 583, 589, 590]{CV93}}]\label{t8.1}
 Start with a matrix $S_1\in T_n^{S^1}(\R)$. 
Conjecturally, the problems in the following two steps
(see the Remarks \ref{t8.2}) can be overcome.

{\bf Step 1:}
Choose in a natural way a continuous path 
$S_\bullet:[0,1]\to T_n^{S^1}(\R)$ 
from $S_0=E_n$ (=the unit matrix) to $S_1$.

{\bf Step 2:}
Choose {\it in a natural way} $n$ continuous maps
$\alpha_{j,\bullet}:[0;1]\to \R,$ $j=1,...,n$, such that 
$\alpha_{j,0}=0$ and $e^{-2\pi i\alpha_{1,r}},...,
e^{-2\pi i\alpha_{n,r}}$ 
are the eigenvalues of $S_r^{-1}S_r^t$ for $r\in[0;1]$.

Then 
$$\Sp(S_1):=(\alpha_{1,1},...,\alpha_{n,1}).$$
\end{conjecture}

\begin{remarks}\label{t8.2}
(i) Problems in Step 1: Is $T_n^{S^1}(\R)$ connected?
If yes, is the choice of the path $S_\bullet$ relevant?
Is $T_n^{S^1}(\R)$ simply connected? Is it contractible?
(Yes in the cases $n=2$ and $n=3$ \cite{BaH20}.)

(ii) Problem in Step 2:
If the path $S_\bullet$ is chosen, how shall one choose 
$\alpha_{j,\bullet}$ near an $r\in ]0,n[$
where $S_r^{-1}S_r^t$ has multiple eigenvalues?

(iii) In the positive (semi)definite cases,
the spectral numbers are in the interval $[-\frac{1}{2},
\frac{1}{2}]$. Therefore we neither need to consider
a path as in Step 1 nor do we have a non-uniqueness problem
as in Step 2. Only the spectral numbers associated to the
eigenvalue $-1$ are not unique, and they are distributed equally
at $\frac{1}{2}$ and $-\frac{1}{2}$.
\end{remarks}

We have a solution \cite{BaH20} of this conjecture for
matrices in certain subsets of $T_n^{S^1}(\R)$. 
The solution is described in Theorem \ref{t8.4} (d).
Definition \ref{t8.3} is needed. 
The solution combines two ideas. 
The first idea is that a subspace
of $T_n^{S^1}(\R)$ of sufficiently good matrices $S$ is
contractible, so paths $S_\bullet$ can be chosen, 
and the precise choice is  irrelevant. 
This solves the problems in Step 1.
The second idea is that for sufficiently good $S_r$ the monodromy
matrix $S_r^{-1}S_r^t$ (or $-S_r^{-1}S_r^t$) 
has a natural $n$-th root, and one can define some 
pseudo-spectral numbers in the interval $[0;1]$ 
from the eigenvalues of this $n$-th root. That they are in the
interval $[0;1]$ solves the problem in step 2. 
Then they induce the spectral numbers $\alpha_{1,r},...,
\alpha_{n,r}$.

\begin{definition}[{\cite[3.1 \& 3.2 \& 4.4]{BaH20}}]\label{t8.3} 
Fix $n\in\N$ and $b\in\{1,2\}$. 

\medskip
(a) Define a subspace 
$T_{HORb}^{pol}(n,\R)\subset\R[x]_{\deg=n}$ by
\begin{eqnarray*}
T_{HORb}^{pol}(n,\R)&:=& \{p(x)=p_nx^n+p_{n-1}x^{n-1}+...+p_1x+p_0\,| \\
&& \hspace*{0.5cm}p_n=1,\quad p_{n-j}=(-1)^{b-1}p_j,\\
&& \hspace*{0.5cm}\textup{all zeros of }p\textup{ are in }S^1\}\\
&\stackrel{!}{=}& 
\{p(x)=\prod_{j=1}^n(x-e^{-2\pi i\beta_j})\, |\, 
0\leq\beta_1\leq...\leq \beta_n\leq 1, \\
&& \hspace*{0.5cm}0=\beta_1\textup{ if }b=2, \quad \beta_j+\beta_{n+b-j}=1\}.
\end{eqnarray*}
$p(x)$ and $(\beta_1,...,\beta_n)$ determine one another.

(b) For $p\in T_{HORb}^{pol}(n,\R)$ define
\begin{eqnarray*}
R_{(b)}(p)&:=& \left(\begin{array}{llll|l}
-p_{n-1} & -p_{n-2} & ... & -p_1 & -p_0 \\ \hline
 & & & & \\
 & &  E_{n-1}& & \\
 & & & & 
 \end{array}\right)
\end{eqnarray*}
It is a companion matrix. Its characteristic polynomial is $p$.
Its eigenvalues are $e^{-2\pi i\beta_j}$. It has for each 
eigenvalue one Jordan block.

(c) For $p\in T_{HORb}^{pol}(n,\R)$ define
\begin{eqnarray*}
S_{(b)}(p)&:=& \left(\begin{array}{lllll}
1 & p_{n-1} & ... & p_2 & p_1 \\
 & \ddots & \ddots & & p_2 \\
 & & \ddots & \ddots & \vdots \\
 & & & \ddots & p_{n-1} \\
 & & & & 1  \end{array}\right)
\end{eqnarray*}
Define 
$$T_{HORb}(n,\R):=S_{(b)}(T_{HORb}^{pol}(n,\R))
\stackrel{\eqref{8.1}}{\subset}T_n^{S^1}(\R).$$
\end{definition}

Part (b) of Theorem \ref{t8.4} says that 
$(-1)^b$ times the monodromy matrix of $S_{(b)}(p)$ 
has a natural $n$-th root, the matrix $R_{(b)}(p)$. 

\begin{theorem}[{\cite[3.2--3.5 \& 4.5]{BaH20}}]\label{t8.4}
Fix $n\in\N$ and $b\in\{1,2\}$.

(a) $T_{HORb}^{pol}(n,\R)$ is a simplex of dimension
$$\begin{array}{lll}
 &  \dim T_{HOR1}^{pol}(n,\R) & \dim T_{HOR2}^{pol}(n,\R) \\
n\textup{ odd } & \frac{n-1}{2} & \frac{n-1}{2} \\
n\textup{ even } & \frac{n}{2} & \frac{n-2}{2} 
 \end{array}$$
For $p$ in the interior of this simplex, the $\beta_j$ are pairwise different.

(b) (Also \cite{Ho17}) For any $p\in T_{HORb}^{pol}(n,\R)$
\begin{eqnarray}\label{8.1}
(-1)^b\cdot S_{(b)}(p)^{-1}S_{(b)}(p)^t &=& R_{(b)}(p)^n.
\end{eqnarray}

(c) $E_n= S_{(b)}(p_0)\in\textup{int}(T_{HORb}(n,\R))$ with 
$$p_0=x^n-(-1)^b,\quad\textup{so}\quad
\beta_{j,0}=\frac{j-\frac{b}{2}}{n}\in [0,1]\textup{ for }
j=1,...,n.$$

(d) (i) The recipe of Cecotti and Vafa works for 
$S_1 =  S_{(b)}(p_1)\in T_{HORb}(n,\R)$ as follows.
Let $(\beta_{1,1},...,\beta_{n,1})$ be associated to $S_1$
as in Definition \ref{t8.3}. Then the spectrum of $S_1$
consists of
\begin{eqnarray}\label{8.2}
\www{\alpha}_{j,1} &=& n\beta_j - (j -\frac{b}{2})
=n(\beta_{j,1}-\beta_{j,0})\quad
\textup{for}\quad j=1,...,n.
\end{eqnarray}

(ii) As $T_{HORb}(n,\R)$ is a simplex, a path 
$S_\bullet:[0,1]\to T_{HORb}(n,\R)$ from $E_n=S_0$ to 
$S_1$ can be chosen, and the precise choice is irrelevant.

The choice of the paths $\www{\alpha}_{j,\bullet}:[0,1]\to\R$ 
is also no problem.
Now $\www{\alpha}_{j,r}$ comes from the eigenvalue 
$e^{-2\pi i\beta_{j,r}}$ of $R_{(b)}(r)$ by formula \eqref{8.2},
and these eigenvalues are pairwise different except at the 
boundary of the simplex $T_{HORb}(n,\R)$.

(iii) The $\www{\alpha}_{j,1}$ are not indexed by size 
(in general), but they satisfy 
\begin{eqnarray*}
\www{\alpha}_{j,1}+\www{\alpha}_{n+1-j,1}
&=& 0 \quad\textup{if }b=1,\\
\www{\alpha}_1=0,\ \www{\alpha}_j+\www{\alpha}_{n+2-j,1}
&=&0\quad\textup{if }b=2.
\end{eqnarray*}
We reorder them by size and call the resulting numbers
$\alpha_1,...,\alpha_n$, so $\Sp(S_1)=(\alpha_1,...,\alpha_n)$
with $\alpha_1\leq ...\leq \alpha_n$. 

(iv) An unordered tuple of $n$ numbers arises as $\Sp(S_1)$ for 
$S_1\in T_{HORb}(n,\R)$ if and only if the numbers can be ordered
as $\www{\alpha}_1,...,\www{\alpha}_n$ such that the symmetries 
in (iii) hold and $\www{\alpha}_{j+1}\geq\www{\alpha}_j-1$, 
and in the case $b=1$ also $\www{\alpha}_1\geq -\frac{1}{2}$.
\end{theorem}

A family of positive cases where everything fits together
are the {\it chain type singularities},
see Definition \ref{t8.5}, Theorem \ref{t8.6}, Theorem 
\ref{t8.7} and the Remarks \ref{t8.8}.

\begin{definition}\label{t8.5}
A {\it chain type singularity} is a quasihomogeneous
singularity of the following shape, 
$$f(x_0,...,x_m)=x_0^{a_0}+\sum_{i=1}^m x_{i-1}x_i^{a_i},$$
where $(a_0,a_1,...,a_m)\in\Z_{\geq 2}\times \N^m$. 
Define $r_{-1}:=1$ and $r_k:=a_0a_1...a_k$ for 
$k=0,...,m$, and define $b\in\{1;2\}$ by $b\equiv m+1\mmod 2$.
\end{definition}

The following theorem was conjectured by Orlik and Randell 
\cite{OR77} and proved recently by Varolgunes \cite{Va23}.

\begin{theorem}[{\cite{Va23}}]\label{t8.6}
Consider a chain type singularity $f$ with Milnor number $n$ and $a_0,...,a_m,
r_{-1},r_0,...,r_m$ and $b$ as in Definition
\ref{t8.5} and consider the polynomial
$$p:=\prod_{l=-1}^m(x^{r_l}-1)^{(-1)^{m-l}}
\in T^{pol}_{HORb}(n,\R).$$
The matrix $S:=S_{(b)}(p)\in T_{HORb}(n,\R)$
is a distinguished matrix of the chain type singularity $f$.
\end{theorem}

\begin{theorem}[{\cite[7.6]{BaH20}}]
\label{t8.7}\cite{BaH20}
Consider a chain type singularity $f$ as in Definition
\ref{t8.5} and the matrix $S_{(b)}(p)$ in Theorem \ref{t8.6}. 
Its spectrum $\Sp(S_{(b)}(p))$ coincides with the normalized
(so symmetric around 0) spectrum $\Sp^{norm}(f)$ 
from the mixed Hodge structure of the singularity $f$. 
\end{theorem}

\begin{remarks}\label{t8.8}
(i) The $b\in\{1,2\}$ in Definition \ref{t8.3} and Theorem
\ref{t8.4} is in the case of an isolated hypersurface singularity
$k+1$ where $k$ is as in \eqref{5.7}. The $b$ here is called
$k$ in \cite{BaH20}.

(ii) In \cite{BaH20} in Remark 7.4 (ii)--(iv) one has to 
replace $m$ by $m+1$, except for $(-1)^{m+1}M$ in (iii).
In \cite{BaH20} in the steps 3 and 4 in the proof of Theorem 7.6,
one has to exchange $m\equiv 0(2)$ and $m\equiv 1(2)$.
\end{remarks}

\begin{examples}\label{t8.9}
(i) As the recipe of Cecotti and Vafa works for matrices
$S\in T_{HORb}(n,\Z)$, one can search for
such matrices $S$ which satisfy the variance inequality 
\eqref{1.4} and also $\tr(S^{-1}S^t)=1$, but which are not
distinguished matrices of singularities. 
By Theorem \ref{t1.5}, then $S+S^t$ is necessarily indefinite.

Uncomfortably, our search was successfull. 
Therefore Theorem \ref{t1.5} does not generalize easily
to the case with $S+S^t$ indefinite. Here are the
first (with respect to $n$) three examples which we found.
(In the following table $\Phi_m$ is the cyclotomic polynomials 
whose roots are the primitve $m$-th roots of unity.)
\begin{eqnarray*}
\begin{array}{l|l|l}
n_1=22 & b_1=2 & p_1=\Phi_1\Phi_2\Phi_3\Phi_4\Phi_8
\Phi_{12}\Phi_{16}=(x^{16}-1)\Phi_3\Phi_{12}\\
n_2=24 & b_2=1 & p_2=\Phi_5\Phi_8\Phi_{15}\Phi_{30}\\
n_3=25 & b_3=2 & p_3=\Phi_1\Phi_6\Phi_{11}\Phi_{28}
\end{array}
\end{eqnarray*}
Then $S_i=S_{(b_i)}(p_i)\in T_{HORb_i}(n_i,\Z)$. 
We discuss the first example in detail in (iii), 
and the other two examples shorter in (iv).

(ii) Instead of $p_i(x)$, one can consider 
$\www{p}_i:=(-1)^{n_i}p_i(-x)$
and suitable $\www{b}_i$. This leads to the matrix
$$\www{S}_i=S_{(\www{b}_i)}(\www{p}_i)=
\diag(1,-1,1,-1,...)S_i\diag(1,-1,1,-1,...),$$ which is
in the $\Br_{n_i}\ltimes\{\pm 1\}^{n_i}$ orbit of $S_i$. 

(iii) The 22 numbers $\beta_1,...,\beta_{22}\in[0;1]$ which are
associated to $p_1$ by Definition \ref{t8.3} (a), are 
$\beta_i=\frac{\www{\beta}_i}{48}$ with $\www{\beta}_i$ as 
follows,
\begin{eqnarray*}
(0,3,4,6,9,12,15,16,18,20,21,24,27,28,30,32,33,36,39,42,44,45).
\end{eqnarray*}
The 22 spectral numbers of $S_1$ are by \eqref{8.2} the 
numbers $\www{\alpha}_i=22\beta_i-(i-1)
=\frac{1}{24}(11\www{\beta}_i-24(i-1))$ for $i=1,...,22$. 
They are as follows,
\begin{eqnarray*}
\frac{1}{24}(0,9,-4,-6,3,12,21,8,6,4,-9,0,9,-4,-6,-8,-21,-12,
-3,6,4,-9).
\end{eqnarray*}
Reordering them to the spectral numbers 
$\alpha_1,...,\alpha_{22}$ with 
$\alpha_1\leq ...\leq \alpha_{22}$
gives
\begin{eqnarray*}
\frac{1}{24}(-21,-12,-9,-9,-8,-6,-6,-4,-4,-3,0,0,3,4,4,6,6,8,9,9,12,21).
\end{eqnarray*}
The characteristic polynomial of $S_1^{-1}S_1^t$ is
\begin{eqnarray*}
p_1(x^{22})=\prod_{j=1}^{22}(x-e^{-2\pi i\alpha_j})
=\Phi_1\Phi_1\Phi_3\Phi_2^2\Phi_4^2\Phi_6^2\Phi_8^2
=(x^8-1)^2\Phi_3\Phi_6^2.
\end{eqnarray*}
Therefore indeed $\tr(S_1^{-1}S_1^t)= 1.$ One also calculates
\begin{eqnarray*}
\Var(\Sp(S_1))&=& \frac{7}{12\cdot 4}
=\frac{\alpha_{22}-\alpha_1}{12}.
\end{eqnarray*}
The distance of $\alpha_1$ to $\alpha_2$
(and of $\alpha_{22}$ to $\alpha_{21}$) is remarkably large.

(iv) In the second example
$$\tr(S^{-1}_2S^t_2)=1\quad\textup{and}\quad 
\Var(\Sp(S_2))=\frac{\alpha_{24}-\alpha_1}{12}.$$
In the third example
$$\tr(S^{-1}_3S^t_3)=1\quad\textup{and}\quad 
\Var(\Sp(S_3))<\frac{\alpha_{25}-\alpha_1}{12}.$$
Also in these two examples the distance of $\alpha_1$ to 
$\alpha_2$ is remarkably large. 
\end{examples}

Part (b) of Conjecture \ref{t1.4} presents an idea how to 
exclude these examples and how to separate the 
distinguished matrices of singularities from the other matrices
in $T_n^{S^1}(\Z)$. It builds on an extension
in \cite{BrH20} of the conjectural variance inequality 
\eqref{1.1} for singularities, see Definition \ref{t8.10} and
Theorem \ref{t8.13} (b).

\begin{definition}[{\cite[ch. 1]{BrH20}}]
\label{t8.10}
Let an abstract spectrum $\Sp=(\alpha_1,...,\alpha_n)\in\R^n$ 
be given with
$$\alpha_1\leq ...\leq \alpha_n,\quad 
\alpha_j+\alpha_{n+1-j}=0.$$

(a) Define the higher moments of this spectrum, 
\begin{eqnarray*}
V_{2k}(\Sp)&:=&\sum_{j=1}^n\alpha_j^{2k} \quad 
2k\textup{\it -th moment for }k\in\Z_{\geq 0},\\
V_2(\Sp)&=&\textup{the non-normalized variance},\\
V_0(\Sp)&=&n.
\end{eqnarray*}
Their generating function is
\begin{eqnarray*}
V(\Sp)&:=& \sum_{k=0}^\infty V_{2k}(\Sp)\cdot\frac{1}{(2k)!}t^{2k} 
= \sum_{j=1}^n\cosh(t\cdot\alpha_j) = \sum_{j=1}^n e^{t\cdot \alpha_j}.
\end{eqnarray*}

(b) For any $\nu\in\R$ and $k\in\Z_{\geq 0}$ 
define the {\it $2k$-th Bernoulli moment}
$$\Gamma^{Ber}_{2k}(\Sp,\nu)\in\sum_{l=0}^k\Q[\nu]_{\deg\leq k-l}
\cdot V_{2l}(\Sp)$$ by the generating function $\Gamma^{Ber}$
for all of them, 
\begin{eqnarray*}
\Gamma^{Ber}(\Sp,\nu)=\sum_{k=0}^\infty 
\Gamma_{2k}^{Ber}(\Sp,\nu)\frac{1}{(2k)!}t^{2k}
= V(\Sp) \cdot \exp\left( \nu\cdot \log\frac{t/2}{\sinh (t/2)}\right).&&
\end{eqnarray*}
\end{definition}

\begin{examples}[{\cite[ch. 1]{BrH20}}]\label{t8.11}
The first four Bernoulli moments are 
\begin{eqnarray*}
\Gamma_0^{Ber}(\Sp,\nu)&=&V_0=n,\\
\Gamma_2^{Ber}(\Sp,\nu)&=&V_2-V_0\cdot(\frac{1}{12}\nu),\\
\Gamma_4^{Ber}(\Sp,\nu)&=&V_4 - V_2\cdot(\frac{1}{2}\nu) + 
V_0\cdot(\frac{1}{120}\nu +\frac{1}{48}\nu^2), \\
\Gamma_6^{Ber}(\Sp,\nu)&=&V_6 - V_4\cdot(\frac{5}{4}\nu) 
+ V_2\cdot(\frac{1}{8}\nu + \frac{5}{16}\nu^2) \nonumber\\
&&- V_0\cdot(\frac{1}{252}\nu + \frac{1}{96}\nu^2 + \frac{5}{576}\nu^3).
\end{eqnarray*}
\end{examples}

The following theorem collects some interesting facts around
the Bernoulli moments. Part (a) implies with little work
part (b). 

\begin{theorem}\label{t8.12}
(a) (Classical, e.g. \cite[Theorem 3.1]{BrH20}) Recall 
$$\log\frac{t/2}{\sinh (t/2)}=\sum_{k=1}^\infty 
\frac{-1}{2k}B_{2k}\frac{1}{(2k)!}t^{2k},$$
with the {\it Bernoulli numbers}
$$B_{2k}=(-1)^{k-1}\frac{2(2k)!}{(2\pi)^{2k}}\zeta(2k)
\in (-1)^{k-1}\Q_{>0}\quad\textup{for}\quad k\geq 1$$
(recall also that $\zeta(2k)\to 1$ fast for $k\to\infty$).

(b) \cite[Lemma 2.4 (b)]{BrH20}
If a number $\nu\in\R$ satisfies 
$$\forall\ k\quad (-1)^k\Gamma^{Ber}_{2k} (\Sp,\nu)\geq 0$$
then also any number $\www \nu$ with $\www \nu\geq \nu$ 
satisfies this.

(c) \cite[Theorem 3.2 (c)]{BrH20}
For any $\nu\in\R-\Z_{\leq 0}$, the sequence of numbers
$$(-1)^k\Gamma^{Ber}_{2k}(\Sp,\nu)\cdot
\frac{(2\pi)^{2k}\cdot\Gamma(\nu)}
{2\cdot(2k)!\cdot (2k)^{\nu-1}}$$
is convergent for $k\to\infty$ with limit  
$$\sum_{j=1}^n\cos(-2\pi\alpha_j)
=\sum_{j=1}^ne^{-2\pi\alpha_j}.$$
\end{theorem}

Part (a) of the following conjecture is the variance
conjecture which Hertling proposed 2000.
Part (b) is the generalization in \cite{BrH20}.

\begin{conjecture}\label{t8.13}
Let $f(z_0,...,z_m)$ be an isolated hypersurface singularity
with Milnor number $n$ and normalized spectrum 
$\Sp^{norm}(f)=(\alpha_1,...,\alpha_n)$ from its mixed Hodge
structure. Write $\nu(f):=\alpha_n-\alpha_1$. 

(a) \cite{He02} Variance conjecture:
\begin{eqnarray*}
&&\Var(\Sp^{norm}(f)) \leq \frac{\nu(f)}{12}.\\
\textup{Equivalent:} 
&&0\leq (-1)\Gamma_2^{Ber}(\Sp^{norm}(f),\nu(f)).
\end{eqnarray*}

(b) \cite[Conjecture 1.2]{BrH20} Generalization for Bernoulli moments:
\begin{eqnarray*}
\forall\ k\in\N\quad 
0\leq (-1)^k\Gamma_{2k}^{Ber}(\Sp^{norm}(f),\nu(f)).
\end{eqnarray*}
\end{conjecture}

This conjecture is open in most cases. But it is true for
quasihomogeneous singularities. In these cases even more is 
known, see Theorem \ref{t8.14}.

\begin{theorem}\label{t8.14}
Consider a quasihomogeneous singularity 
$f(z_0,...,z_m)$ with weights 
$(w_0,...,w_m)\in(0;\frac{1}{2}]\cap\Q$.
Write $\nu(f):=\alpha_n-\alpha_1\stackrel{!}{=}
m+1-2\sum_{i=0}^mw_i$. 

(a) \cite[Theorem 5.4]{BrH20}
Conjecture \ref{t8.13} (b) is true.

(b) \cite[Theorem 5.4]{BrH20}
\begin{eqnarray}\label{8.3}
\Gamma^{Ber}_4(\Sp^{norm}(f),\nu(f))
=\frac{1}{30}\cdot n\cdot\sum_{i=0}^m(\frac{1}{2}-w_i)w_i(1-w_i).
\end{eqnarray}

(c) (New)
\begin{eqnarray}\label{8.4}
0\leq \Gamma^{Ber}_4(\Sp^{norm}(f),\nu(f))
\leq \frac{1}{240}\cdot n\cdot \nu(f),\\
n\cdot\frac{\nu(f)}{48}\cdot(\nu(f)-\frac{2}{5})
\leq V_4(\Sp^{norm}(f)) \leq 
n\cdot\frac{\nu(f)}{48}\cdot(\nu(f)-\frac{1}{5}).\label{8.5}
\end{eqnarray}
\end{theorem}

{\bf Proof: of part (c):}
Because $w_i\in(0,\frac{1}{2}]$, $w_i(1-w_i)\leq\frac{1}{4}$.
Part (b) gives
\begin{eqnarray*}
\Gamma^{Ber}_4(\Sp^{norm}(f),\nu(f))
\leq \frac{1}{30}\cdot n\cdot\sum_{i=0}^m(\frac{1}{2}-w_j)
\cdot \frac{1}{4}
=\frac{1}{240}\cdot n\cdot \nu(f),
\end{eqnarray*}
so \eqref{8.4}. Together with the formula for 
$\Gamma^{Ber}_4(V,\nu)$ in Example \ref{t8.11}
and the fact $V_2(\Sp^{norm}(f))=n\cdot \frac{\nu(f)}{12}$
for a quasihomogeneous singularity, this implies \eqref{8.5}.
\hfill$\Box$

\begin{examples}\label{t8.15}
We reconsider the Examples \ref{t8.9}.
Write $\nu:=\alpha_n-\alpha_1$. 
In all three examples the large distance between $\alpha_1$ 
and $\alpha_2$ leads to a large value 
$\Gamma^{Ber}_4(\Sp(S_1),\nu)$, which is in fact larger
than $\frac{1}{240}\cdot n\cdot \nu$. 

In the first example $S_1$
$$\Gamma^{Ber}_4(\Sp(S_1),\nu)=0,3375..,\quad
\frac{1}{240}\cdot n\cdot \nu=0,1604.. .$$
So here the right inequality in \eqref{8.4} does not hold.

In the first and second example $\Gamma^{Ber}_2(\Sp(S_i),\nu)=0$,
so the comparison with quasihomogeneous singularities is
reasonable. In the third example
$\Gamma^{Ber}_2(\Sp(S_i),\nu)<0$, and the comparison only
with quasihomogeneous singularities is not reasonable.
\end{examples}

\begin{examples}\label{t8.16}
(i) The examples in \ref{t8.9} suggest that we need upper
estimates for 
$(-1)^k\Gamma^{Ber}_{2k}(\Sp(S),\alpha_n-\alpha_1)$
in order to separate the distinguished matrices of singularities
from the other matrices. 
Part (c) of Theorem \ref{t8.12} might lead 
to the hope that the limit might also be an upper bound, so 
$$1\stackrel{?}{\geq}|\Gamma^{Ber}_{2k}(V,\nu)|\cdot
\frac{(2\pi)^{2k}\cdot\Gamma(\nu)}
{2\cdot(2k)!\cdot (2k)^{\nu-1}}
=|\Gamma^{Ber}_{2k}(V,\nu)|\cdot
\frac{\zeta(2k)\cdot\Gamma(\nu)}{|B_{2k}|\cdot (2k)^{\nu-1}}.$$
But this hope is too optimistic, as the family of examples
$A_2^{m+1}$ for $m\in\Z_{\geq 0}$ in part (ii) shows
in the case $k=2$.

(ii) The singularity $A_2^{m+1}$ for $m\in\Z_{\geq 0}$ is the 
quasihomogeneous polynomial $\sum_{i=0}^m z_i^3$ with
weights $(w_0,...,w_m)=(\frac{1}{3},...,\frac{1}{3})$,
$n=2^{m+1}$ and 
\begin{eqnarray*}
\nu(f)&=&\alpha_n-\alpha_1=m+1-2\sum_{i=0}^mw_i=\frac{m+1}{3},\\
\Gamma^{Ber}_4(\Sp^{norm}(f),\nu(f))
&=&\frac{1}{30}\cdot n\cdot 
\sum_{i=0}^m(\frac{1}{2}-w_i)w_i(1-w_i)\\
&=&\frac{1}{30}\frac{2^{m+1}(m+1)}{27}
=B_4\frac{2^{m+1}(m+1)}{27}.
\end{eqnarray*}
\end{examples}

We do not have good candidates for upper bounds for the higher Bernoulli moments. The following conjecture is thus not
as precise as we wish.

\begin{conjecture}\label{t8.17}
A family of functions $r(k,n,.):\R_{>0}\to \R_{>0}$ for
$k,n\in\N$ exists which satisfies the following properties.

(a) For any singularity $f$ with Milnor number $n$ 
\begin{eqnarray*}
(-1)^k\Gamma^{Ber}_{2k}(\Sp^{norm}(f),\alpha_n-\alpha_1)
\leq r(k,n,\alpha_n-\alpha_1).
\end{eqnarray*}

(b) For any matrix $S\in T_n^{S^1}(\Z)$
hopefully the Cecotti-Vafa recipe works and gives a
spectrum $\Sp(S)=(\alpha_1,...,\alpha_n)$ with 
$\alpha_1\leq ...\leq\alpha_n$. Then $S$ is a distinguished
matrix of a singularity if and only if its 
Coxeter-Dynkin diagram ${\rm CDD}(S)$ is connected and for all
$k\in\N$
\begin{eqnarray*}
0\leq (-1)^k\Gamma^{Ber}_{2k}(\Sp(S),\alpha_n-\alpha_1)
\leq r(k,n,\alpha_n-\alpha_1).
\end{eqnarray*}

(c) For any $n\in\N$ and any $\nu\in\R_{>0}$ the sequence 
$$\Bigl(r(k,n,\nu)\cdot \frac{\Gamma(\nu)}{|B_{2k}|\cdot (2k)^{\nu-1}}\Bigr)_{k\in\N}$$
is convergent for $k\to\infty$ with limit $1$. 
\end{conjecture}

\end{document}